\documentclass[11pt]{article}
 \usepackage{benstyle}
 \usepackage{cite}
 
\usepackage[font=small]{caption}

\title{Compressive Hermite interpolation: sparse, high-dimensional approximation from gradient-augmented measurements}
\author{Ben Adcock \\ Department of Mathematics \\ Simon Fraser University \\ Canada \and Yi Sui \\ Department of Mathematics \\ Simon Fraser University \\ Canada}

\begin{document}

\maketitle

\begin{abstract}
We consider the sparse polynomial approximation of a multivariate function on a tensor product domain from samples of both the function and its gradient.  When only function samples are prescribed, weighted $\ell^1$ minimization has recently been shown to be an effective procedure for computing such approximations.  We extend this work to the gradient-augmented case.  Our main results show that for the same asymptotic sample complexity, gradient-augmented measurements achieve an approximation error bound in a stronger Sobolev norm, as opposed to the $L^2$-norm in the unaugmented case.  For Chebyshev and Legendre polynomial approximations, this sample complexity estimate is algebraic in the sparsity $s$ and at most logarithmic in the dimension $d$, thus mitigating the curse of dimensionality to a substantial extent.
We also present several experiments numerically illustrating the benefits of gradient information over an equivalent number of function samples only.
\end{abstract}

\section{Introduction}\label{s:introduction}

The concern of this paper is the approximation of a smooth, high-dimensional function $f : (-1,1)^d \rightarrow \bbR$ using multivariate polynomials.  Recent years have seen an increasing focus on this problem, due to its applications in Uncertainty Quantification (UQ), where the function $f$ is typically a solution of a parametric PDE.

In a typical setup, which we shall also consider in this paper, $f$ is expressed as an expansion in an orthogonal basis of polynomials according to some tensor-product probability measure, often referred to as a \textit{Polynomial Chaos Expansion}.  Samples are drawn randomly and independently according to this measure, and then the objective is to compute the $s$ expansion coefficients in some finite index set accurately from the corresponding measurements of $f$.  Least-squares fitting has often been used to effect this approximation \cite{ChkifaDLS,CohenLSA,MiglioratiJAT,MiglioratiFCM,MiglioratiSIAM,Cohenpr,Hadigolpr,Seshadri}.  However, in last several years there has been an increasing focus on the use of sparse regularization procedures for this task, based on the principles of compressed sensing \cite{RauhutWard,RauhutWardWeighted,BAptwise,AdcockCSFunInterp,Chkifalower,Pengpce,Tangthesis}.  The efficacy of such procedures has recently been theoretically established.  Specifically, it has been shown that suitable weighted $\ell^1$ minimization procedures achieve quasi-optimal error decay rates for approximations in so-called lower sets.  The corresponding sample complexities are algebraic in the number of coefficients $s$ sought and only (poly)logarithmic in the dimension $d$ \cite{Chkifalower, AdcockCSFunInterp, BASBCWMatheon}.  Hence the curse of dimensionality is significantly ameliorated.

In this paper, we consider the extension and analysis of sparse regularization procedures for the modified problem where both $f$ and its gradient $\nabla f$ are measured at the sample points.  This can be viewed as a multivariate extension of the classical Hermite interpolation problem in numerical analysis.  Yet this problem is increasingly encountered in UQ applications (see, for example, \cite{Penggradient} and references therein), where gradient measurements can be computed relatively inexpensively via, for example, adjoint sensitivity analysis \cite{Komkovsen}.  As is typical, our objective is to use this additional information to enhance the accuracy of the computed approximation to $f$.

\subsection{Contributions}

In \cite{AdcockCSFunInterp} it was shown that a certain weighted $\ell^1$ minimization procedure produces a quasi-optimal best $s$-term approximation in lower sets using a number of measurements that polynomial in $s$ and logarithmic in $d$.  Up to the logarithmic factors, these sample complexity bounds are identical to the best known estimates for oracle least-squares estimators based on \textit{a priori} knowledge of the support set.  We review these results in more detail in \S \ref{s:prelims}.

The primary contribution of this paper is to extend this work to the case of gradient-augmented measurements.  Our main result shows that recovery from gradient-enhanced samples can be achieved under the same sufficient condition on the sample complexity, up to minor variations in the logarithmic factor.  However, the approximation error -- which in \cite{AdcockCSFunInterp} is evaluated in an $L^2$-norm -- is for the gradient-enhanced problem evaluated in a stronger $H^1$-type norm.  In other words, by sampling both $f$ and $\nabla f$ one guarantees an error bound in a stronger norm, under the same asymptotic measurement condition.  

The analysis in \cite{AdcockCSFunInterp} is considered for Legendre and Chebyshev polynomial approximations.  Our work extends this to Jacobi polynomials, and furthermore, to any orthonormal basis of functions (not necessarily polynomials) arising as eigenfunctions of a singular Sturm--Liouville problem.  We also briefly discuss the case of regular Sturm--Liouville problems; in particular, the Fourier basis (i.e.\ multivariate trigonometric polynomial approximation).

Our analysis provides a theoretical insight into the advantage conveyed by gradient information.  We also present a series of numerical results to compare gradient-augmented measurements with function samples only when the error is measured in the same norm (specifically, the $L^\infty$-norm).  Using the cost model that $\nabla f$ can be computed in roughly the same time as $f$ (which is realistic in some applications), these results show that the former can achieve a smaller error for a comparable computational cost; another advantage of using gradient information.

Finally, we discuss several variations on the setup.  For instance, the problem where $\nabla f$ is only evaluated at a fraction of the sample points,
and when $f$ and $\nabla f$ are sampled at different points.

\subsection{Previous work}

Sparse Legendre approximations from gradient-enhanced measurements was first investigated empirically in \cite{Tangthesis}.  In \cite{Penggradient}, the authors made a first theoretical analysis using compressed sensing techniques with Hermite polynomials.  Specifically, for unweighted $\ell^1$-minimization it was shown that gradient-enhancement leads to a better null space property and a smaller coherence, both of which are sufficient conditions for recovery.  Related analysis of $\ell^1$-minimization has been given in \cite{GuoEtAlGradient} and \cite{XuZhouGradient}, with the latter considering the case of Fourier expansions.  We note in passing, however, that unweighted $\ell^1$-minimization does not overcome the curse of dimensionality in high-dimensional approximation.  The best known sample complexity estimates all involve factors that are exponentially-large in the dimension $d$ or the degree of the polynomial space, and therefore significantly worse than those of oracle estimators.  Conversely, as mentioned, weighted $\ell^1$ minimization has sample complexities that agree with those of oracle estimators, up to logarithmic factors.

In this paper we use gradient measurements to effect a Hermite polynomial interpolant, i.e.\ a polynomial which interpolates both $f$ and $\nabla f$ at the nodes\footnote{This is not to be confused with expansions in Hermite polynomials, which we do not address in this paper.  See \cite{Penggradient} for some work in this direction.}.  We note in passing that gradient information can  also be used in other ways, for instance as part of dimensionality reduction techniques \cite{ConstantineBook}.  We make no attempt to compare these procedures in this paper, as they address quite fundamentally different function classes (e.g.\ ridge functions).  Finally, for applications of gradient-enhanced measurements to UQ problems, we refer to \cite{Penggradient,LockwoodUQ,AlekseevUQ,Li}.

\subsection{Outline}

The outline of this paper is as follows.  In \S \ref{s:prelims} we introduce the polynomial approximation problem, and define a number of key concepts, including lower sets.  The gradient-augmented problem is formulated in \S \ref{s:gradaug}, along with the relevant weighted Sobolev spaces.  With this in hand, the main results of the paper are given in \S \ref{s:mainres}.  Next in \S \ref{s:numexp} we present numerical experiments, and finally, in \S \ref{s:proofs} we give the proofs of the main results.

\section{Background}\label{s:prelims}

In this section, we review the main aspects of polynomial approximation of high-dimensional functions without gradient enhancement using weighted $\ell^1$ minimization.  We follow the setup of \cite{AdcockCSFunInterp}.

\subsection{Notation}

We first require some notation.  Throughout $y \in (-1,1)$ and $\bm{y} = (y_1,\ldots,y_d) \in D$ denote the one- and $d$-dimensional variables respectively, where $D = (-1,1)^d$ is the $d$-dimensional domain.  The function to recover is denoted by $f : D \rightarrow \bbC$.  We write $\nu(y)$ for a probability density function on $(-1,1)$ and $\nu(\bm{y})  = \prod^{d}_{i=1} \nu(y_i)$ for the corresponding tensor-product probability density function on $D$.  The spaces of square-integrable functions with respect to $\nu$ are denoted by $L^2_{\nu}(-1,1)$ and $L^2_{\nu}(D)$ respectively.  We write $\nm{\cdot}_{L^2(D)}$ and $\ip{\cdot}{\cdot}_{L^2(D)}$ for the corresponding norm and inner product.

We consider approximations in orthonormal bases on these spaces, which are typically (but not necessarily) of polynomial type.  We write $\{ \phi_n \}^{\infty}_{n=0}$ for a one-dimensional orthonormal basis of $L^2_{\nu}(-1,1)$ and $\{ \phi_{\bm{n}} \}_{\bm{n} \in \bbN^d_0}$ for the corresponding tensor-product orthonormal basis of $L^2_{\nu}(D)$, i.e.\
\bes{
\phi_{\bm{n}}(\bm{y}) = \prod^{d}_{i=1} \phi_{n_i}(y_i) ,\qquad \bm{n} = (n_1,\ldots,n_d) \in \bbN^d_0.
}
Here and throughout, $\bm{n} = (n_1,\ldots,n_d)$ is a multi-index in $\bbN^d_0$.  We write $\Lambda$ for the finite set of multi-indices from which the approximation to $f$ is sought, and $N = | \Lambda |$ for its cardinality.  We also use $\Delta$ to denote a finite multi-index set, typically of size $|\Delta| = s$, corresponding to the coefficients of $f$ that give the best or quasi-best $s$-term approximation, or more frequently, the best or quasi-best $s$-term approximation in lower sets.

The norm $\nm{\cdot}_{2}$ and inner product $\ip{\cdot}{\cdot}$ denote the $\ell^2$-norm and inner product on either $\bbC^N$ or $\ell^2(\bbN^d_0)$.  Given an infinite vector of positive weights $\bm{w} = (\bm{w}_{\bm{n}})_{\bm{n} \in \bbN^d_0}$ we write $\nm{\cdot}_{1,\bm{w}}$ for the norm on the weighted space $\ell^1_{\bm{w}}(\bbN^d_0)$
\bes{
\nm{\bm{x}}_{1,\bm{w}} = \sum_{\bm{n} \in \bbN^d_0} w_{\bm{n}} | x_{\bm{n}} |,
}
and likewise for finite vectors of  positive weights in $\bbR^N$.  

We consider approximating $f$ from samples taken at $m$ points denoted by $\bm{y}_1,\ldots,\bm{y}_m$.  As discussed, these will be chosen randomly according to some measure.  To this end, we let  $\mu(y)$ be a probability density function on $(-1,1)$ and $\mu(\bm y) = \prod^{d}_{i=1} \mu(y_i)$ be the corresponding tensor-product probability measure.  Typically, but not always, we have $\mu = \nu$.

Finally, for $k=1,\ldots,d$ we let $\partial_{k}$ be the partial derivative operator with respect to $y_{k}$, i.e.\ $\partial / \partial y_{k}$.  For convenience, we also write $\partial_{0}$ to mean the identity operator, i.e. $\partial_0 f = f$.

\subsection{Weighted $\ell^1$ minimization}\label{ss:weightedl1}

Let $\{ \phi_{\bm{n}} \}_{\bm{n} \in \bbN^d_0}$ be a tensor-product orthonormal basis of $L^2_{\nu}(D)$, where $\nu$ is a tensor-product probability density function.  Then we can write any $f \in L^2_{\nu}(D)$ as
\bes{
f = \sum_{\bm{n} \in \bbN^d_0} x_{\bm{n}} \phi_{\bm{n}} ,\qquad \bm{x}_n = \ip{f}{\phi_{\bm{n}}}_{L^2_{\nu}(D)}.
}
In order to approximate $f$ we first truncate this expansion using the multi-index set $\Lambda$.  Write
\be{
\label{eLambda}
f = f_{\Lambda} + e_{\Lambda} = \sum_{\bm{n} \in \Lambda} x_{\bm{n}} \phi_{\bm{n}} + \sum_{\bm{n} \notin \Lambda} x_{\bm{n}} \phi_{\bm{n}},
}
and let $\bm{x} = \left ( \bm{x}_{\bm{n}} \right )_{\bm{n} \in \bbN^d_0} \in \ell^2(\bbN^d_0)$ be the infinite vector of coefficients.  For reasons discussed in \S \ref{ss:lowersets}, given $s \geq 1$ we choose $\Lambda$ as the \textit{hyperbolic cross} index set of degree $s$:
\be{
\label{HCset}
\Lambda = \Lambda^{\mathrm{HC}}_{s} = \left \{ \bm{n} \in \bbN^d_0 : \prod^{d}_{k=1} (n_k+1) \leq s  +1  \right \}.
}
Let
\be{
\label{ordering}
\bm{n}_1,\ldots,\bm{n}_N,
}
be an ordering of the multi-indices in $\Lambda$.  Then we write $\bm{x}_{\Lambda} = (x_{\bm{n}})_{\bm{n} \in \Lambda} = (x_{\bm{n}_i})^{N}_{i=1} \in \bbC^{N}$ for the corresponding finite vector of coefficients.  Here and through the paper we shall index over the multi-index set $\Lambda$ or the index set $\{1,\ldots,N\}$ (using \R{ordering}) interchangeably.  The meaning will be clear from the context.

Let $\mu$ be another tensor-product probability density function on $D$.  For technical reasons, we assume throughout that
\be{
\label{mucond}
\sup_{\bm{y} \in D} \sqrt{\nu(\bm{y})/\mu(\bm{y})} |\phi_{\bm{n}}(\bm{y})|  < \infty,\qquad \forall \bm{n} \in \bbN^d_0.
} 
Note that this condition holds in particular when $\mu = \nu$ and the $\phi_{\bm{n}}$ are polynomials.  Let $\bm{y}_{1},\ldots,\bm{y}_{m} \in D$ be sample points, drawn independently and randomly according to $\mu$.  If
\be{
\label{Apoly}
A = \frac{1}{\sqrt{m}} \left ( \phi_{\bm{n}_j}(\bm{y}_i) \right )^{m,N}_{i,j=1} \in \bbC^{m \times N},
}
is the resulting measurement matrix, then we have the linear system of equations
\be{
\label{unaug}
\bm{f} = A \bm{x}_{\Lambda} + \bm{e},
\qquad
\mbox{where}\ \bm{f} = \frac{1}{\sqrt{m}} \left ( f(\bm{y}_i) \right )^{m}_{i=1},\ \bm{e} = \frac{1}{\sqrt{m}} \left ( e_{\Lambda}(\bm{y}_i) \right )^{m}_{i=1}.
}
Suppose now that $\bm{e}$ satisfies 
\be{
\label{eeta}
\nm{\bm{e}}_2 \leq \eta,
}
for some known $\eta \geq 0$ (see Remark \ref{r:eta} below).  Then, given weights $\bm{w} = \left ( w_{\bm{n}} \right )_{\bm{n} \in \bbN^d_0}$ with $w_{\bm{n}} \geq 1$, $\forall \bm{n}$, we consider the weighted $\ell^1$ minimization problem
\be{
\label{weightedl1}
\min_{\bm{z} \in \bbC^N} \nm{\bm{z}}_{1,\bm{w}}\ \mbox{subject to $\nm{A \bm{z} - \bm{y}}_2 \leq \eta$}.
}
If $\bm{\hat{x}} \in \bbC^N$ is a minimizer of this problem, then the resulting approximation to $f$ is given by
\be{
\label{weightedl1fn}
\hat{f} = \sum_{\bm{n} \in \Lambda} \hat{x}_{\bm{n}} \phi_{\bm{n}}.
}

\rem{
\label{r:eta}
In practice, a bound such as \R{eeta} may not be available, since $\bm{e}$ depends on the unknown function $f$.  Recovery guarantees for sparse regularization under \textit{unknown errors} have been considered in \cite{BASBCSmodel} and \cite{ABBCorrecting}.  In particular, \cite{ABBCorrecting} shows that a weighted version of the square-root LASSO optimization problem can successfully avoid the \textit{a priori} bound \R{eeta}.  For simplicity, we shall not consider this in this paper, although we expect a similar result to hold in this case as well.
}

\subsection{Lower sets}\label{ss:lowersets}

Standard compressed sensing \cite{FoucartRauhut,CandesWakin} concerns the recovery of a vector of coefficients $\bm{x} \in \bbC^N$ that is approximately \textit{sparse}; that is, well-approximated by its best $s$-term approximation.  Its signature results show recovery of $\bm{x}$ up to its best $s$-term approximation error from a suitable measurement matrix $A$ with a number of measurements $m$ that is linear in $s$ and logarithmic in $N$.  This recovery can be effected using constrained $\ell^1$ minimization, for example.  

Unfortunately, the measurement matrices \R{Apoly} arising in multivariate polynomial approximation do not give optimal guarantees for the recovery of approximately sparse polynomial coefficients via $\ell^1$ minimization.  The best known estimates involve exponentially-large factors in either $d$ or the polynomial degree $s$ \cite{BASBCWMatheon,HamptonDoostanCSPCE,YanGuoXui_l1UQ}, and therefore suffer from the curse of dimensionality.  

However, recent work \cite{AdcockCSFunInterp,Chkifalower} has shown that such estimates are not sharp, and that polynomial coefficients can be recovered with much lower (and nearly-optimal) sample complexities.  The key is to exploit the additional structure that polynomial coefficients of smooth, high-dimensional functions possess; specifically, \textit{lower set} structure:

\defn{
A set $\Delta \subseteq \mathbb{N}^d_0$ is lower if whenever $\bm{n} = (n_1, \dots, n_d) \in \Delta$ and 
$\bm{n}' = (n_1', \dots, n'_d) \in \mathbb{N}_0^d$ satisfies $n_k' \leq n_k$, $k =1, \dots, d$, then $\bm{n}' \in \Delta$. 
}

Lower sets (also known as monotone or downward closed sets) have been studied extensively in the context of multivariate polynomial approximation \cite{Cohenlower,ChkifaDLS, ChkifaHD}.  In particular, for functions arising as solutions of a broad class of parametric PDEs it has been shown that there exist sequences of lower sets of cardinality $s$ which achieve the same approximation error bounds as those of the best $s$-term approximation \cite{Chkifapde}. 

In tandem with these results, a series of works \cite{AdcockCSFunInterp, BASBCWMatheon, Chkifalower} have shown that quasi-best $s$-term approximations in lower sets can be obtained by solving the weighted $\ell^1$ minimization problem \R{weightedl1} with a suitable choice of weights.  Since the union of all lower sets of size $s$ is precisely the hyperbolic cross index set
\be{
\label{unionlower_HC}
\bigcup \left \{ \Delta : | \Delta |\leq s,\ \mbox{\small $\Delta$ lower} \right \} = \Lambda^{\mathrm{HC}}_{s},
}
the approach developed in \cite{BASBCWMatheon} computes an approximation $\hat{f}$ to $f$ via \R{weightedl1}, using this choice of truncated index set.  Due to the additional structure imposed by lower sets, and the promotion of this structure via the weights, the sample complexity estimates transpire to be at most logarithmic in the dimension $d$, and polynomial in $s$ for large classes of polynomial bases.  Moreover, these estimates agree (up to possible log factors) with the best known estimates for oracle estimators based on lower sets.  We refer to \S \ref{ss:discussion} for the specific estimates.

The main results of this paper extend this analysis to the gradient-augmented setting.  Correspondingly, we derive conditions on $m$ under which the approximation error $f - \tilde{f}$ (measured in a suitable Sobolev norm) can be estimated in terms of the $\ell^1_{\bm{w}}$-norm error of the best lower $s$-term approximation of $\bm{x}$:
\be{
\label{sigma_sL}
\sigma_{s,L}(\bm{x})_{1,\bm{w}} = \inf \left \{ \nm{\bm{x} - \bm{z} }_{1,\bm{w}} : \ \bm{z} \in \ell^1_{\bm{w}}(\bbN^d_0),\ | \supp(\bm{z}) | \leq s,\ \mbox{$ \supp(\bm{z})$ lower} \right \}.
}
Here $\supp(\bm{z}) = \{ i : z_{i} \neq 0 \}$ is the set of indices where $\bm{z}$ is nonzero.  As mentioned above for functions arising as solutions of parametric PDEs $\sigma_{s,L}(\bm{x})_{1,\bm{w}}$ is a reasonable surrogate for the true best $s$-term approximation
\bes{
\sigma_{s}(\bm{x})_{1,\bm{w}} = \inf \left \{ \nm{\bm{x} - \bm{z} }_{1,\bm{w}} : \ \bm{z} \in \ell^1_{\bm{w}}(\bbN^d_0),\ | \supp(\bm{z}) | \leq s \right \}.
}

\section{Recovery from gradient-augmented measurements}\label{s:gradaug}

Having reviewed weighted $\ell^1$ minimization for polynomial approximation, we now extend it to the gradient-augmented setting.  Our main tool to do so will be Sturm--Liouville theory, described next.

\subsection{Sturm--Liouville eigenfunctions}\label{s:SL}

Recall that a Sturm--Liouville problem is an eigenvalue problem of the form
\be{
\label{SLprob}
-(\chi u')' + \zeta u = \lambda \nu u,
}
where $\chi$ is continuously differentiable and positive in $(-1,1)$ and continuous in $[-1,1]$, $\zeta$ is continuous in $[-1,1]$ and $\nu$ is continuous and nonnegative in $(-1,1)$ and integrable.   The problem is singular if $\chi(\pm 1) = 0$.    Such a problem has a countable set of eigenvalues $0 \leq \lambda_0 < \lambda_1 < \ldots$ and eigenfunctions $\{ \phi_n \}_{n \in \bbN_0}$, with the latter constituting an orthogonal basis of $L^2_{\nu}(-1,1)$.  

Of relevance to this paper, the classical orthogonal polynomials are all singular Sturm--Liouville eigenfuntions:

\pbk
\textit{Legendre polynomials.} These are Sturm--Liouville eigenfunctions corresponding to
\bes{
\chi(y) = \frac12 (1-y^2),\qquad \nu(y) = \frac12,\qquad \zeta(y) = 0.
}
The corresponding eigenvalues are $\lambda_{n} = n(n+1)$.
Note that it is customary to write $\chi(y) = 1-y^2$ and $\nu(y) = 1$ here.  We have normalized by $1/2$ so that $\nu$ is a probability density function.

\pbk
\textit{Chebyshev polynomials.} These are Sturm--Liouville eigenfunctions corresponding to
\bes{
\chi(y) = \frac{\sqrt{1-y^2}}{\pi},\qquad \nu(y) = \frac{1}{\pi \sqrt{1-y^2}},\qquad \zeta(y) = 0.
}
The corresponding eigenvalues are $\lambda_{n} = n^2$.

\pbk
\textit{Jacobi polynomials.} These are Sturm--Liouville eigenfunctions corresponding to
\be{
\label{Jacobi_pw}
\chi(y) = \frac{1}{c^{(\alpha,\beta)}}(1-y)^{\alpha+1}(1+y)^{\beta+1},\qquad 
\nu(y) = \frac{(1-y)^{\alpha}(1+y)^{\beta}}{c^{(\alpha,\beta)}},\qquad \zeta(y) = 0.
}
where $\alpha,\beta > -1$ and $c^{(\alpha,\beta)} = \int^{1}_{-1} (1-y)^{\alpha}(1+y)^{\beta} \D y$.  The corresponding eigenvalues are
\be{
\label{lambda_Jacobi}
\lambda^{(\alpha,\beta)}_{n} = n ( n+\alpha+\beta + 1).
}
Note that Jacobi polynomials include both Legendre and Chebyshev polynomials as the special cases $\alpha = \beta = 0$ and $\alpha = \beta = -1/2$ respectively.

\pbk
Throughout the paper, we assume that the orthonormal basis $\{ \phi_n \}^{\infty}_{n=0}$ introduced in \S \ref{s:prelims} arises as the eigenfunctions of a singular Sturm--Liouville problem \R{SLprob}. For convenience we also assume that
\be{
\label{zetazero}
\zeta(y) = 0.
}
This is not strictly necessary for what follows.  However, it holds for all cases relevant to this paper; specifically, the classical orthogonal polynomials discussed above.

\subsection{Sobolev orthogonality}
The main advantage of this setup for the gradient-augmented problem is that the derivatives of Sturm--Liouville eigenfunctions are also orthogonal in a particular weighted $L^2$ space.  We now formalize this notion.  Note that this space does not usually coincide with the original weighted space $L^2_{\nu}(-1,1)$.  Two exceptions are the Fourier and Hermite bases, studied in \cite{XuZhouGradient} and \cite{Penggradient} respectively.  The change of weight that occurs in the general case requires some additional effort when deriving the gradient-enhanced system.  See \S \ref{ss:gradenhance}. 

Consider equation \R{SLprob}.  Multiplying by $\overline{\phi_m}$, integrating by parts and using the fact that $\chi(\pm 1) = 0$ since the problem is assumed to be singular, we get

\eas{
 \int^{1}_{-1} -(\chi(y) \phi'_n(y))' \overline{\phi_m(y)} \D y &= - \chi(y)\phi'_n(y)\overline{\phi_m(y)} \Big|_{-1}^1 + 
  \int^{1}_{-1} \chi(y) \phi'_n(y) \overline{\phi'_m(y)} \D y 
  \\
& =  
  \lambda_n\int^{1}_{-1}\nu(y)\phi_n(y)\overline{\phi_m(y)} \D y ,
}
Hence the derivatives $\phi'_n$ are orthogonal in $L^2_\chi(-1,1)$:
\be{
\label{derivIP}
 \int^{1}_{-1} \chi(y) \phi'_n(y) \overline{\phi'_m(y)} \D y = \lambda_{n} \delta_{n,m},\qquad n,m = 0,1,\ldots.
}

Now define the weighted Sobolev space
\bes{
\tilde{H}^1(-1,1) = \left \{ f \in L^2_{\nu}(-1,1) : f' \in L^2_\chi(-1,1) \right \},
}
with norm and inner product
\bes{
\nm{f}^2_{\tilde{H}^1(-1,1)} = \nm{f}^2_{L^2_\nu(-1,1)} + \nm{f'}^2_{L^2_\chi(-1,1)},
\qquad
\ip{f}{g}_{\tilde{H}^1(-1,1)} = \ip{f}{g}_{L^2_\nu(-1,1)} + \ip{f'}{g'}_{L^2_\chi(-1,1)}.
}
It follows from \R{derivIP} that the functions
\bes{
\psi_{n}(y) = \frac{1}{\sqrt{1+\lambda_n}} \phi_n(y),\qquad n= 0 ,1,2\ldots ,
}
are an orthonormal system $\tilde{H}^1(-1,1)$, and moreover, are an orthonormal basis.

Now consider the case of $d \geq 2$ dimensions.  Define the weighted Sobolev space
\be{
\label{tildeH1_def}
\tilde{H}^1(D) = \left \{ f \in L^2_{\nu}(D) : \partial_k f \in L^2_{\nu_k}(D),\ k=0,\ldots,d \right \},
}
where $\nu_k(\bm{y})$ is the weight function given by
\bes{
\nu_0(\bm{y}) = \nu(\bm{y}) = \prod^{d}_{j=1} \nu(y_j),\qquad \nu_{k}(\bm{y}) = \chi(y_k) \prod^{d}_{\substack{j=1 \\ j \neq k}} \nu(y_j),\quad k=1,\ldots,d.
}
The associated norm and inner product are
\bes{
\nm{f}^2_{\tilde{H}^1(D)} = \sum^{d}_{k=0} \nm{\partial_k f}^2_{L^2_{\nu_k}(D)},
\qquad
\ip{f}{g}_{\tilde{H}^1(D)}= \sum^{d}_{k=0} \ip{\partial_k f}{\partial_k g}_{L^2_{\nu_k}(D)},
}
respectively.  Furthermore, the functions
\bes{
\psi_{\bm{n}}(\bm{y}) = \frac{1}{\sqrt{1+\lambda_{\bm{n}}}} \phi_{\bm{n}}(\bm{y}),\quad \bm{n} \in \bbN^d_0,
}
where
\be{
\label{lambda_def}
\lambda_{\bm{n}} = \sum^{d}_{k=1} \lambda_{n_k},
}
constitute an orthonormal basis of $\tilde{H}^1(D)$.

Since it will be useful later, we now make one further observation.  Let $g \in \tilde{H}^1(D)$.  Since $g \in L^2_{\nu}(D)$ by assumption, we may write
\bes{
g = \sum_{\bm{n} \in \bbN^d_0} x_{\bm{n}} \phi_{\bm{n}},
\qquad x_{\bm{n}} = \ip{g}{\phi_{\bm{n}}}_{L^2_{\nu}(D)},
}
so that
\bes{
\nm{g}^2_{L^2_{\nu}(D)} =\sum_{\bm{n} \in \bbN^d_0} | x_{\bm{n}} |^2.
}
However, due to the orthogonality relations, the coefficients of $g$ with respect to the basis $\psi_{\bm{n}}$ are
\bes{
\ip{g}{\psi_{\bm{n}}}_{\tilde{H}^1(D)} = \sqrt{1+\lambda_{\bm{n}}} x_{\bm{n}}.
}
In particular,
\bes{
\nm{g}^2_{\tilde{H}^1(D)} = \sum_{\bm{n} \in \bbN^d_0} (1+\lambda_{\bm{n}}) | x_{\bm{n}} |^2. 
}

\subsection{The gradient-enhanced linear system}\label{ss:gradenhance}
We are now in a position to formulate the gradient-enhanced recovery problem.  First, following the notation of \S \ref{s:prelims}, we define the matrices
\bes{
A_{k} = \frac{1}{\sqrt{m}} \left ( \frac{\partial \phi_{\bm{n}_j}(\bm{y}_i)}{\partial y_k}  \right )^{m,N}_{i=1,j=1} \in \bbC^{m \times N},\qquad k=0,\ldots,d.
}
Here and elsewhere, when $k = 0$ we mean that no partial derivative is taken, i.e.
\bes{
A_0 =  \frac{1}{\sqrt{m}} \left ( \phi_{\bm{n}_j}(\bm{y}_i)  \right )^{m,N}_{i=1,j=1}  \in \bbC^{m \times N}.
}
Recall that $\bm{x}_{\Lambda}$ denotes the vector of coefficients of $f$ corresponding to the index set $\Lambda$.  Therefore
\bes{
 \frac{1}{\sqrt{m}} \left ( \partial_k f(\bm{y}_i) \right )^{m}_{i=1} = A_k \bm{x}_{\Lambda} +  \frac{1}{\sqrt{m}}  \left ( \partial_k e_{\Lambda}(\bm{y}_i)  \right )^{m}_{i=1},
}
where $e_{\Lambda}$ is as in \R{eLambda}.  For reasons that will become clear in a moment, we let
\bes{
\bar{A} =  \begin{bmatrix} T_0 A_0 \\ T_1 A_1 \\ \vdots \\ T_d A_d \end{bmatrix} \in \bbC^{(d+1) m \times N},
}
where  $T_k = \diag \left( \left ( \sqrt{\tau_k(\bm{y}_i)} \right )^{m}_{i=1} \right ) \in \bbC^{m \times m}$,
are diagonal scaling matrices, and the $\tau_k$ are given by
\bes{
\tau_0(\bm{y}) = \frac{\prod^{d}_{j=1} \nu(y_j)}{\prod^{d}_{j=1} \mu(y_j)}=\frac{\nu_0(\bm{y})}{\mu(\bm{y})},\qquad \tau_k(\bm{y}) = \frac{ \chi(y_k) \prod^{d}_{{j=1, j \neq k}} \nu(y_j)}{\prod^{d}_{j=1} \mu(y_j)} = \frac{\nu_k(\bm{y})}{\mu(\bm{y})},\quad k=1,\ldots,d.
}
As we will show in \S \ref{matrescale}, the diagonal scaling matrices $T_k$ are used to ensure that $\bar{A}^*\bar{A}$ is diagonal in expectation, which is important for the subsequent analysis.
With this in hand, we can write the linear system of the gradient-augmented recovery problem as
\be{
\label{GElinsys}
 \bm{f} = \bar{A} \bm{x}_{\Lambda} + \bm{e},
}
where
\be{
\label{y_def}
\bm{f} = \begin{bmatrix} \bm{f}_0 \\ \vdots \\ \bm{f}_d \end{bmatrix},\qquad \bm{f}_k = \frac{1}{\sqrt{m}} \left (\sqrt{ \tau_k(\bm{y}_i)}\partial_k f(\bm{y}_i) \right )^{m}_{i=1},
}
and
\bes{
\bm{e} = \begin{bmatrix} \bm{e}_0 \\ \vdots \\ \bm{e}_d \end{bmatrix},\qquad \bm{e}_k =  \frac{1}{\sqrt{m}} \left ( \sqrt{\tau_k(\bm{y}_i)} \partial_k e_{\Lambda}(\bm{y}_i) \right )^{m}_{i=1}.
}
As in \S \ref{ss:weightedl1}, we shall assume that the tail error satisfies
\be{
\label{eta_def}
\nm{\bm{e}}_{2} \leq \eta,
}
for some known $\eta \geq 0$.  Note that this is implied by the condition
\bes{
\sup_{y \in D} \sum^{d}_{k=0} \tau_k(y) | \partial_k e_{\Lambda}(y) |^2 \leq \eta^2.
}

\subsection{Matrix scaling, problem formulation and Sobolev norm error bounds}\label{matrescale}

Recall that the points $\bm{y}_1,\ldots,\bm{y}_m$ are independently and identically distributed according to $\mu$.  Due to the diagonal scaling matrices and the Sobolev orthogonality of the basis functions, we have
\eas{
 \bbE \left ( \bar{A}^* \bar{A} \right )_{\bm{n},\bm{n'}} &= \sum^{d}_{k=0} \int_{D} \partial_{k} \phi_{\bm{n}}(\bm{y}) \overline{ \partial_{k} \phi_{\bm{n'}}(\bm{y}) } \tau_k(\bm{y}) \mu(\bm{y}) \D y
 \\
 &= \sum^{d}_{k=0} \int_{D} \partial_{k} \phi_{\bm{n}}(\bm{y}) \overline{ \partial_{k} \phi_{\bm{n'}}(\bm{y}) } \nu_k(\bm{y}) \D y 
 = \left ( 1 + \lambda_{\bm{n}} \right ) \delta_{\bm{n},\bm{n'}}.
}
For this reason, we introduce the diagonal scaling matrix
$
Q = \diag \left ( \sqrt{1+\lambda_{\bm{n}}} \right)_{\bm{n} \in \Lambda},
$
so that the scaled matrix
\be{
\label{A_deriv_prob}
A = \bar{A} Q^{-1},
}
satisfies $\bbE(A^*A) = I$.
With this in hand, we are now in a position to formulate the gradient-augmented weighted $\ell^1$ minimization problem:
\be{
\label{l1u_min_deriv}
\min_{\bm{z} \in \bbC^N} \nm{\bm{z}}_{1,\bm{w}}\ \mbox{subject to $\nm{A \bm{z} - \bm{y}}_2 \leq \eta$}.
}
Note that if $\bm{\hat{z}}$ is a minimizer of this problem, then we define $\hat{\bm{x}} = Q^{-1} \hat{\bm{z}}$ as the approximation to the true coefficients $\bm{x}_{\Lambda}$, and let
\bes{
\hat{f} = \sum_{\bm{n} \in \Lambda} \hat{x}_{\bm{n}} \phi_{\bm{n}},
}
be the corresponding approximation to $f$.

Finally, we note the following.  If $f_{\Lambda}$ is as in \R{eLambda}, then, due to the Sobolev orthogonality,
\eas{
\nm{f_{\Lambda} - \hat{f}}_{\tilde{H}^1(D)} = \nm{Q(\bm{x}_{\Lambda} - \bm{\hat{x}})}_{2} = \nm{\bm{z}_{\Lambda} - \bm{\hat{z}}}_{2},
}
where $\bm{z}_{\Lambda} = Q \bm{x}_{\Lambda}$ are the coefficients of $f$ with respect to the Sobolev-orthogonal basis $\{ \psi_{\bm{n}} \}_{\bm{n} \in \bbN^d_0}$.  Thus, since the analysis of the problem \R{l1u_min_deriv} will provide a bound for $\nm{\bm{z}_{\Lambda} - \bm{\hat{z}}}_{2}$, we correspondingly obtain a bound for the approximation error in the Sobolev-type norm $\tilde{H}^1(D)$.

\section{Main results}\label{s:mainres}

In order to state our main results, we require several additional definitions.  First, given weights $\bm{w} = (w_{\bm{n}})_{\bm{n} \in \bbN^d_0}$ and a set $\Delta \subset \bbN^d_0$ 
we define the weighted cardinality of $\Delta$ as 
\be{
|\Delta|_{\bm{w}} = \sum_{\bm{n}\in \Delta} w_{\bm{n}}^2.
}
Second, given $\nu$, $\mu$ and $\{ \phi_{\bm{n}} \}_{\bm{n} \in \bbN^d_0}$ as in \S \ref{s:prelims} and \S \ref{s:gradaug}, we define the \textit{intrinsic weights} $\bm{u} = (u_{\bm{n}})_{\bm{n} \in \mathbb{N}_0^d}$ as
 \be{
 \label{u_def}
 u_{\bm{n}} = \sup_{\bm{y} \in D} \sqrt{\nu(\bm{y})/\mu(\bm{y})} |\phi_{\bm{n}}(\bm{y})|.
}
Third, we let
\be{
\label{kappa_def_1D}
\kappa_{n} = u^{-2}_{n} \sup_{y \in (-1,1)} \frac{\chi(y)}{\mu(y)} | \phi'_{n}(y) |^2,\qquad n = 0,1,2,\ldots,
}
and for $\bm{n} \in \bbN^d_0$, we set
\be{
\label{kappa_def}
\kappa_{\bm{n}} = \sum^{d}_{k=1} \kappa_{n_k}.
}
Finally, given $\bm{x}_{\Lambda}$ we let $\bm{x}_{\Delta} \in \bbC^{N}$ be the vector obtained from $\bm{x}_{\Lambda}$ by setting all terms corresponding to indices $\bm{n} \in \Lambda \backslash \Delta$ to zero.

\subsection{General recovery guarantees}

Our first result is as follows:

\thm{
\label{t:main_res_derivs_identical}
Let $\Lambda \subset \bbN^d_0$ with $|\Lambda| = N \geq 2$, $0 < \epsilon < 1$,  $\bm{w} \in \bbR^N$ be a vector of weights with $w_{\bm{n}} \geq 1$, $\forall \bm{n}$, $\Delta \subset \Lambda$, $|\Delta | \geq 2$ and $f = \sum_{\bm{n} \in \bbN^d_0} x_{\bm{n}} \phi_{\bm{n}} \in \tilde{H}^1(D)$, where $D=(-1,1)^d$ and $\tilde{H}^1(D)$ is as in \R{tildeH1_def}.  Let 
\be{
\label{sampcomp1}
m \gtrsim \max_{\bm{n} \in \Lambda} \left \{ \frac{1+\kappa_{\bm{n}}}{1+\lambda_{\bm{n}} } \right \} \cdot \left ( | \Delta |_{\bm{u}} + \max_{\bm{n} \in \Lambda} \left \{ \frac{u^2_{\bm{n}} }{w^2_{\bm{n}} } \right \} | \Delta |_{\bm{w}} \right ) \cdot L,
}
where
\bes{
L = \log(N/\epsilon) + \log(|\Delta|_{\bm{w}}) \cdot \log(|\Delta|_{\bm{w}} / \epsilon ),
}
draw $\bm{y}_{1},\ldots,\bm{y}_m$ independently according to the density $\mu$, and let $A$, $\bm{f}$ and $\eta$ be as in \R{A_deriv_prob}, \R{y_def} and \R{eta_def} respectively.  Then, if $\hat{\bm{z}}$ is any minimizer of \R{l1u_min_deriv} and $\hat{\bm{x}} = Q^{-1} \hat{\bm{z}}$, the approximation $\hat{f} = \sum_{\bm{n} \in \Lambda} \hat{x}_{\bm{n}} \phi_{\bm{n}}$ satisfies
\bes{
\nmu{f - \hat{f} }_{\tilde{H}^1(D)}  \lesssim \nmu{f - f_{\Lambda}}_{\tilde{H}^1(D)} +  \nm{\bm{x}_{\Lambda} - \bm{x}_{\Delta}}_{1,\bm{v}} + \sqrt{|\Delta|_{\bm{w}}} \eta,
}
with probability at least $1-\epsilon$, where $v_{\bm{n}} = \sqrt{1+\lambda_{\bm{n}}} w_{\bm{n}}$, $\bm{n} \in \bbN^d_0$.
}

This result is understood as follows.  For a fixed function $f$ with coefficients $\bm{x}$, and a fixed set $\Delta$, by drawing $m$ samples according to $\mu$, with $m$ given by \R{sampcomp1}, we can recover $f$ up to an error (measured in a Sobolev norm) depending on how well $\bm{x}$ is approximated by its coefficients with indices in $\Delta$ (the term $\nm{\bm{x}_{\Lambda} - \bm{x}_{\Delta}}_{1,\bm{v}}$).  As with the other results in this section, this is a type of \textit{nonuniform} recovery guarantee; see \S \ref{ss:discussion}.

Note that this result makes no assumptions on $\Delta$.  In a moment however, we shall specialize it to the case of lower sets (recall \S \ref{ss:lowersets}).
First, however, we note an immediate consequence of Theorem \ref{t:main_res_derivs_identical}.  Namely, in order to minimize the right-hand side of \R{sampcomp1}, the weights $\bm{w}$ should be chosen as
\bes{
\bm{w} = \bm{u}.
}
That is, the best optimization weights are precisely the intrinsic weights \R{u_def}.  This is identical to a conclusion reached in \cite{AdcockCSFunInterp} for the unaugmented problem. 

With this in hand, we now consider recovery in lower sets:

\cor{
\label{c:main_cor}
Let $s \geq 2$, $\Lambda = \Lambda^{\mathrm{HC}}_{s}$ be the hyperbolic cross index set \R{HCset}, $0 < \epsilon < 1$ and $f = \sum_{\bm{n} \in \bbN^d_0} x_{\bm{n}} \phi_{\bm{n}} \in \tilde{H}^1(D)$, where $D=(-1,1)^d$ and $\tilde{H}^1(D)$ is as in \R{tildeH1_def}.  Suppose that
\be{
\label{sampcomp2}
m \gtrsim \max_{\bm{n} \in \Lambda} \left \{ \frac{1+\kappa_{\bm{n}}}{1+\lambda_{\bm{n}} } \right \}  \cdot K(s) \cdot L.
}
where $L = \left ( \min \{ d+ \log(s/\epsilon) , \log(2d) \log(s/\epsilon) \} + \log(K(s)) \cdot \log(K(s) / \epsilon ) \right )$,
\be{
\label{Ksdef}
K(s) =  \max \left \{ | \Delta |_{\bm{u}} : \mbox{$|\Delta | \leq s$ and $\Delta$ is lower} \right \},
}
and $\bm{u}$ are the weights defined in \R{u_def}.  Draw $\bm{y}_{1},\ldots,\bm{y}_m$ independently according to the density $\mu$, let $A$, $\bm{f}$ and $\eta$ be as in \R{A_deriv_prob}, \R{y_def} and \R{eta_def} respectively and set $\bm{w} = \bm{u}$.  Then, if $\hat{\bm{z}}$ is any minimizer of \R{l1u_min_deriv} and $\hat{\bm{x}} = Q^{-1} \hat{\bm{z}}$, the approximation $\hat{f} = \sum_{\bm{n} \in \Lambda} \hat{x}_{\bm{n}} \phi_{\bm{n}}$ satisfies
\bes{
\nmu{f - \hat{f} }_{\tilde{H}^1(D)} \lesssim \nmu{f - f_{\Lambda}}_{\tilde{H}^1(D)} + \sigma_{s,L}(\bm{x}_{\Lambda})_{1,\bm{v}} + \sqrt{K(s)} \eta,
}
with probability at least $1-\epsilon$, where $\sigma_{s,L}(\cdot)_{1,\bm{v}}$ is as in \R{sigma_sL} and $v_{\bm{n}} = \sqrt{1+\lambda_{\bm{n}}} u_{\bm{n}}$, $\bm{n} \in \bbN^d_0$.
}

Specializing Theorem \ref{t:main_res_derivs_identical}, the error estimate in this result is given in terms of the best $s$-term approximation error in lower sets $\sigma_{s,L}(\cdot)_{1,\bm{v}}$.  However, the sample complexity estimate \R{sampcomp2} is not given completely explicitly in terms of the sparsity $s$ and dimension $d$.  For this, we need estimate the quantities $(1+\kappa_{\bm{n}})/(1+\lambda_{\bm{n}})$, $\bm{n} \in \Lambda$, and $K(s)$, and this requires the basis $\{ \phi_{\bm{n}} \}$ and sampling density $\mu$ to be specified.  We do this next.

\subsection{The case of Jacobi polynomials with $\mu = \nu$}

Consider the Jacobi polynomial basis (recall \S \ref{s:SL}) and sampling density $\mu = \nu$.  We have

\cor{
\label{c:main_res_jacobi_identical_lower}
Consider setup of Corollary \ref{c:main_cor}, where $\{ \phi_{\bm{n}} \}_{\bm{n} \in \bbN^d_0}$ is the tensor-product Jacobi polynomial basis with parameters $\alpha,\beta \geq - 1/2$ and $\mu = \nu$.  Suppose that
\bes{
m \gtrsim K(s) \cdot \left ( \min \{ d+ \log(s/\epsilon) , \log(2d) \log(s/\epsilon) \} + \log(K(s)) \cdot \log(K(s) / \epsilon ) \right ),
}
where 
\be{
\label{Kdef}
K(s) = K^{(\alpha,\beta)}(s) = \max \left \{ | \Delta |_{\bm{u}} : \mbox{$|\Delta | \leq s$ and $\Delta$ is lower} \right \}.
}
Then, if $\hat{\bm{z}}$ is any minimizer of \R{l1u_min_deriv} and $\hat{\bm{x}} = Q^{-1} \hat{\bm{z}}$, the approximation $\hat{f} = \sum_{\bm{n} \in \Lambda} \hat{x}_{\bm{n}} \phi_{\bm{n}}$ satisfies
\bes{
\nmu{f - \hat{f} }_{\tilde{H}^1(D)} \lesssim \nmu{f - f_{\Lambda}}_{\tilde{H}^1(D)}  + \sigma_{s,L}(\bm{x}_{\Lambda})_{1,\bm{v}} + \sqrt{K(s)} \eta,
}
with probability at least $1-\epsilon$, where $\sigma_{s,L}(\cdot)_{1,\bm{v}}$ is as in \R{sigma_sL} and $v_{\bm{n}} = \sqrt{1+\lambda_{\bm{n}}} u_{\bm{n}}$, $\bm{n} \in \bbN^d_0$.
}

The proof of this corollary involves showing that $\kappa_{\bm{n}} \lesssim \lambda_{\bm{n}}$ for the Jacobi polynomials.  See \S \ref{ss:corproof} for details. Having done this, the sample complexity is determined up to magnitude of $K(s)$, which depends on the indices $\alpha,\beta$ of the Jacobi polynomials.  For certain values of 
$\alpha$ and $\beta$, we have the following result (see \cite{MiglioratiJAT}):

\thm{
\label{K_bnd}
Let $K(s) = K^{(\alpha,\beta)}(s)$ be as in \R{Kdef}. Then the following hold:
\begin{enumerate}
\item[(i)] if $\alpha,\beta \in \bbN_0$  then $K(s) \leq s^{2 \max \{ \alpha,\beta \} + 2 }$ ,
\item[(ii)]  if $\beta = \alpha$ and $2 \alpha + 1 \in \bbN$ then $K(s) \leq s^{2 \alpha + 2 }$,
\item[(iii)] if $\alpha = \beta = -1/2$ then $K(s) \leq s^{\log(3)/\log(2)}$ .
\end{enumerate}
In particular, $K(s) \leq s^2$ for Legendre polynomials ($\alpha = \beta = 0$) and $K(s) \leq s^{\log(3)/\log(2)}$ for Chebyshev polynomials ($\alpha = \beta = -1/2$).
}

This result implies that, for values of $\alpha, \beta$ satisfying Theorem \ref{K_bnd}, the sample complexity reduces to an estimate of the form
\be{
\label{jacobiSC}
m \gtrsim s^{\gamma} \cdot \log(2d)\cdot \log^2(s/\epsilon),
}
where $\gamma \geq 1$ depends on $\alpha$ and $\beta$ -- in other words, polynomial in $s$ and logarithmic in the dimension $d$.  Hence the curse of dimensionality is mitigated to a substantial extent.  Up to constants and log factors, this is the same as the unaugmented case.  See \S \ref{ss:discussion} for further discussion.

%\subsection{Sampling from a measure other than the orthogonality measure, i.e. $\mu(\bm{y})  \neq \nu_0(\bm{y})$} 

\subsection{Legendre polynomials and preconditioning}

In the previous section, we considered sampling with the same density as the orthogonality density $\nu$.  A number of settings call for the use of a different sampling density $\mu$.  In particular, the case where $\phi_{\bm{n}}$ are the Legendre polynomials and $\mu$ is the Chebyshev density has been studied in \cite{AdcockCSFunInterp,RauhutWard,YanGuoXui_l1UQ}, where it is referred to as \textit{preconditioning}.    For this case we have the following:

\cor{
\label{c:main_res_leg_lower}
Let $\mu$ be the tensor Chebyshev density, $\nu$ be the uniform density, $s \geq 2$, $\Lambda = \Lambda^{\mathrm{HC}}_{s}$ be the hyperbolic cross index set \R{HCset}, $0 < \epsilon < 1$, $\bm{u}$ be the weights defined in \R{u_def} and $f = \sum_{\bm{n} \in \bbN^d_0} x_{\bm{n}} \phi_{\bm{n}} \in \tilde{H}^1(D)$, where $D=(-1,1)^d$ and $\tilde{H}^1(D)$ is as in \R{tildeH1_def}.  Suppose that
\be{
\label{sampcompLeg}
m \gtrsim \min \left\{2^ds, (\pi/2)^d s^{\log(1+4/\pi)/\log(2)} \right \} \cdot (d + \log(s))\cdot (d+\log(s/\epsilon) ).
}
Draw $\bm{y}_{1},\ldots,\bm{y}_m$ independently according to the density $\mu$, let $A$, $\bm{f}$ and $\eta$ be as in \R{A_deriv_prob}, \R{y_def} and \R{eta_def} respectively and set $\bm{w} = \bm{u}$.  Then, if $\hat{\bm{z}}$ is any minimizer of \R{l1u_min_deriv} and $\hat{\bm{x}} = Q^{-1} \hat{\bm{z}}$, the approximation $\hat{f} = \sum_{\bm{n} \in \Lambda} \hat{x}_{\bm{n}} \phi_{\bm{n}}$ satisfies
\bes{
\nmu{f - \hat{f} }_{\tilde{H}^1(D)} \lesssim \nmu{f - f_{\Lambda}}_{\tilde{H}^1(D)} + \sigma_{s,L}(\bm{x}_{\Lambda})_{1,\bm{v}} + \sqrt{K(s)} \eta,
}
with probability at least $1-\epsilon$, where $\sigma_{s,L}(\cdot)_{1,\bm{v}}$ is as in \R{sigma_sL} and $v_{\bm{n}} = \sqrt{1+\lambda_{\bm{n}}} u_{\bm{n}}$, $\bm{n} \in \bbN^d_0$.
}

\subsection{Discussion}\label{ss:discussion}

We now compare our results to those obtained in \cite{AdcockCSFunInterp} for the problem of recovery from function samples only.  Using the same setup and notation, in \cite{AdcockCSFunInterp} it was proved that if
\be{
\label{nogradSC}
m \gtrsim K(s) \cdot \log(\epsilon^{-1}) \cdot L,\qquad L = \left( \min \{ \log(2s) + d , \log(2d) \log(2s) \} + \log(K(s)) \right ),
}
where $K(s)$ is as in \R{Ksdef}, then the recovery error satisfies
\be{
\label{nograderr}
\| f - \hat{f} \|_{L^2(D)} \lesssim \sigma_{s,L}(\bm{x}_{\Lambda})_{1,\bm{u}} + \| f - f_{\Lambda} \|_{L^2(D)} + \lambda \sqrt{K(s)} \eta ,
}
with high probability, where $\lambda = 1 + \sqrt{\log(\epsilon^{-1})}/L$ (see Theorem 6.1 and Remark 7.9 of \cite{AdcockCSFunInterp}).  The main point is that the sample complexity estimates in Corollaries \ref{c:main_res_jacobi_identical_lower} and \ref{c:main_res_leg_lower} are identical, up to minor changes in the log factor\footnote{As we discuss in \S \ref{s:proofs}, we use a slightly different method of proof to remove the factor $\lambda$ in the error bound, at the expense of a slightly increased log factor.}, to those obtained in \cite{AdcockCSFunInterp}.  In particular, \R{nogradSC} reduces to
\bes{
m \gtrsim s^{\gamma} \cdot \log(\epsilon^{-1}) \cdot \log(2d) \cdot \log(2s),
}
in the case of Jacobi polynomials, as in \R{jacobiSC} (a similar statement can be made concerning Corollary \ref{c:main_res_leg_lower}).  However, the error in the gradient-augmented case is bounded in the stronger Sobolev norm, as opposed to the $L^2(D)$ norm in \R{nograderr}.

Similar to those of \cite{AdcockCSFunInterp}, the results of this section are \textit{nonuniform} recovery guarantees: they ensure recovery of a single $f$ from a random draw of sample points.  For the 
 unaugmented case, \textit{uniform} recovery guarantees for Chebyshev and Legendre polynomials (with $\mu = \nu$) have been proved in \cite{BASBCWMatheon,Chkifalower}.  The corresponding sample complexity estimates are similar to \R{nogradSC}, except with higher log factors.  Conversely, the error bound \R{nograderr} is improved by a factor of $1/\sqrt{K(s)}$.  This is typical for uniform recovery guarantees in compressed sensing.  We expect a similar uniform recovery guarantee is possible for the gradient-augmented setting, but we leave this as future work.

\rem{
The error bounds in the gradient-augmented setting measure the error in the $\ell^1_{\bm{v}}$-norm, for modified weights $v_{\bm{n}} = \sqrt{1+\lambda_{\bm{n}}} u_{\bm{n}}$, as opposed to the $\ell^1_{\bm{u}}$-norm for the case of function samples only.  This is quite natural.  First, we recall that error estimates in $\ell^2$-type norms are not generally possible in compressed sensing under optimal sample complexities (see, for example, \cite[Chpt.\ 11]{FoucartRauhut}).  Second, note that
\be{
\label{unaugerrbd}
\sup_{y \in D} | g(y) | \leq \| \bm{x} \|_{1,\bm{u}},\qquad g = \sum_{\bm{n} \in \bbN^d_0} x_{\bm{n}} \phi_{\bm{n}}.
}
Hence the $\ell^1_{\bm{u}}$-norm of $\bm{x}$ provides an upper bound for $\| g \|_{L^\infty(D)}$.  Similarly, one can show that
\be{
\label{augerrbd}
\sup_{y \in D} \sqrt{\sum^{d}_{k=0} \tau_k(y) | \partial_k g(y) |^2 } \leq \| \bm{x} \|_{1,\bm{v}}.
}
Hence the $\ell^1_{\bm{v}}$-norm of the coefficients provides an upper bound on a particular weighted $L^\infty$-type Sobolev norm.  
Now recall that in order to formulate the various optimization problems we introduce an error vector $\bm{e}$ which includes the expansion tail (see \R{unaug} and \R{GElinsys}).  In the unaugmented case, \R{unaugerrbd} gives that this vector satisfies the bound
\bes{
\| \bm{e} \|_{2} \leq \| \bm{x} - \bm{x}_{\Lambda} \|_{1,\bm{u}},
}
and in the augmented case \R{augerrbd} gives
\bes{
\| \bm{e} \|_{2} \leq \| \bm{x} - \bm{x}_{\Lambda} \|_{1,\bm{v}}.
}
In other words, the $\ell^1_{\bm{u}}$- and $\ell^1_{\bm{v}}$-norms are tight weighted $\ell^1$-norm bounds for the error vector in terms of the expansion coefficients.
}

\subsection{Sparse trigonometric polynomial approximations}

To complete this section, we note that this approach can be easily extended to other related Sturm--Liouville eigenfunctions, of both singular and regular types.  Of particular importance is the case of trigonometric polynomial expansions, equivalent to approximations in the Fourier basis
\be{
\label{Fourier}
\phi_{\bm{n}}(\bm{y}) = \exp \left (\I \pi (n_1 y_1 + \ldots + n_d y_d ) \right ),\qquad \bm{y} \in D = \bbT^d,\ \bm{n} \in \bbZ^d,
}
where $\bbT^d = [-1,1)^d$ is the unit $d$-torus.
These one-dimensional basis consists of eigenfunctions of the regular Sturm--Liouville problem with periodic boundary conditions and $\nu(y) = \chi(y) = 1/2$, $\zeta(y) = 0$.  The eigenvalues are $\lambda_{n} = n^2 \pi^2$.  Correspondingly, the scaled functions $\phi_{\bm{n}}(\bm{y}) / \sqrt{1+\lambda_{\bm{n}} }$ are an orthonormal basis of the periodic Sobolev space $H^1(\bbT)$.

With this in hand, the following is a straightforward consequence of Theorem \ref{t:main_res_derivs_identical}:

\cor{
Let $\Lambda \subset \bbZ^d$ with $| \Lambda | = N \geq 2$, $2 \leq s \leq N$, $0 < \epsilon < 1$ and $f = \sum_{\bm{n} \in \bbZ^d} x_{\bm{n}} \phi_{\bm{n}} \in H^1(\bbT)$, where $\phi_{\bm{n}}$ is the Fourier basis \R{Fourier}.  Let
\bes{
m \gtrsim s \cdot L,\qquad L = \log(N/\epsilon) + \log(s) \cdot \log(s/\epsilon).
}
draw $\bm{y}_1,\ldots,\bm{y}_m$ independently according to the uniform density $\mu(\bm{y}) = 2^{-d}$ and let $A$, $\bm{f}$ and $\eta$ be as in \R{A_deriv_prob}, \R{y_def} and \R{eta_def} respectively.  Then, if $\bm{\hat{x}} = Q^{-1} \hat{\bm{z}}$ where $\hat{\bm{z}}$ is any minimizer of \R{l1u_min_deriv} with weights $w_{\bm{n}} = 1$, $\forall \bm{n} \in \bbZ^d$, the approximation $\hat{f} = \sum_{\bm{n} \in \Lambda} \hat{x}_{\bm{n}} \phi_{\bm{n}}$ satisfies
\bes{
\nmu{f - \hat{f} }_{H^1(\bbT)}  \lesssim \nmu{f - f_{\Lambda}}_{H^1(\bbT)} + \sigma_{s}(\bm{x}_{\Lambda})_{1,\bm{v}} + \sqrt{s} \eta,
}
with probability at least $1-\epsilon$, where $v_{\bm{n}} = \sqrt{1+\lambda_{\bm{n}}}$, $\bm{n} \in \bbZ^d$.
}

Note the Fourier basis is uniformly bounded with $\| \phi_{\bm{n}} \|_{L^\infty(\bbT)} = 1$.  Hence in this case no lower set structure is required.  The corresponding sample complexity estimate scales linearly (and therefore optimally) in $s$.

\section{Numerical experiments} \label{s:numexp}
In this section, we wish to demonstrate the benefits of gradient-augmented sampling numerically for tensor Legendre and Chebyshev polynomials.
In order to do this, we shall assume that the computational cost of computing the gradient is roughly the same as the cost of computing function values.  This is reasonable in certain UQ applications, where $f$ is a quantity of interest of a parametric PDE and the gradient samples are computed via adjoint sensitivity analysis (for example).  See \cite{Penggradient} for further information.  For this reason, we model the total cost of computing the gradient-augmented measurements by
\be{
\label{cost}
\tilde m = m_o + m_g,
}
where $m_o$ is the number of function samples and $m_g$ is the number gradient samples. For the unaugmented 
problem, the computational cost is just $\tilde m = m_o$.

\begin{figure}[ht]
\centering
  \begin{tabular}{cc}
     \includegraphics[width=1.8in]{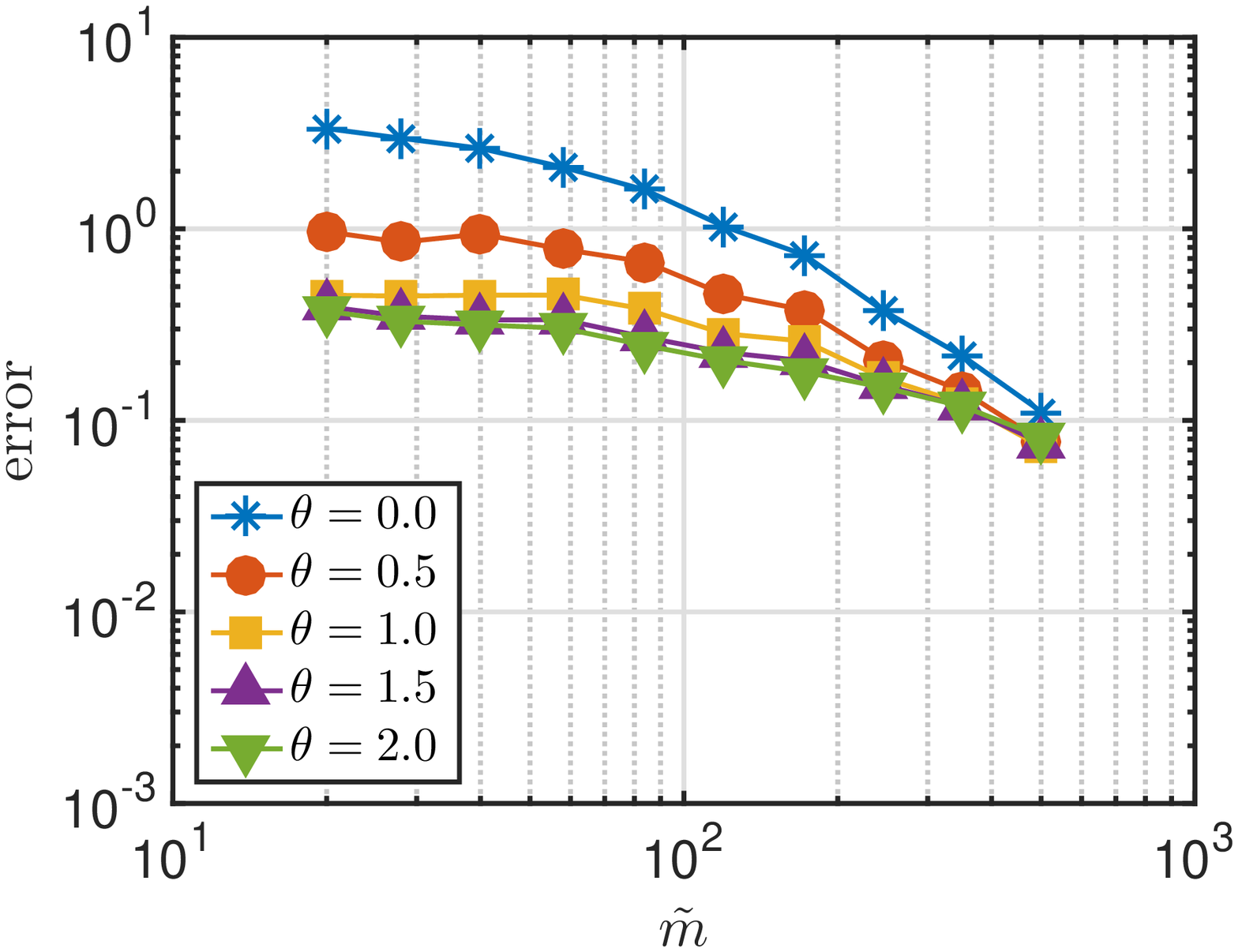} 
     \includegraphics[width=1.8in]{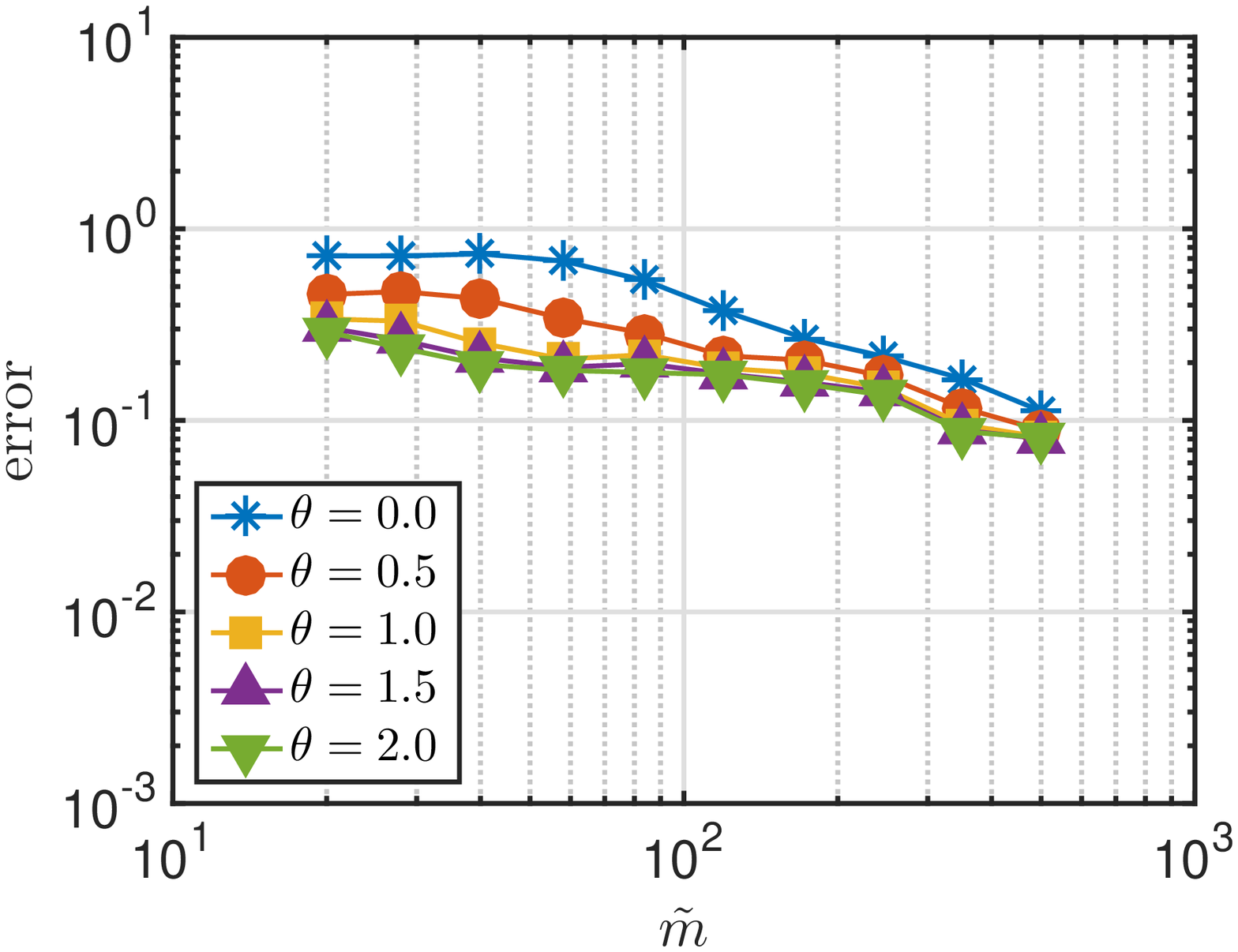}
     \includegraphics[width=1.8in]{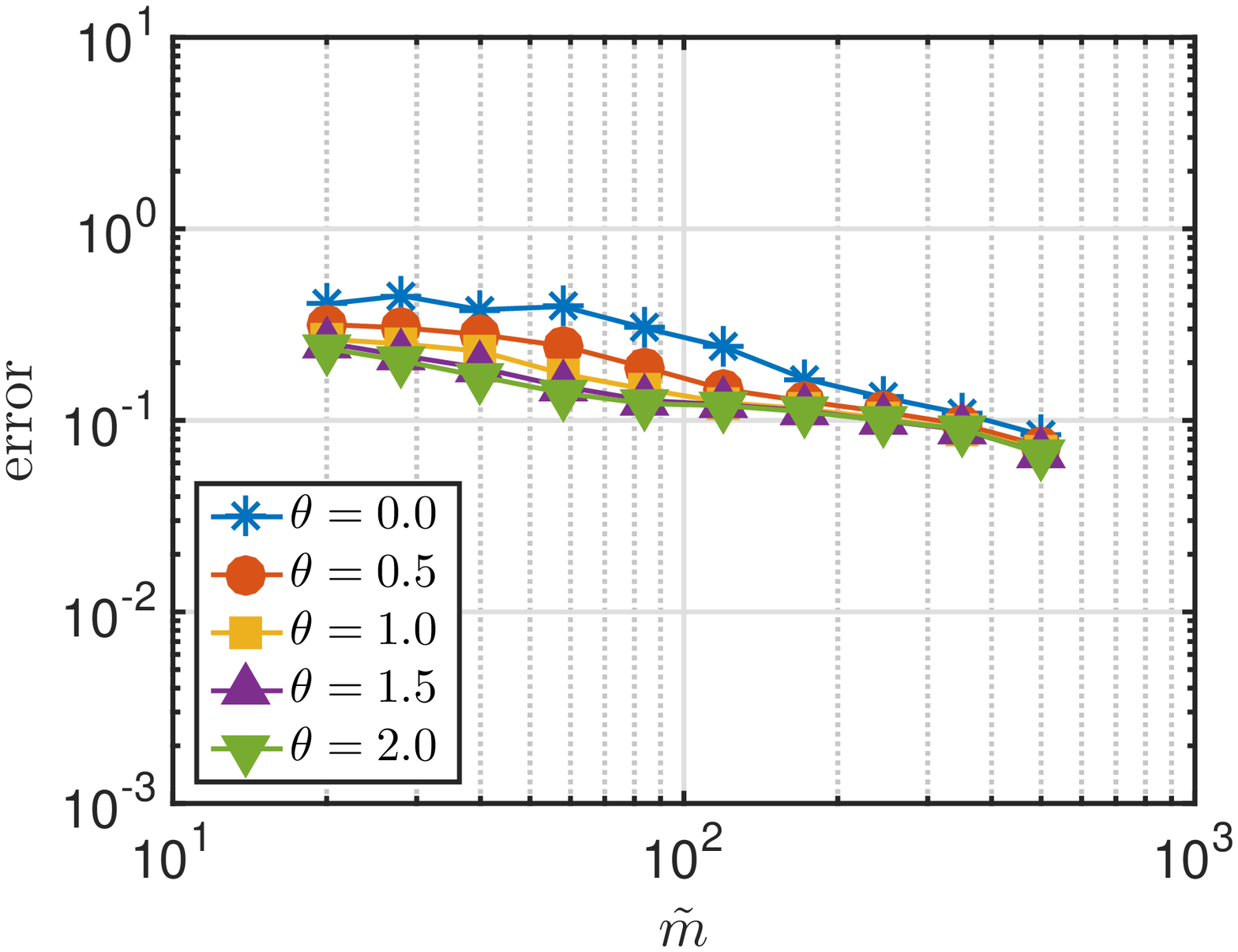}\\
      \includegraphics[width=1.8in]{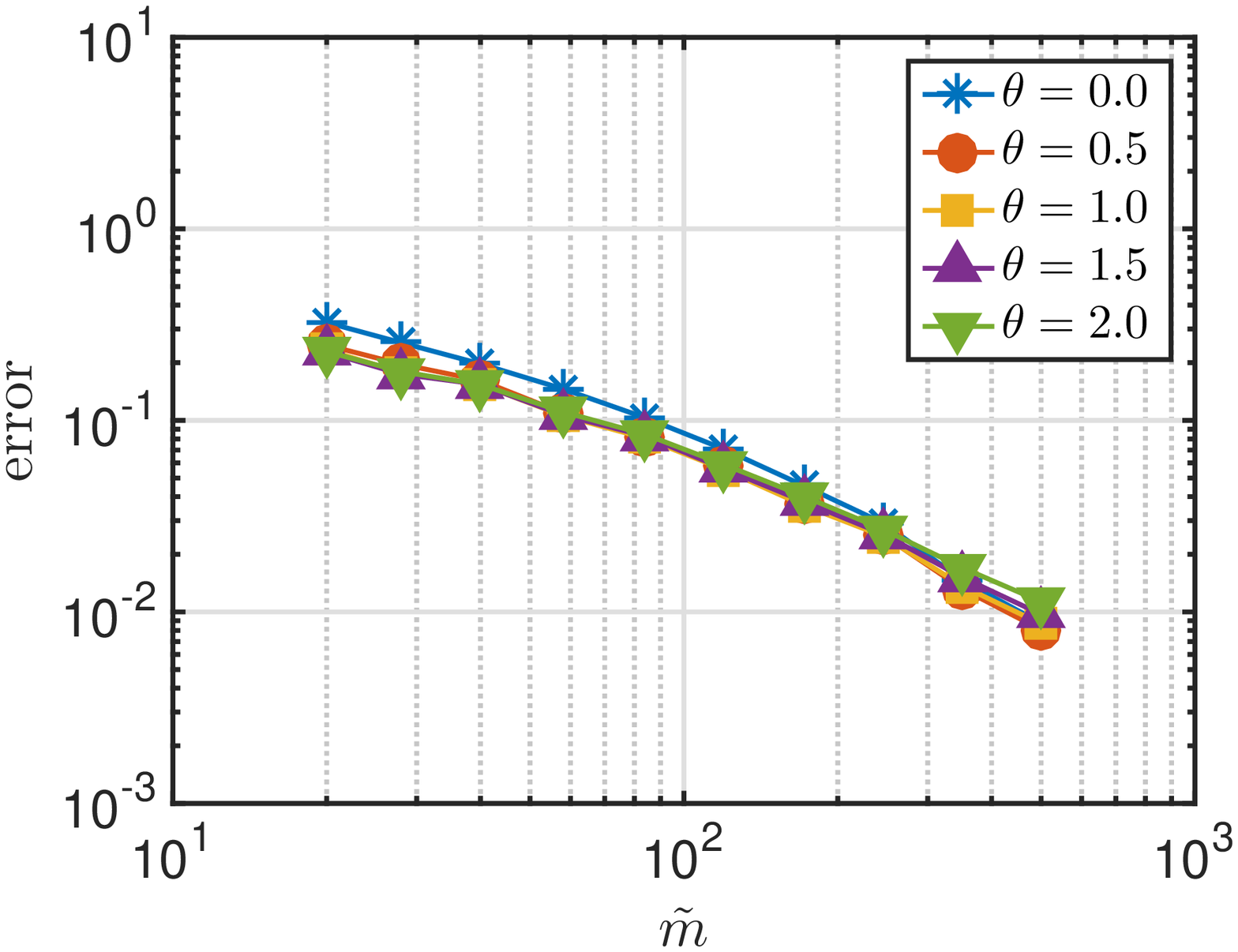} 
     \includegraphics[width=1.8in]{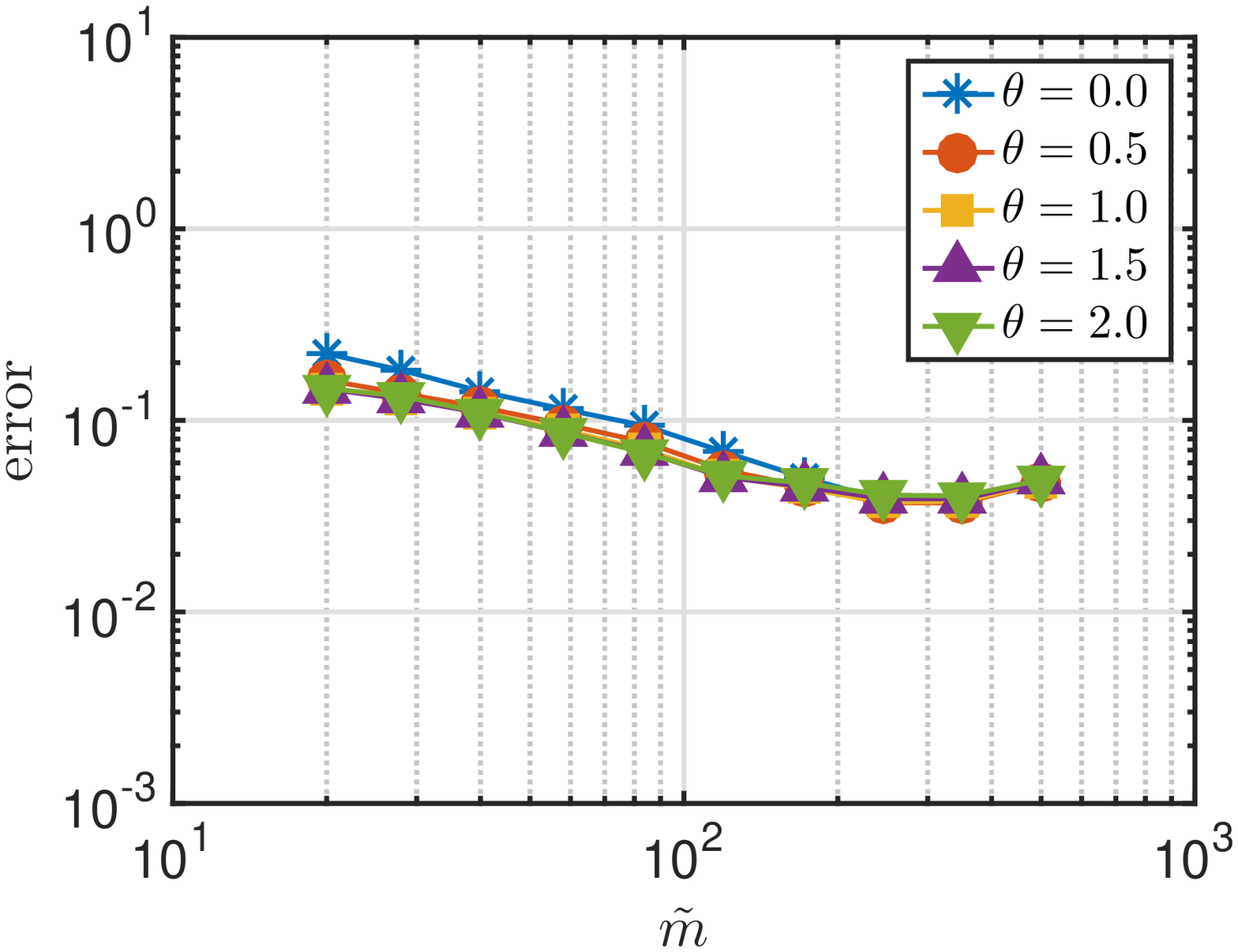}
      \includegraphics[width=1.8in]{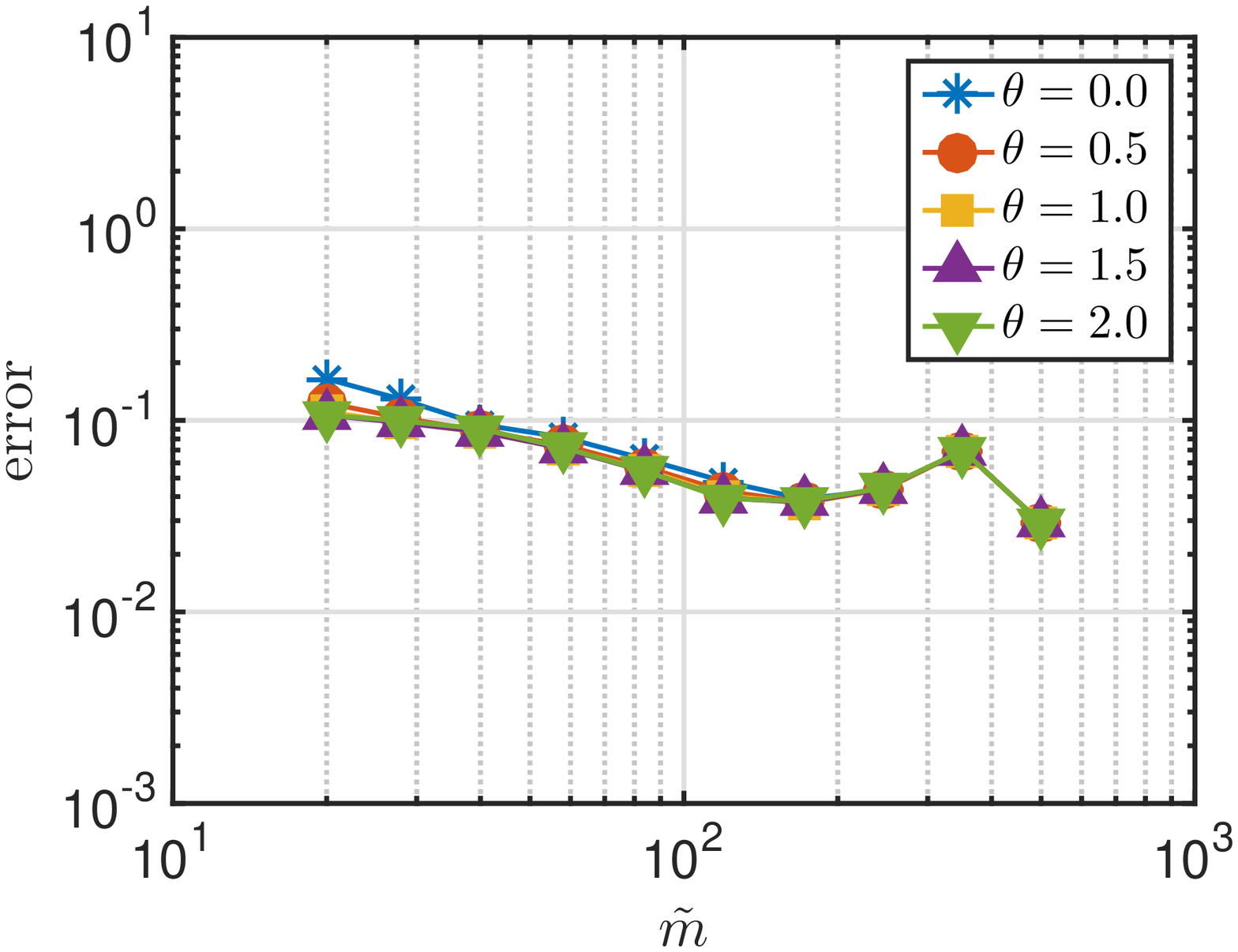}
    \end{tabular}
\caption{The error $\nmu{f_1-\tilde{f}_1}_{\tilde{H}^1(D)}$ against $\tilde{m}$ for Legendre polynomials with points drawn 
from the uniform density. From left to right, the values $(d, s) = (4, 72), (8, 23), (12, 14)$ 
were used. The unaugmented case is shown on the top row and the gradient-augmented case is shown on the bottom row. 
}
\label{LUoptpeak}
\end{figure}

\begin{figure}[h!]
\centering
  \begin{tabular}{cc}
     \includegraphics[width=1.8in]{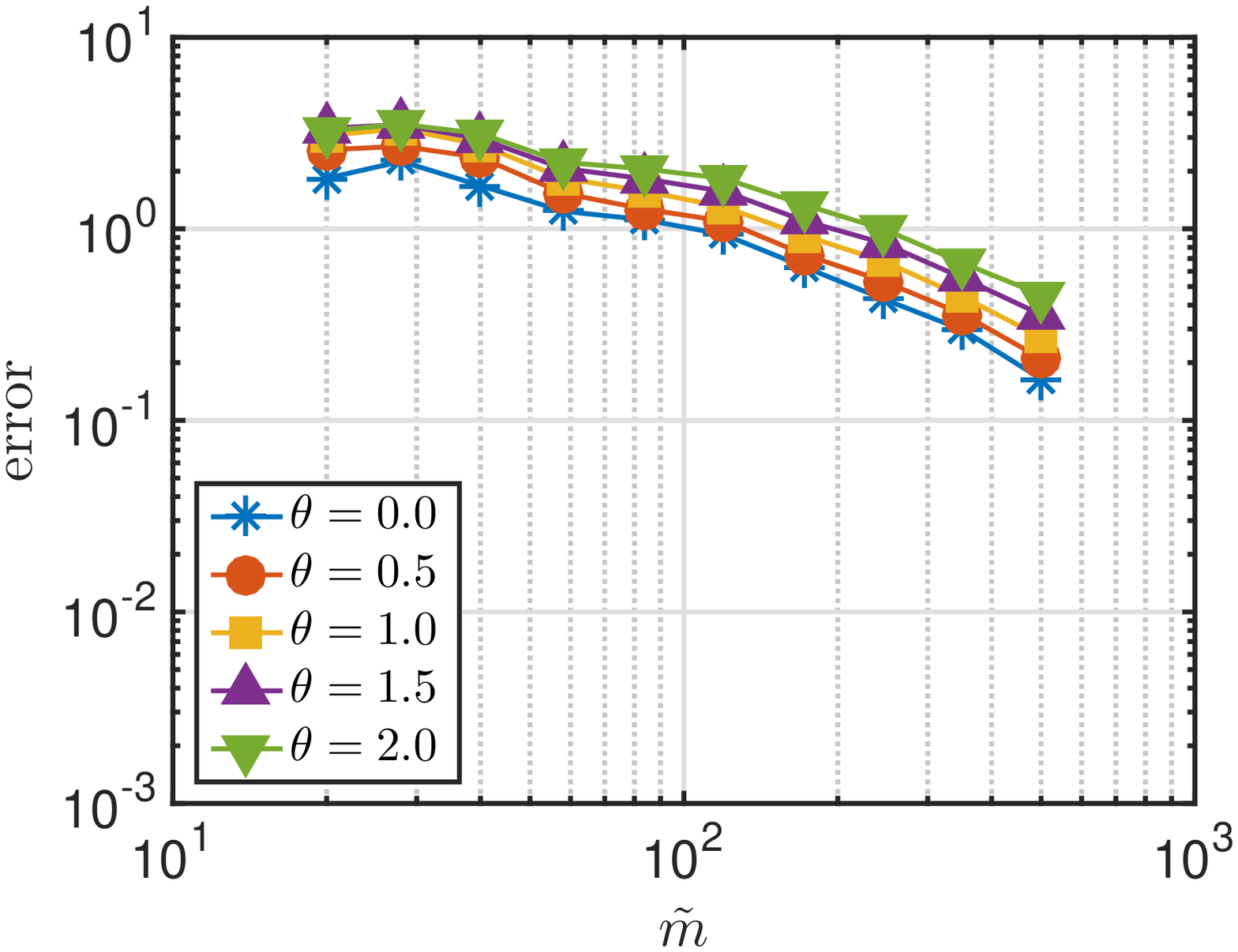} 
     \includegraphics[width=1.8in]{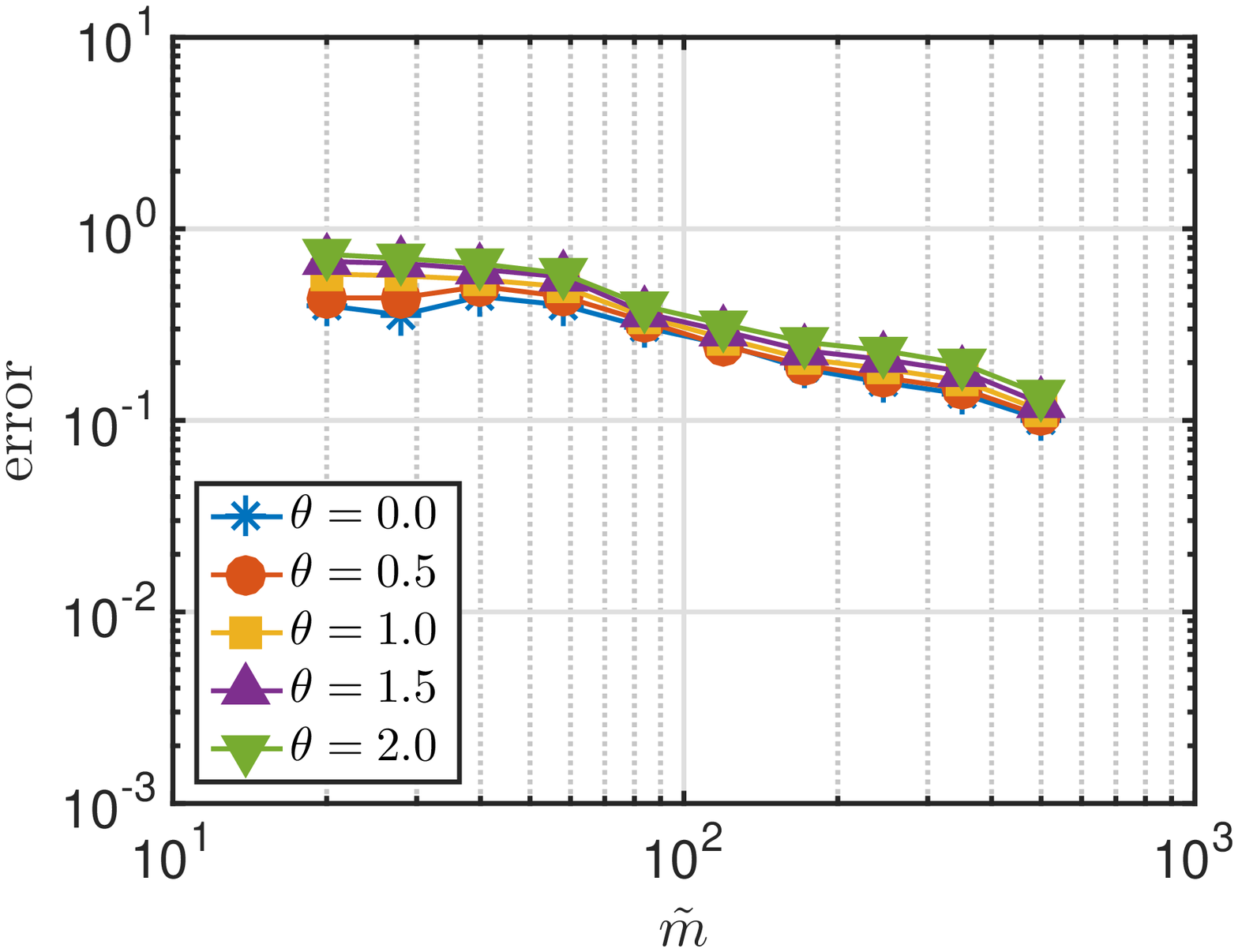}
     \includegraphics[width=1.8in]{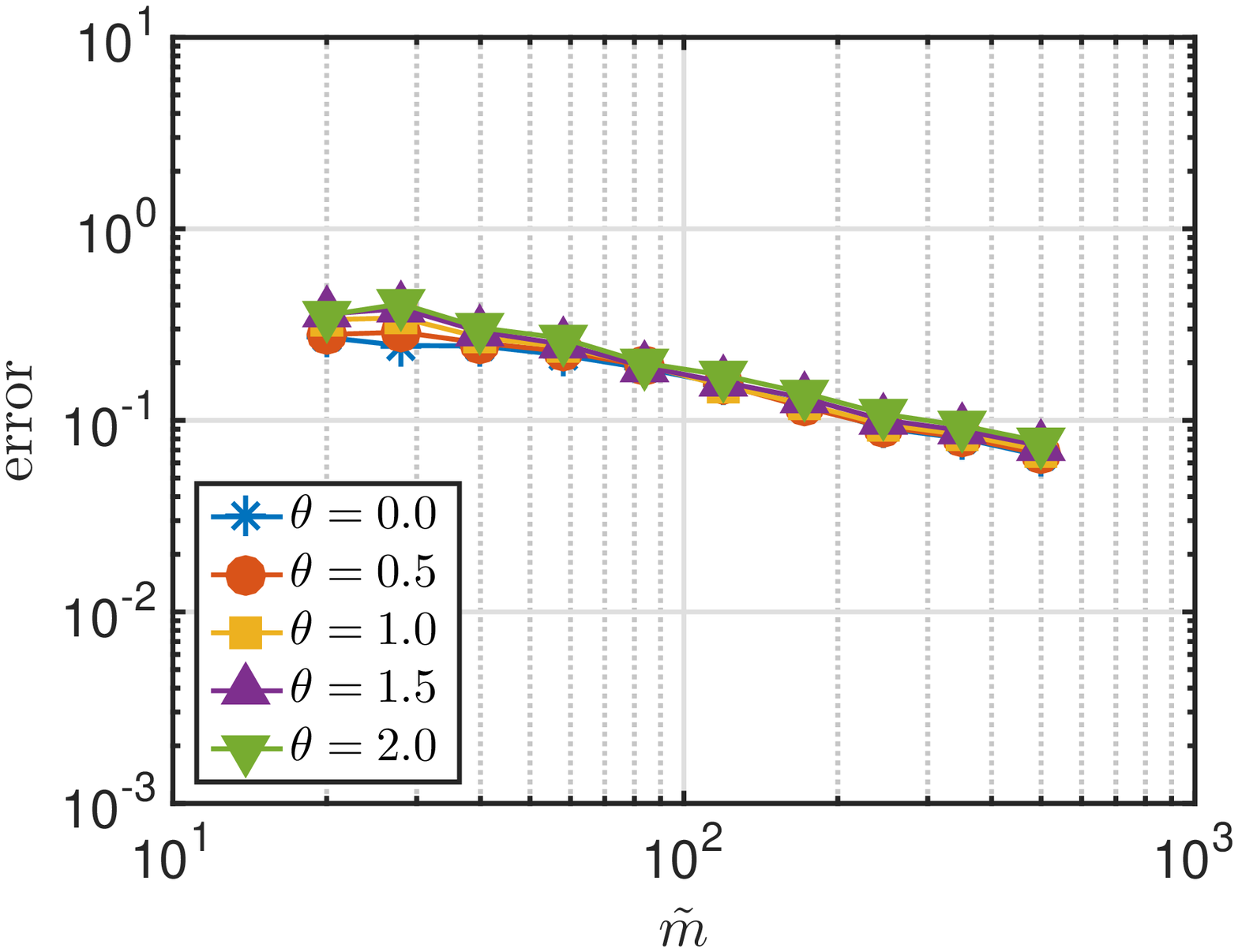}\\
      \includegraphics[width=1.8in]{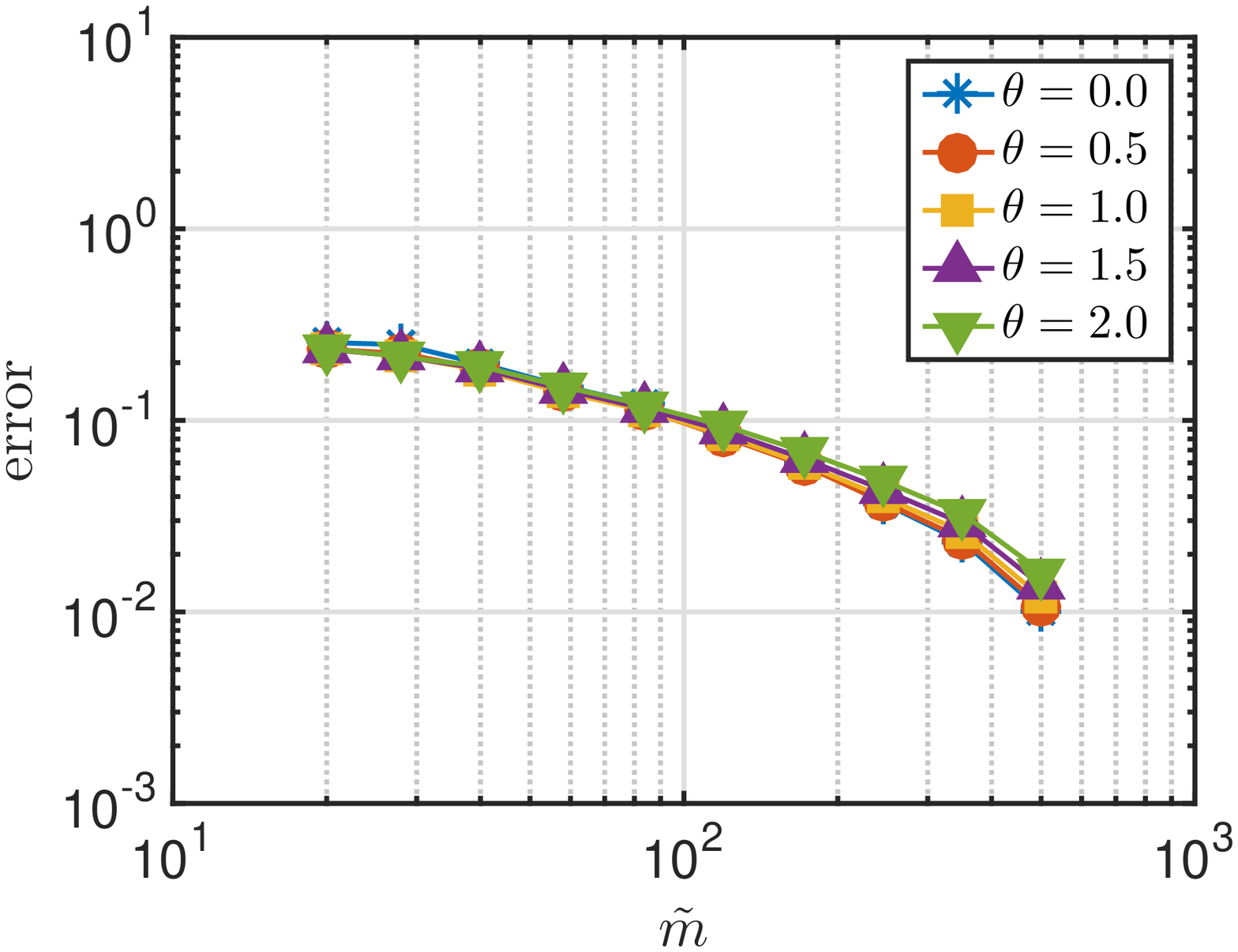} 
     \includegraphics[width=1.8in]{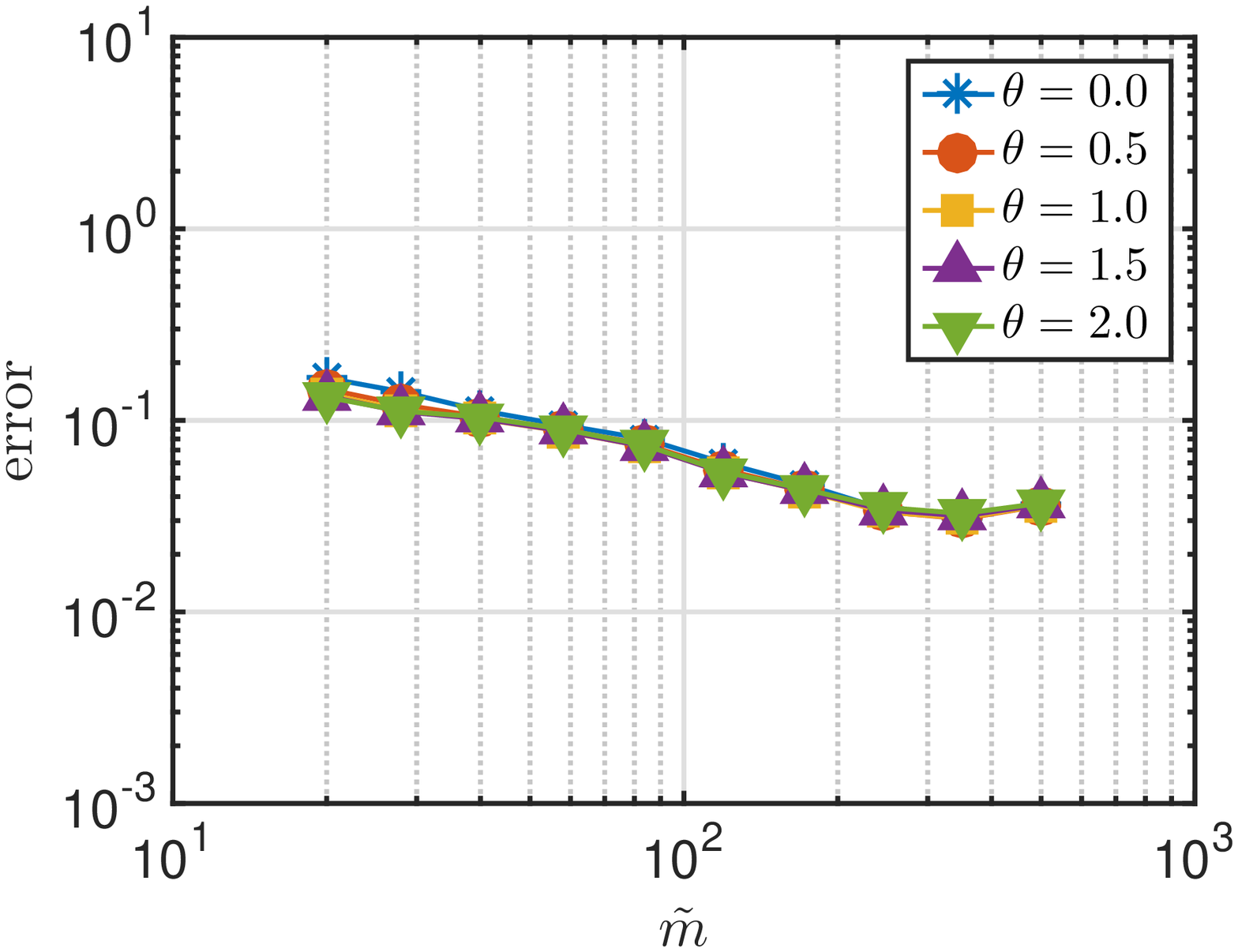}
      \includegraphics[width=1.8in]{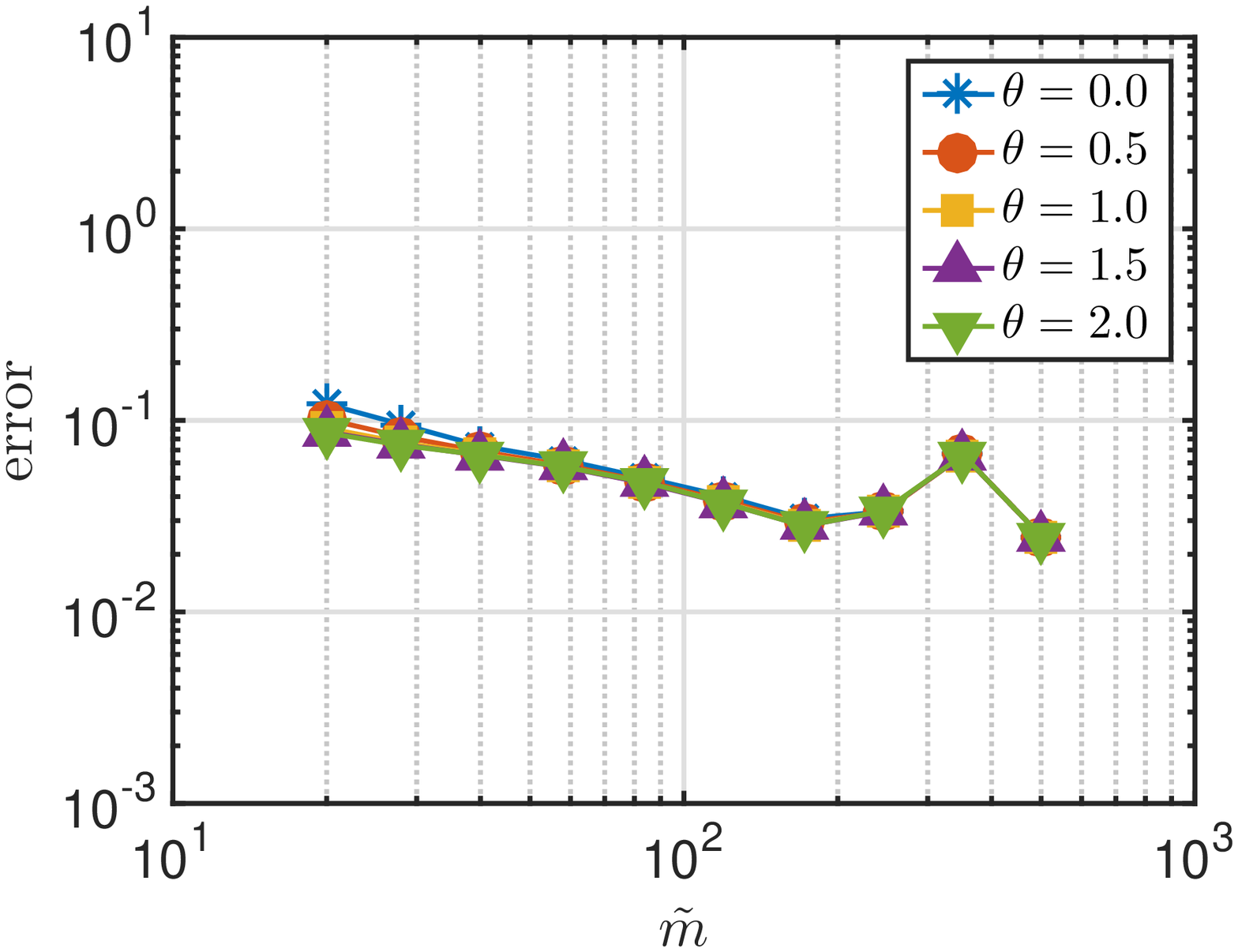}
    \end{tabular}
\caption{The same as Fig.\ \ref{LUoptpeak} but for Chebyshev polynomials with points drawn 
from the Chebyshev density.}
\label{CCoptpeak}
\end{figure}

\begin{figure}[ht]
\centering
  \begin{tabular}{cc}
     \includegraphics[width=1.8in]{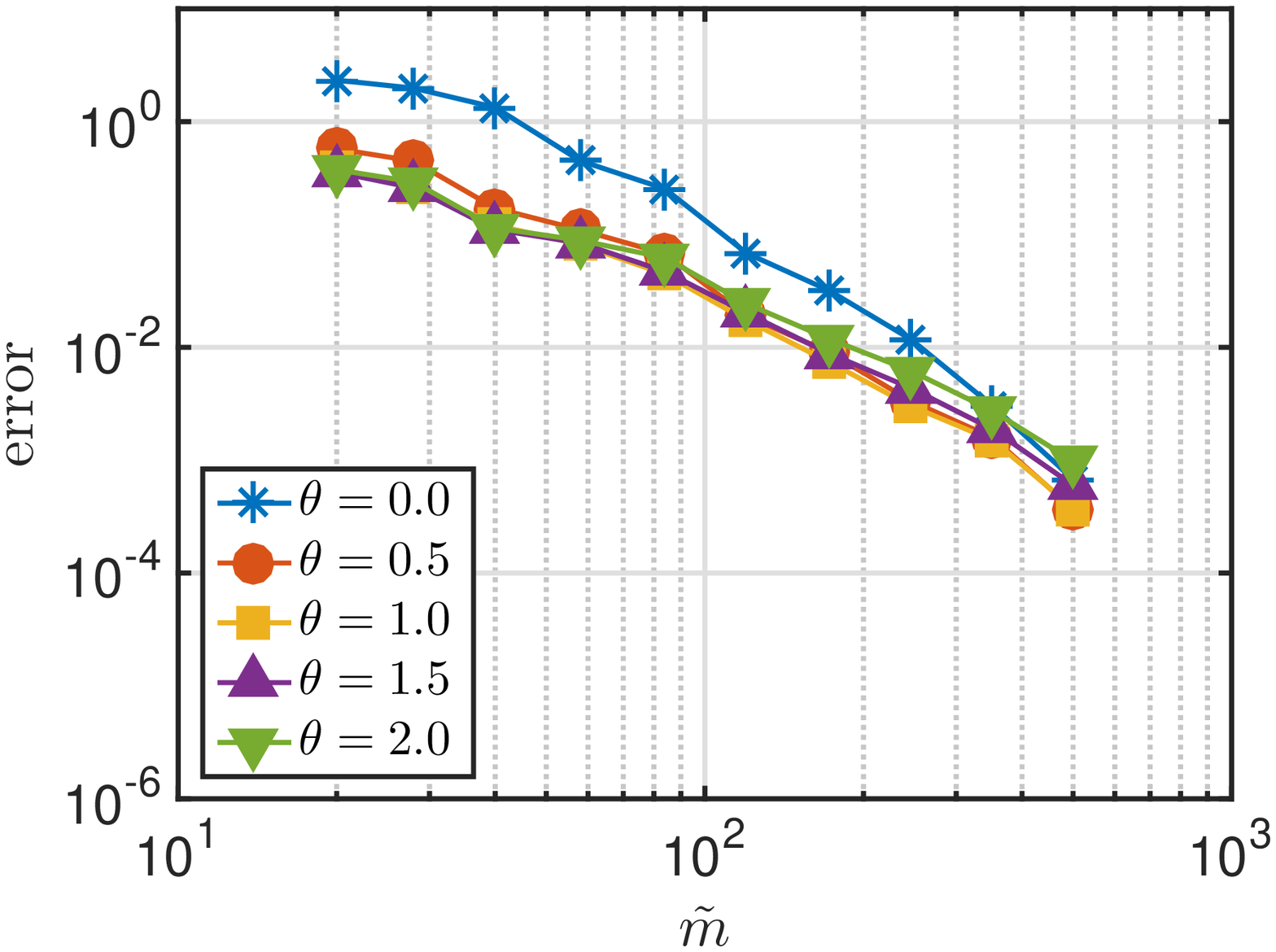} 
     \includegraphics[width=1.8in]{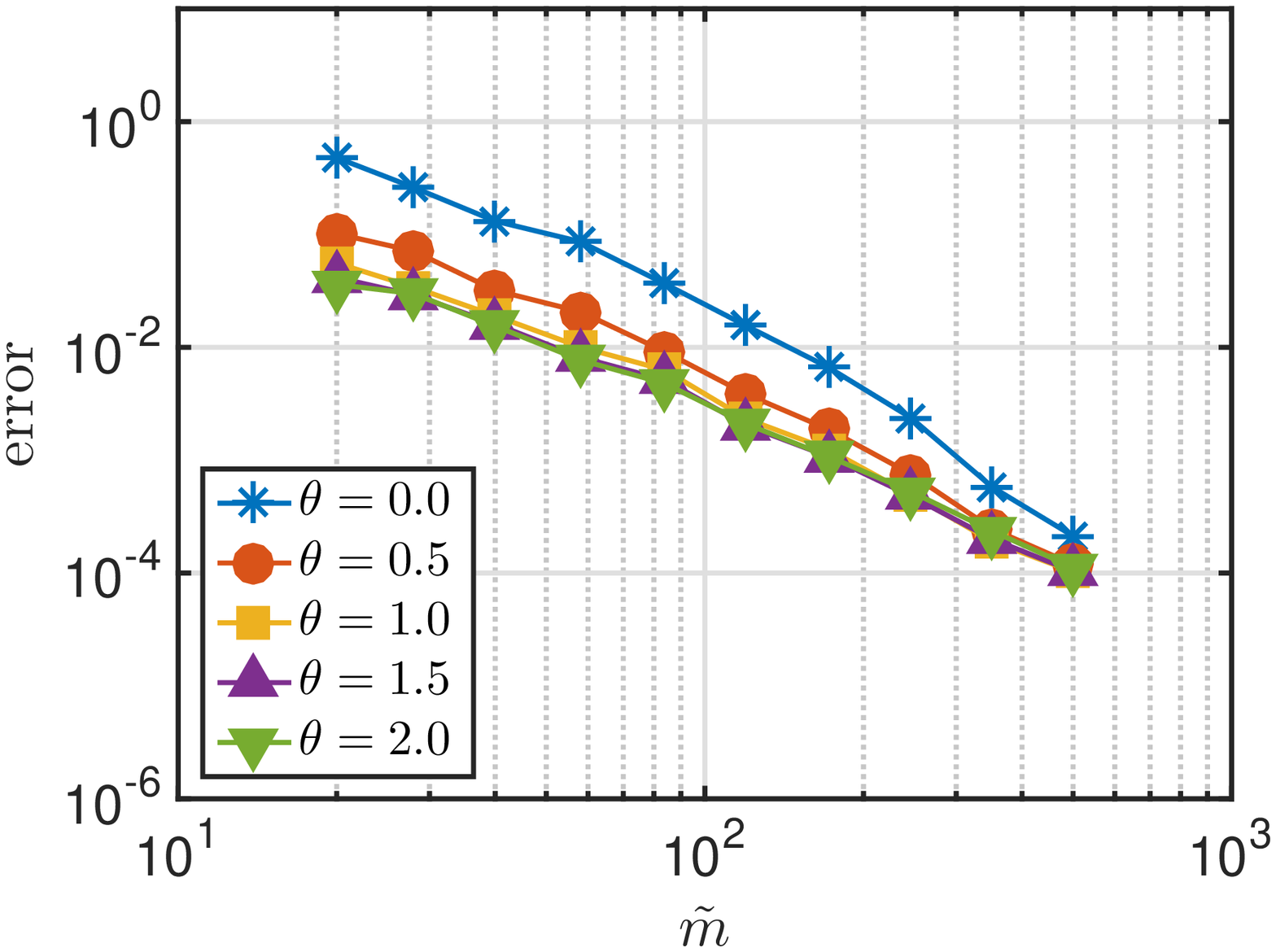}
     \includegraphics[width=1.8in]{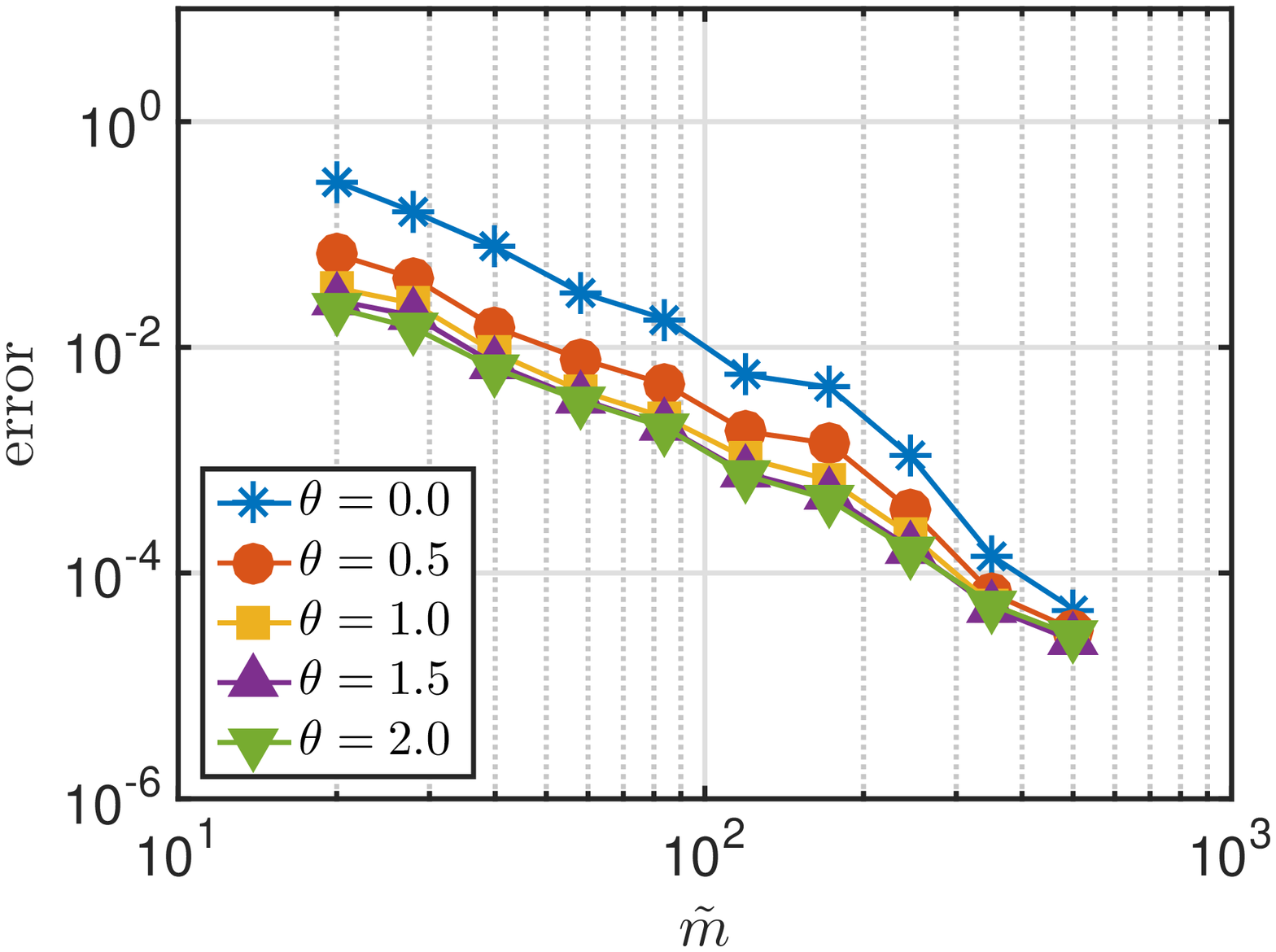}\\
      \includegraphics[width=1.8in]{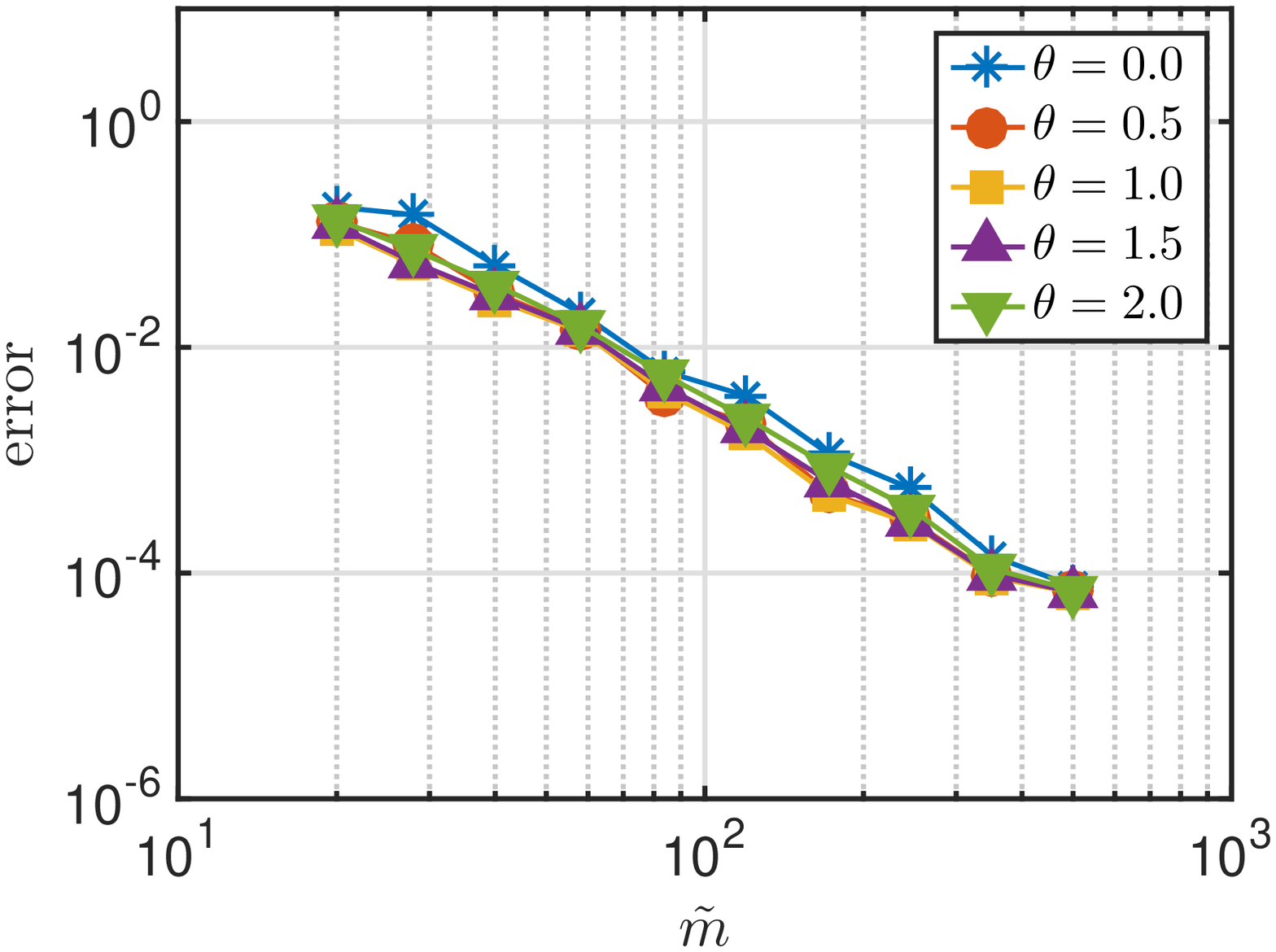} 
     \includegraphics[width=1.8in]{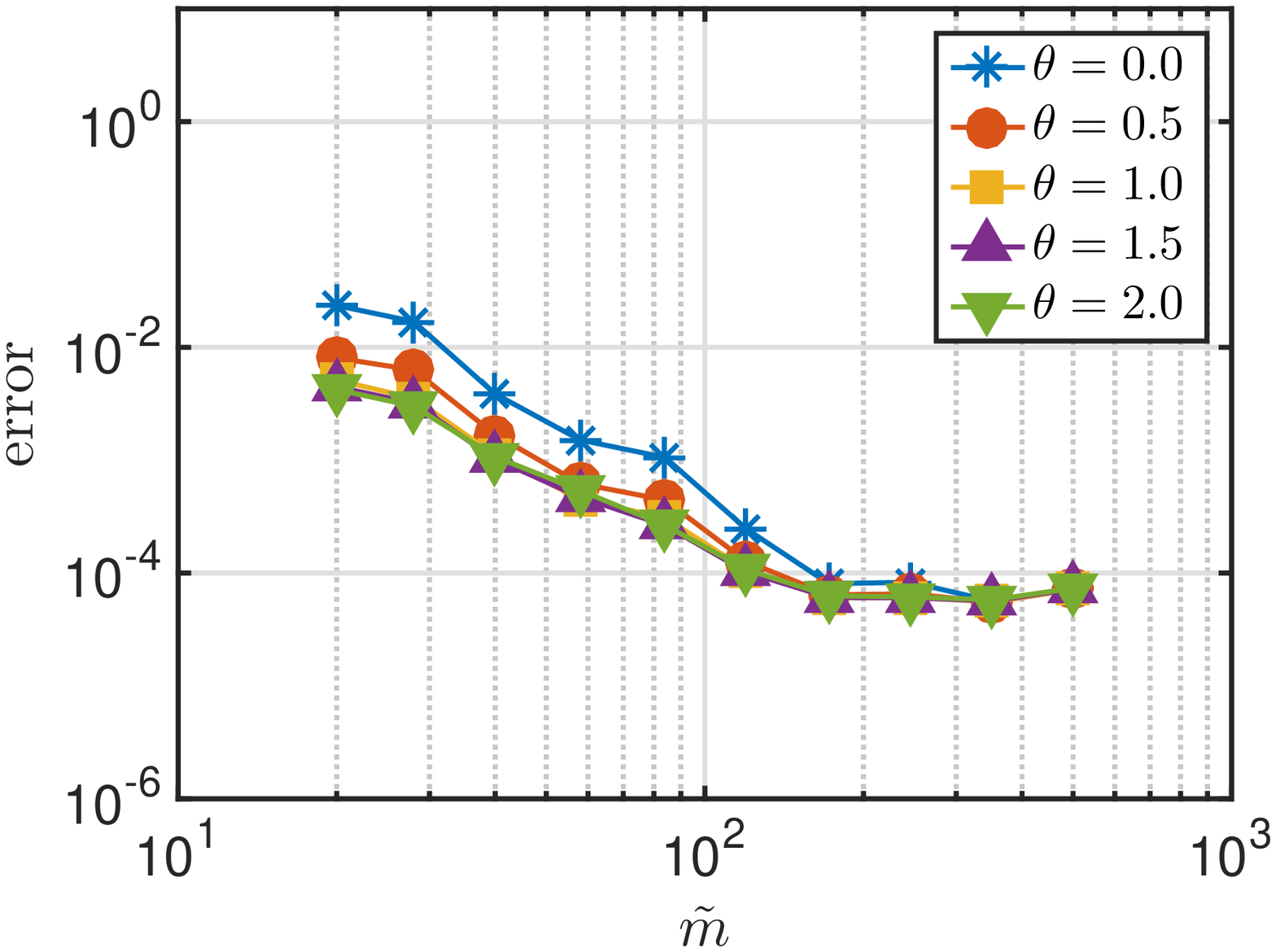}
      \includegraphics[width=1.8in]{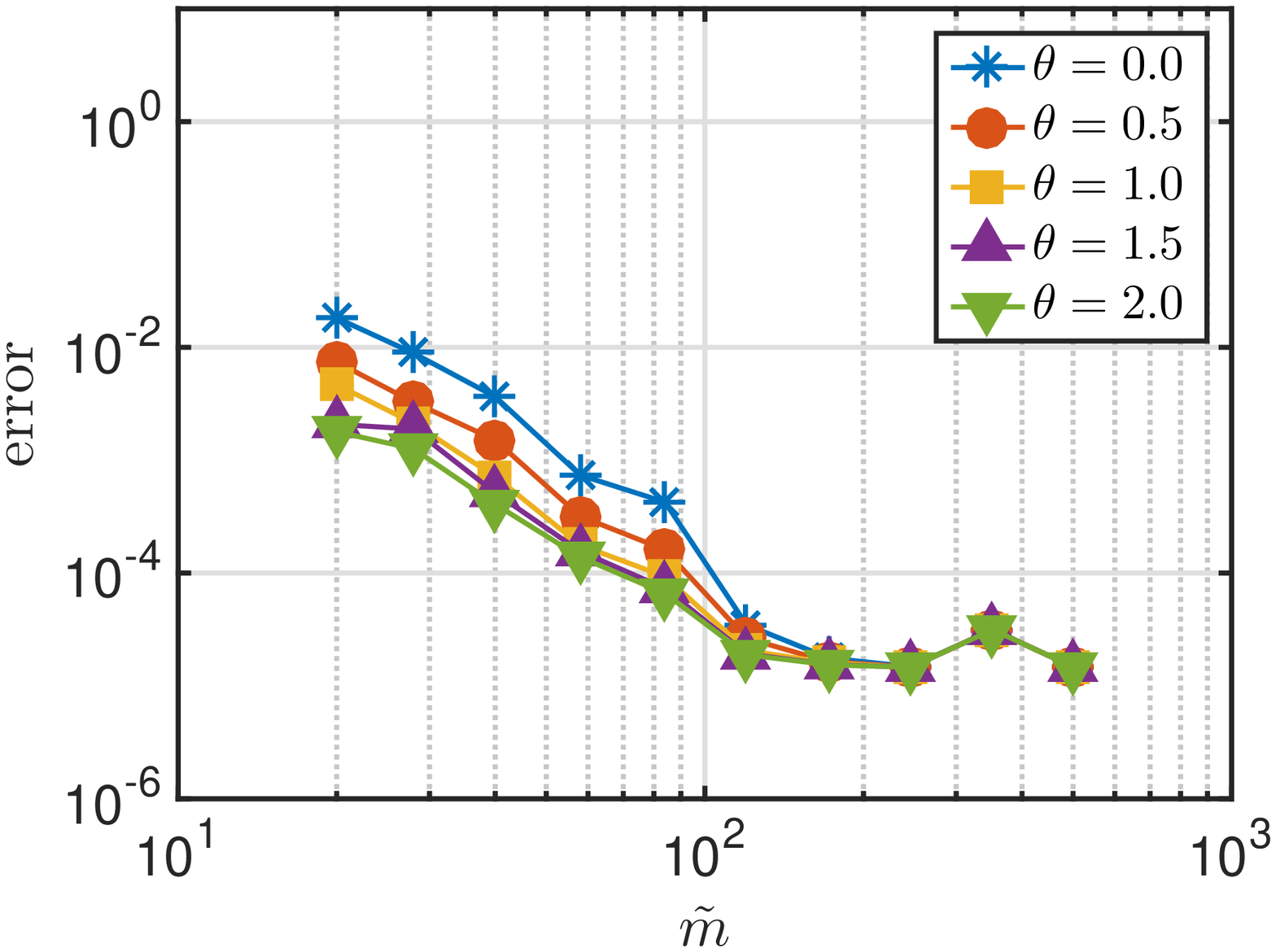}
    \end{tabular}
\caption{The same as Fig.\ \ref{LUoptpeak} but for $f_2$.}
\label{LUoptcos}
\end{figure}

\begin{figure}[h!]
\centering
  \begin{tabular}{cc}
     \includegraphics[width=1.8in]{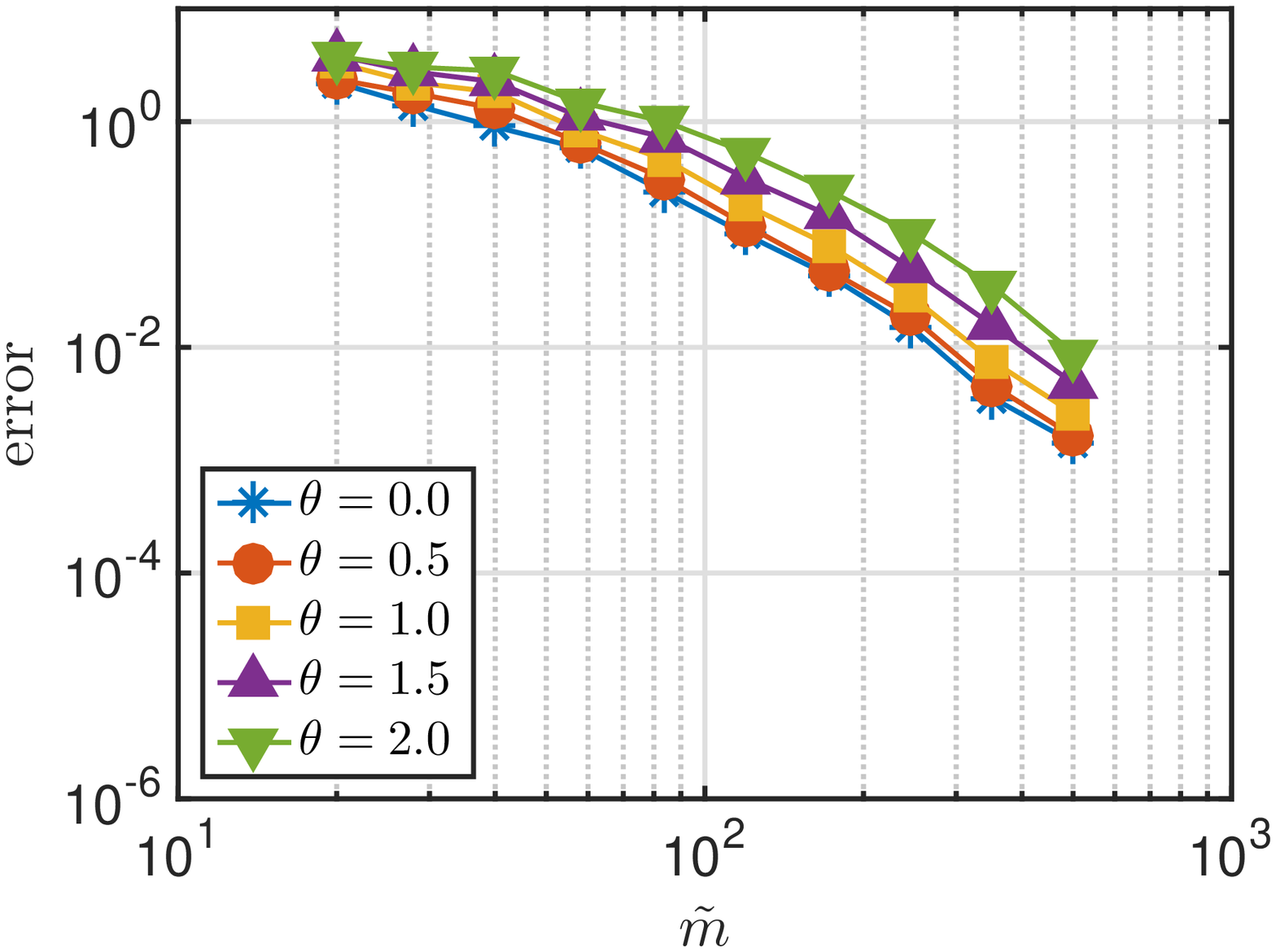} 
     \includegraphics[width=1.8in]{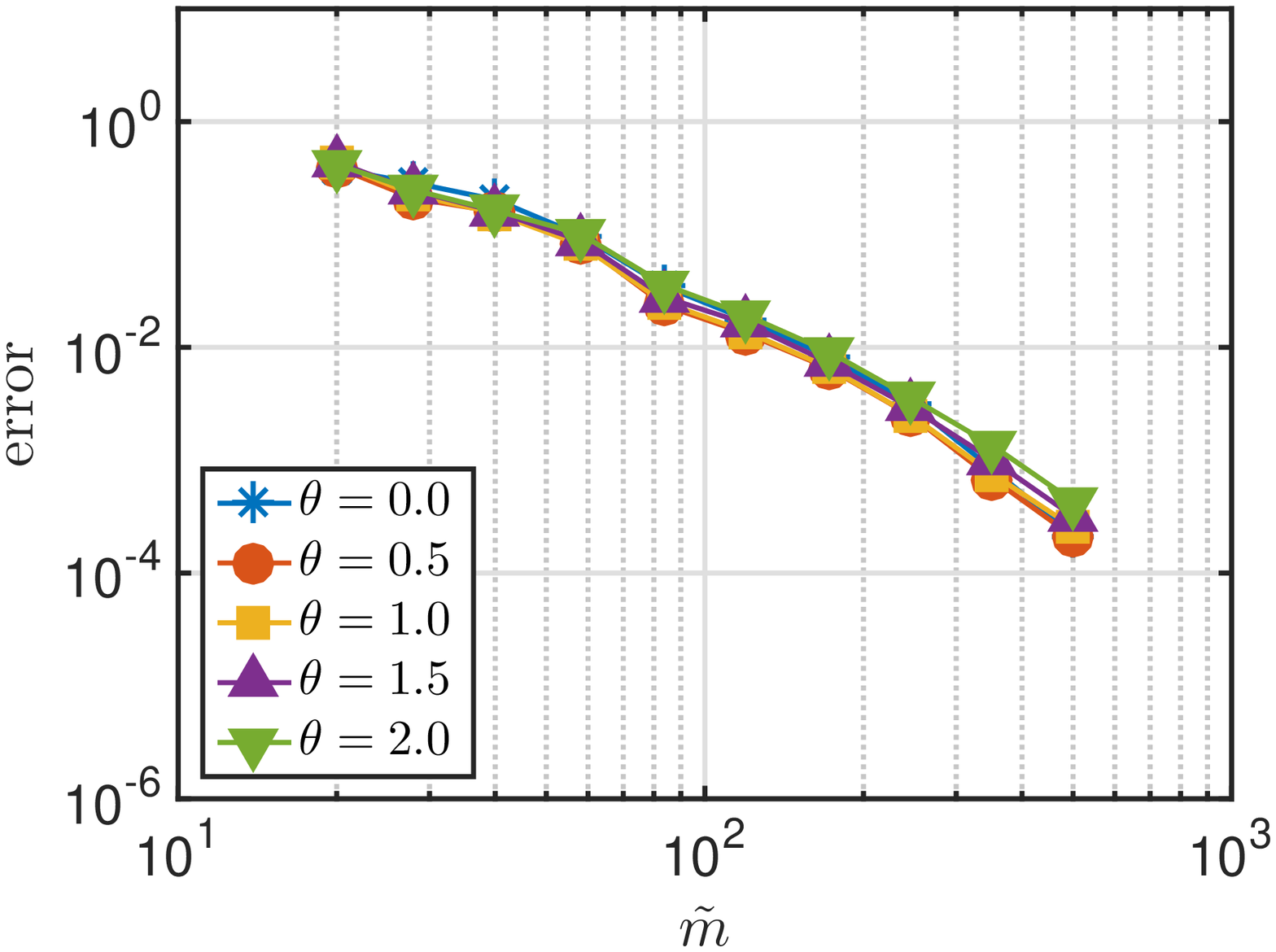}
     \includegraphics[width=1.8in]{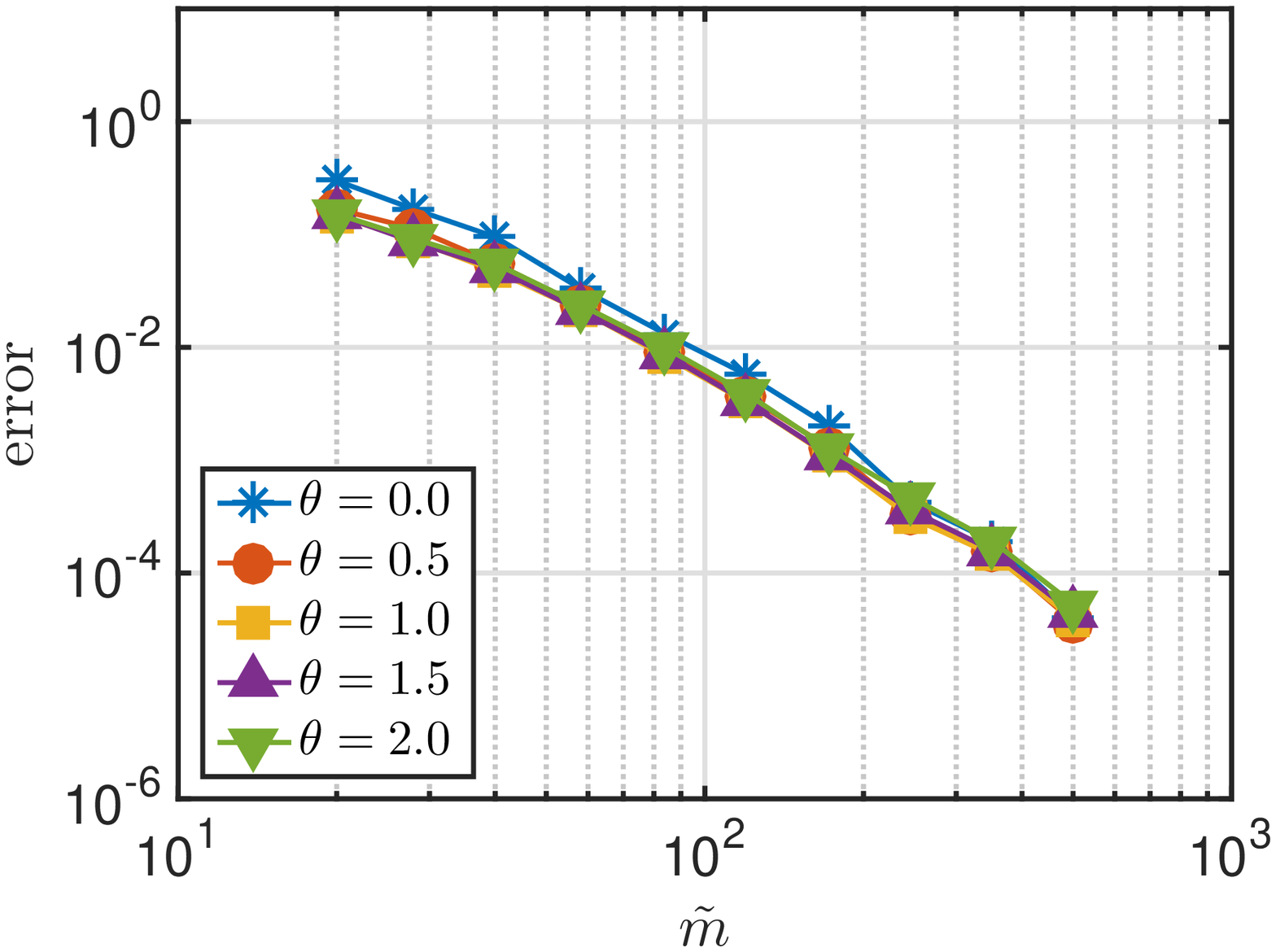}\\
      \includegraphics[width=1.8in]{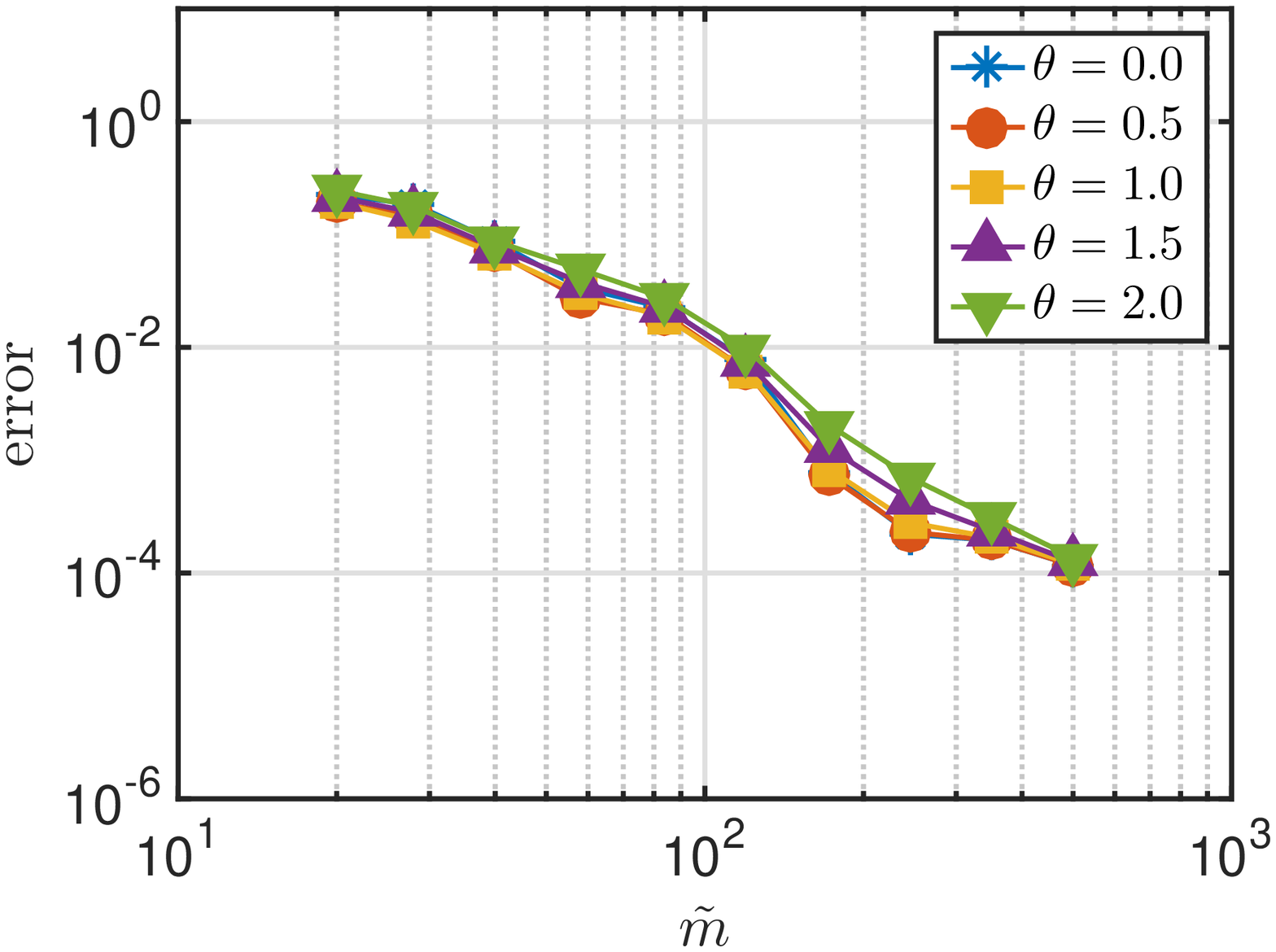} 
     \includegraphics[width=1.8in]{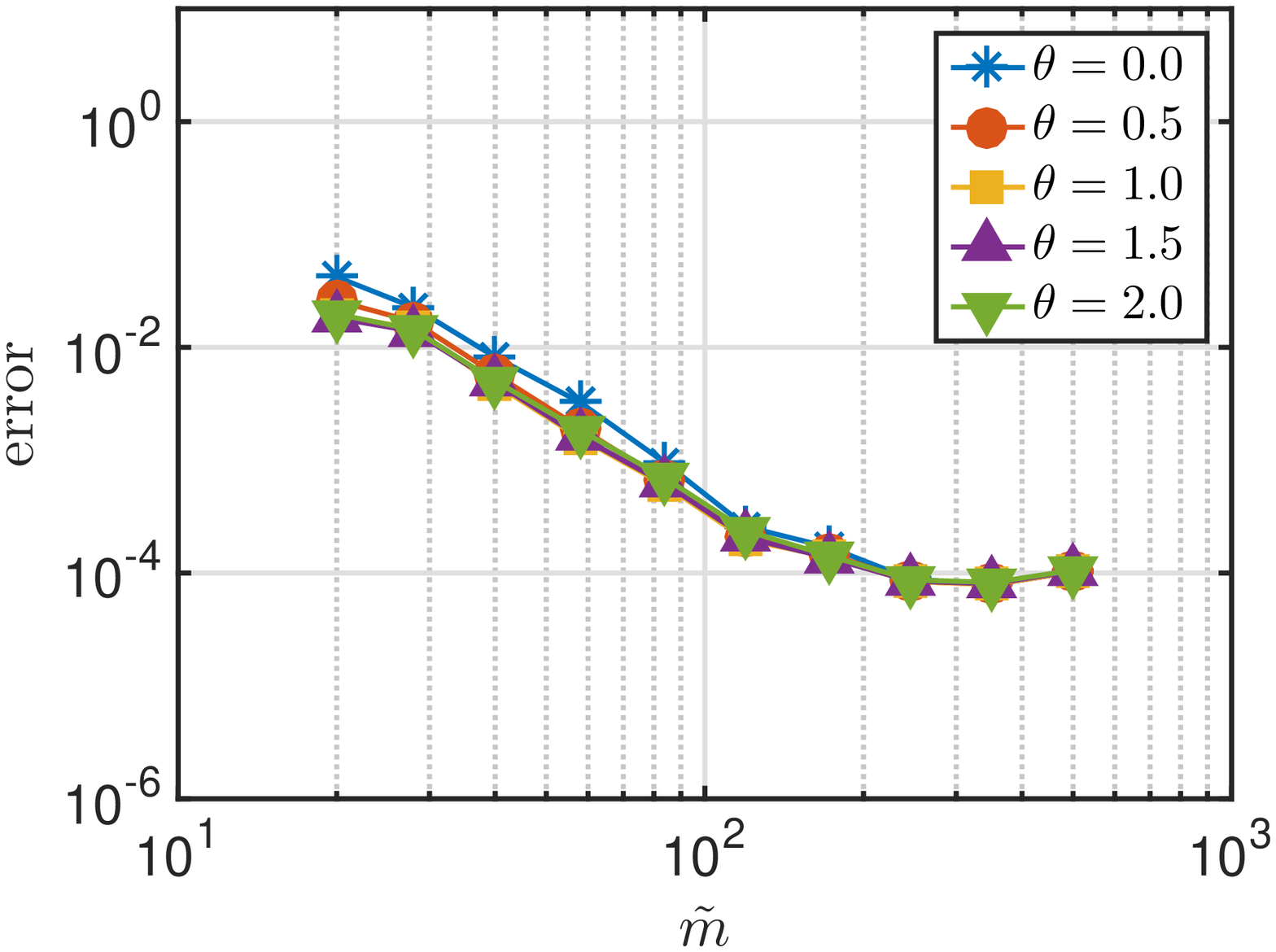}
      \includegraphics[width=1.8in]{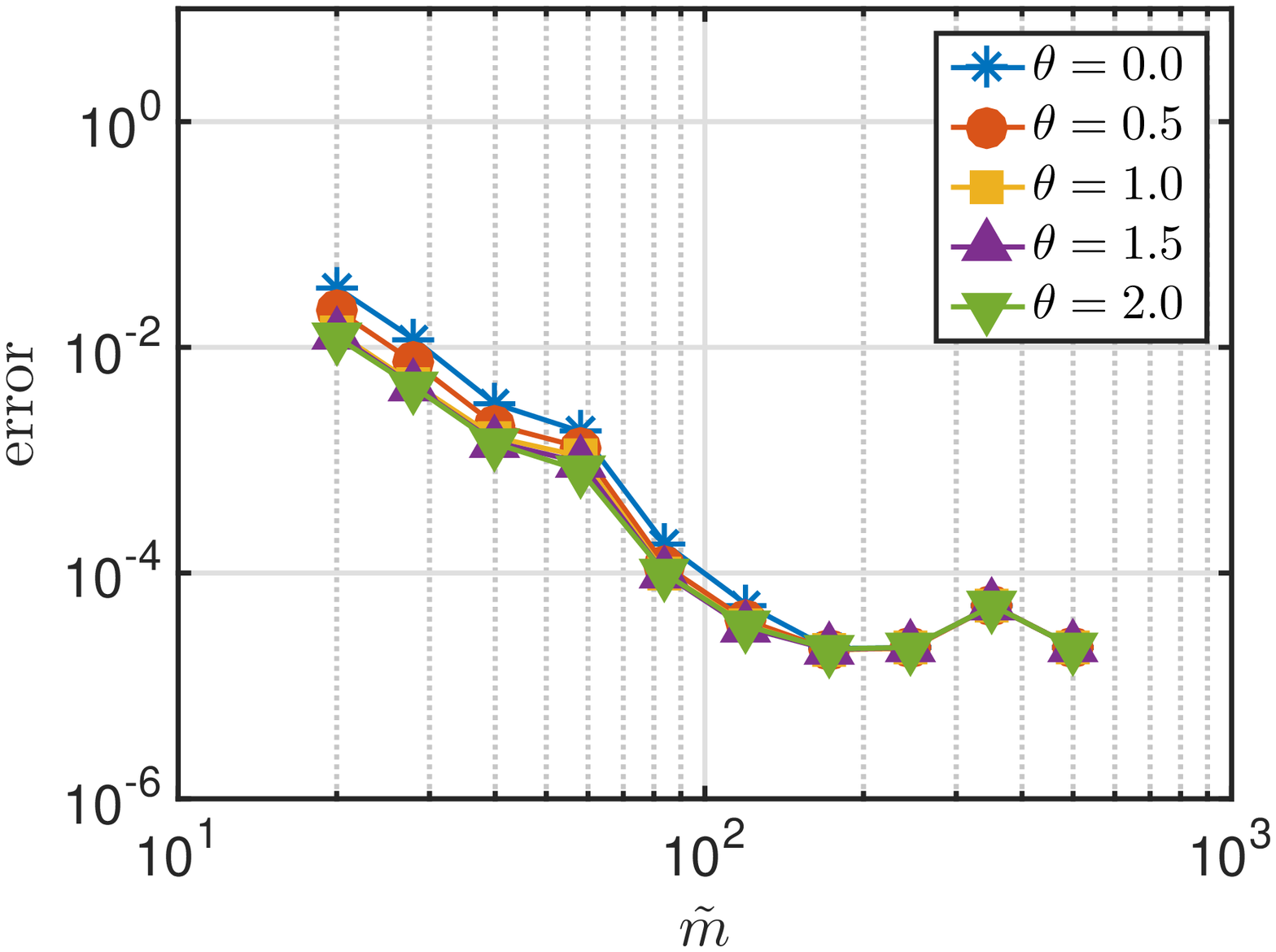}
    \end{tabular}
\caption{The same as Fig.\ \ref{LUoptcos} but for Chebyshev polynomials with points drawn 
from the Chebyshev density.}
\label{CCoptcos}
\end{figure}

\begin{figure}[ht]
\centering
  \begin{tabular}{cc}
     \includegraphics[width=1.8in]{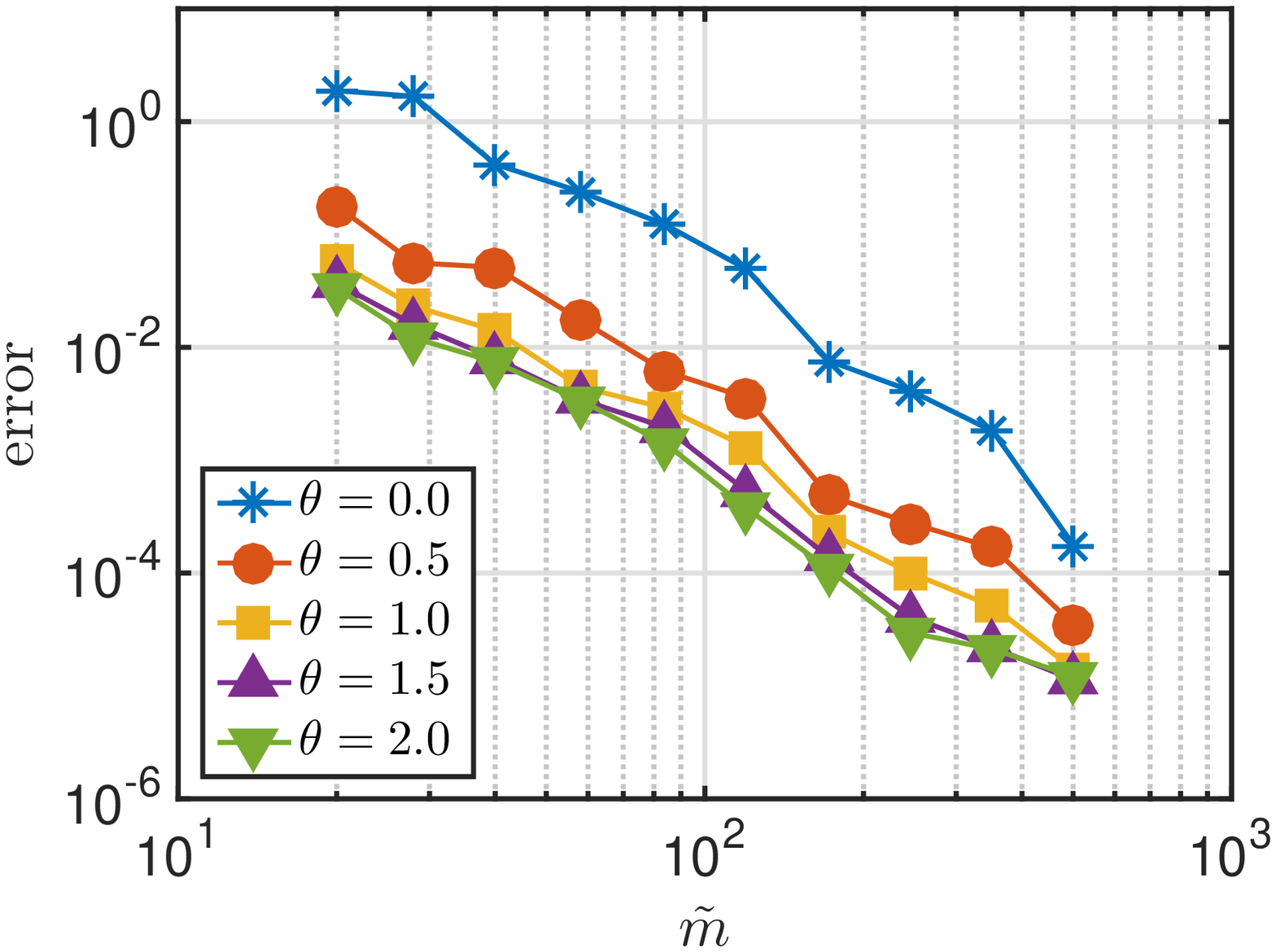} 
     \includegraphics[width=1.8in]{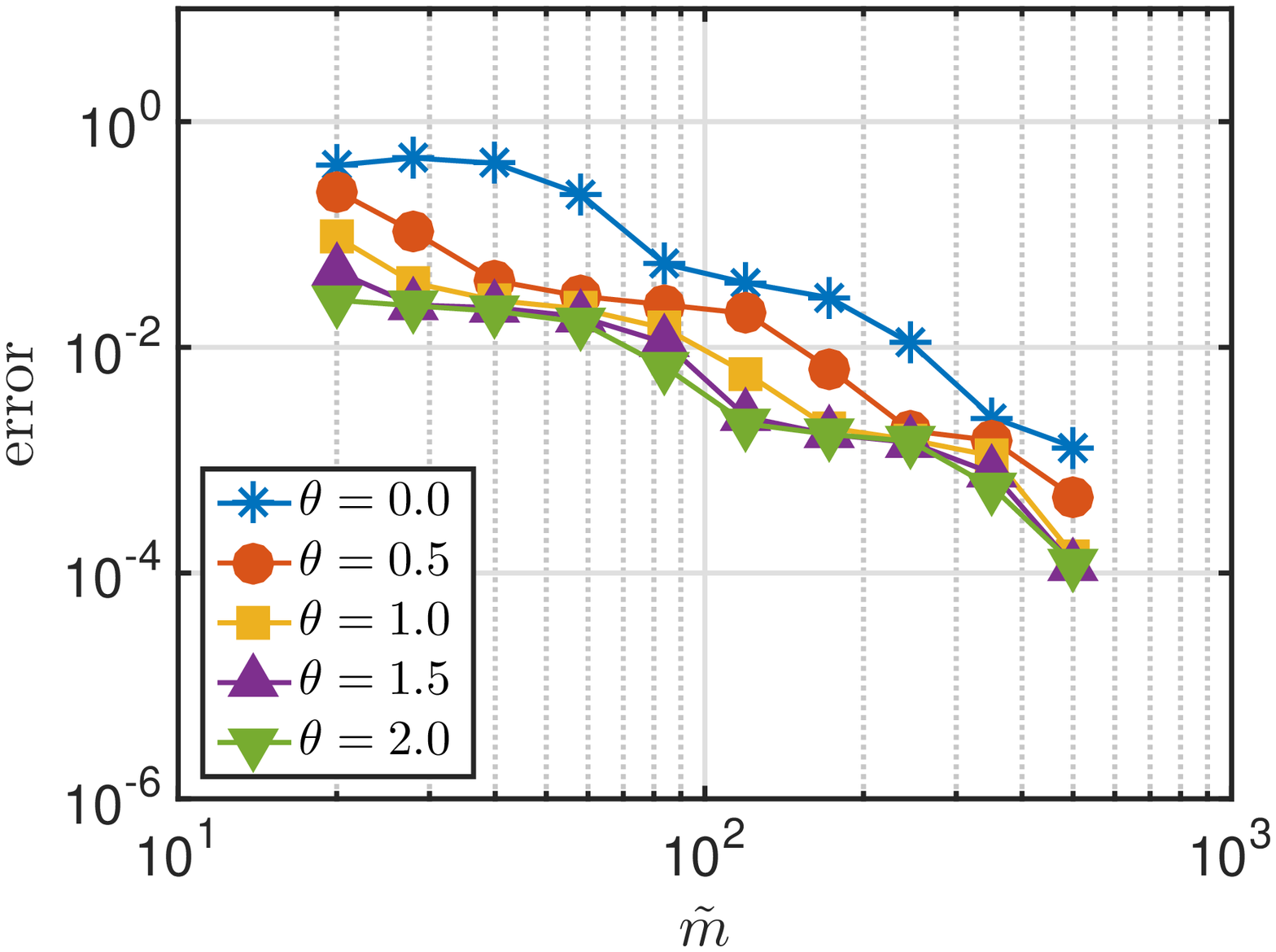}
     \includegraphics[width=1.8in]{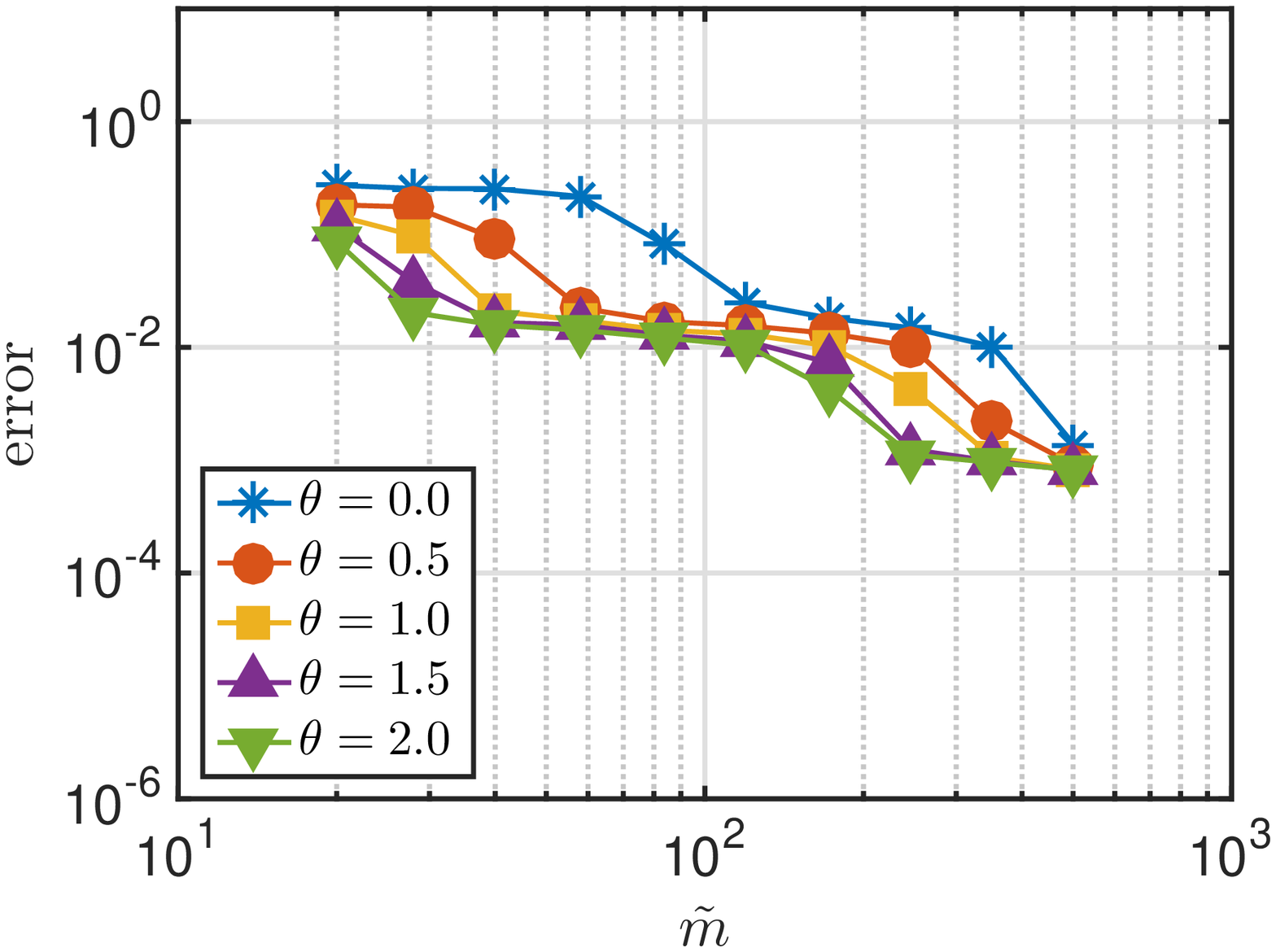}\\
      \includegraphics[width=1.8in]{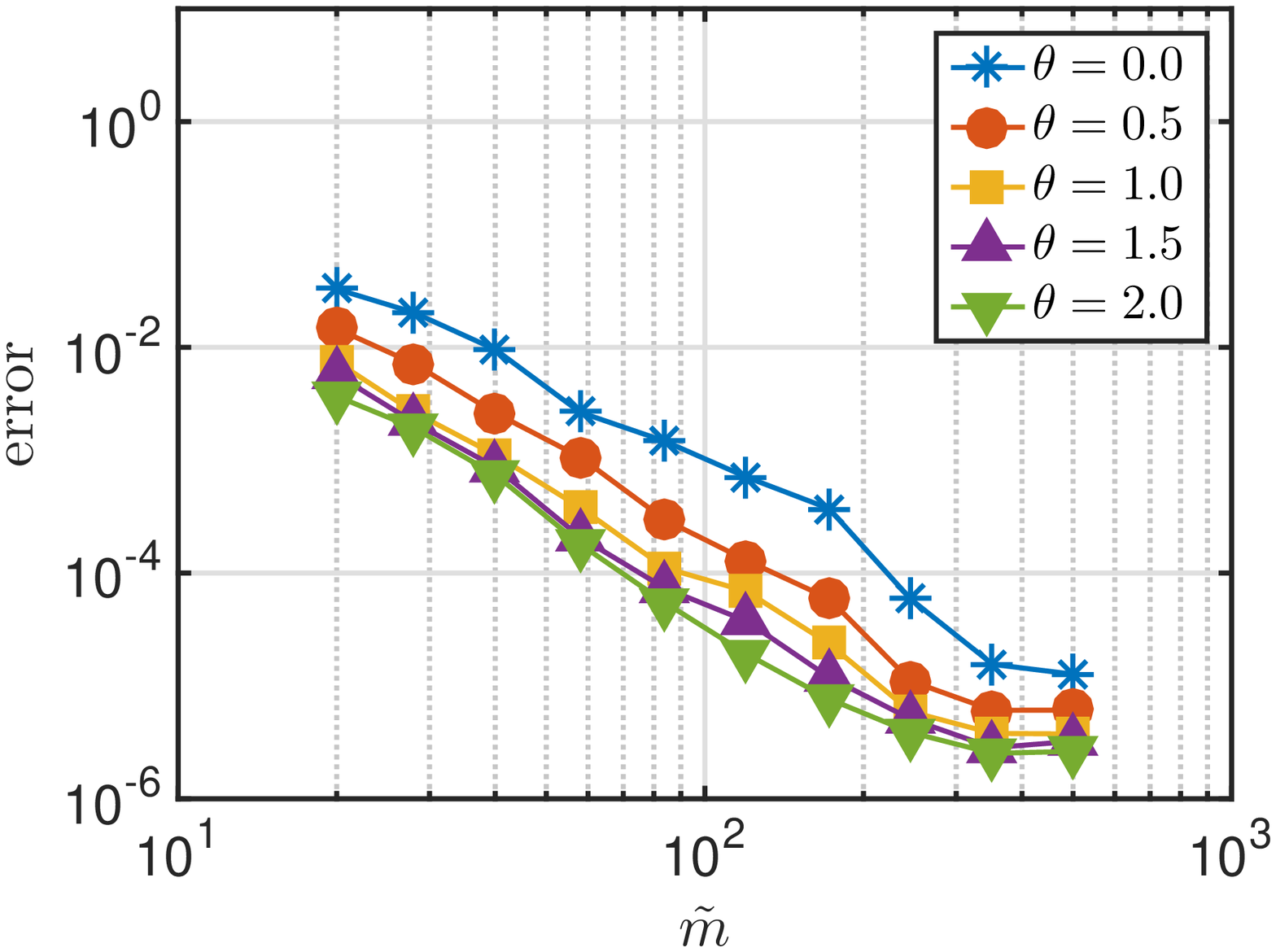} 
     \includegraphics[width=1.8in]{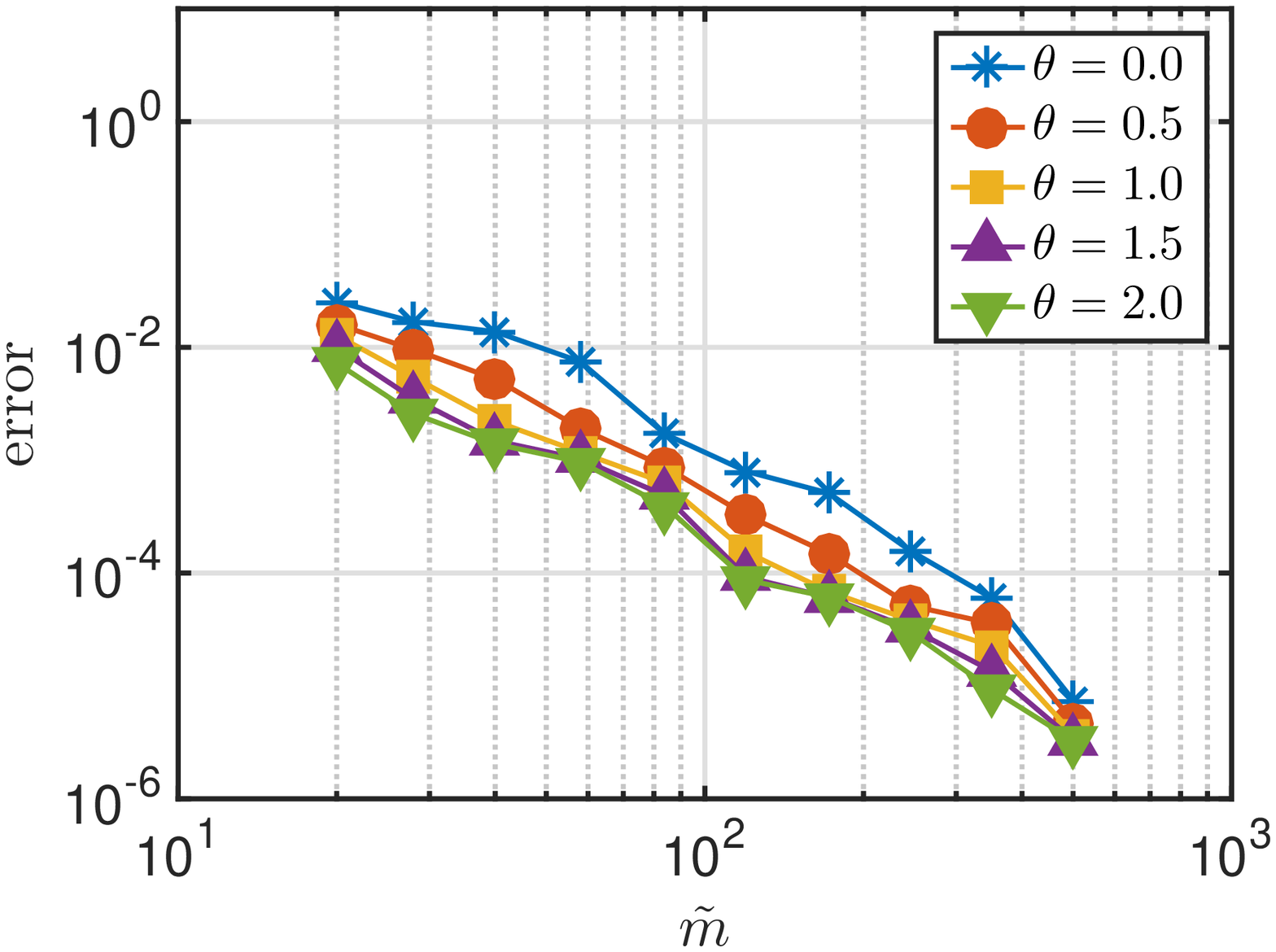}
      \includegraphics[width=1.8in]{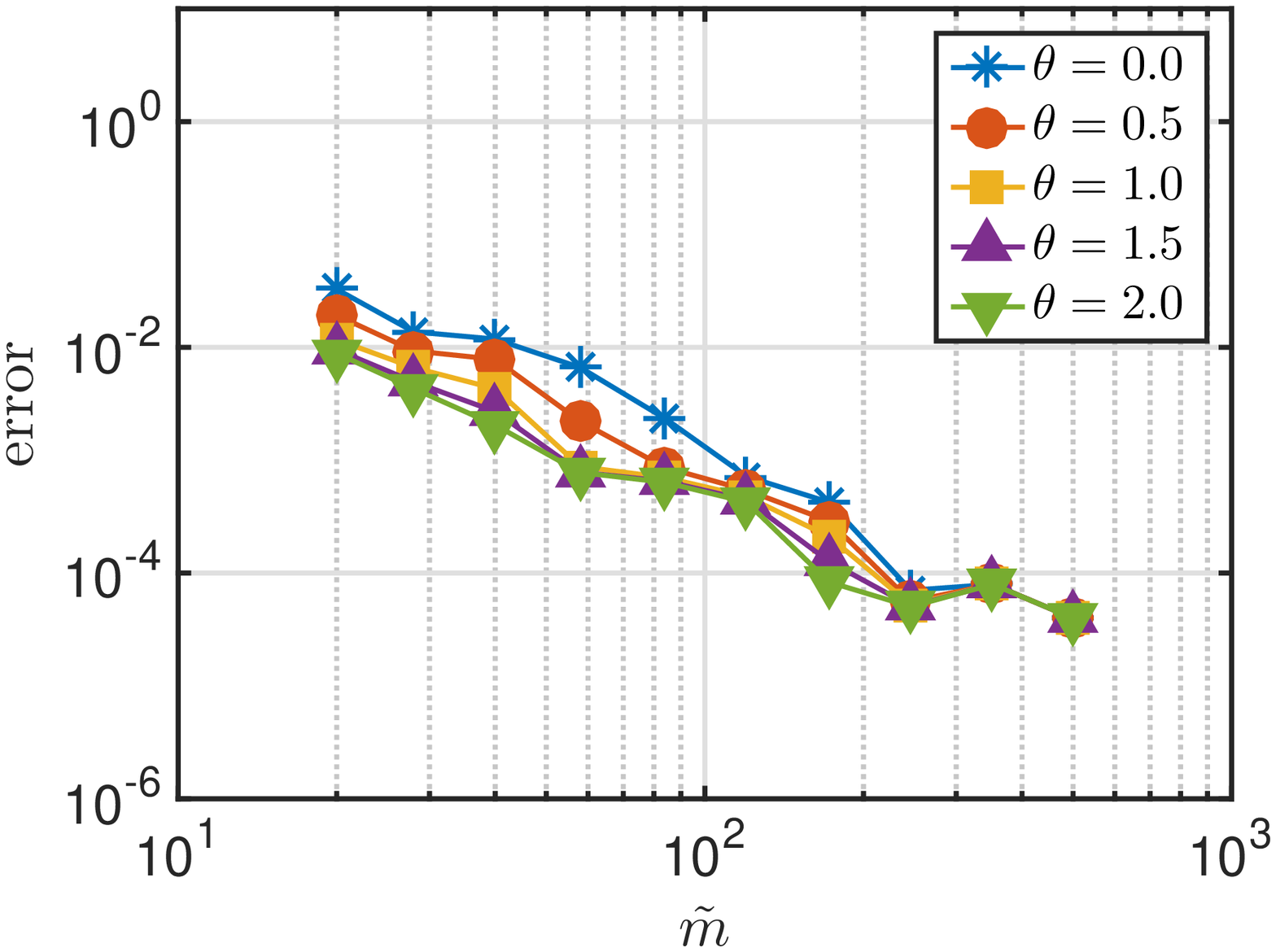}
    \end{tabular}
\caption{The same as Fig.\ \ref{LUoptpeak} but for $f_3$.}
\label{LUoptexp}
\end{figure}

\begin{figure}[h!]
\centering
  \begin{tabular}{cc}
     \includegraphics[width=1.8in]{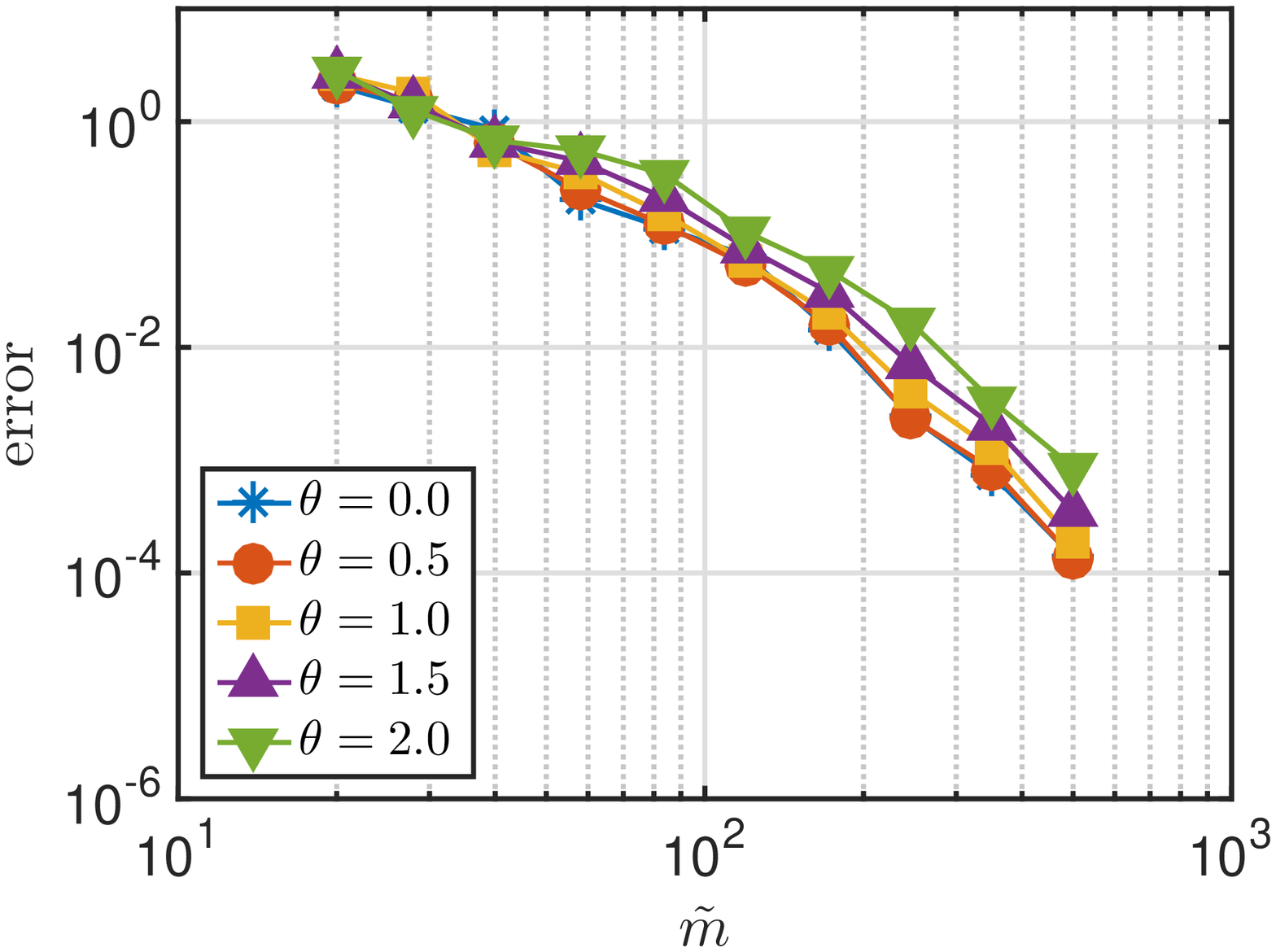} 
     \includegraphics[width=1.8in]{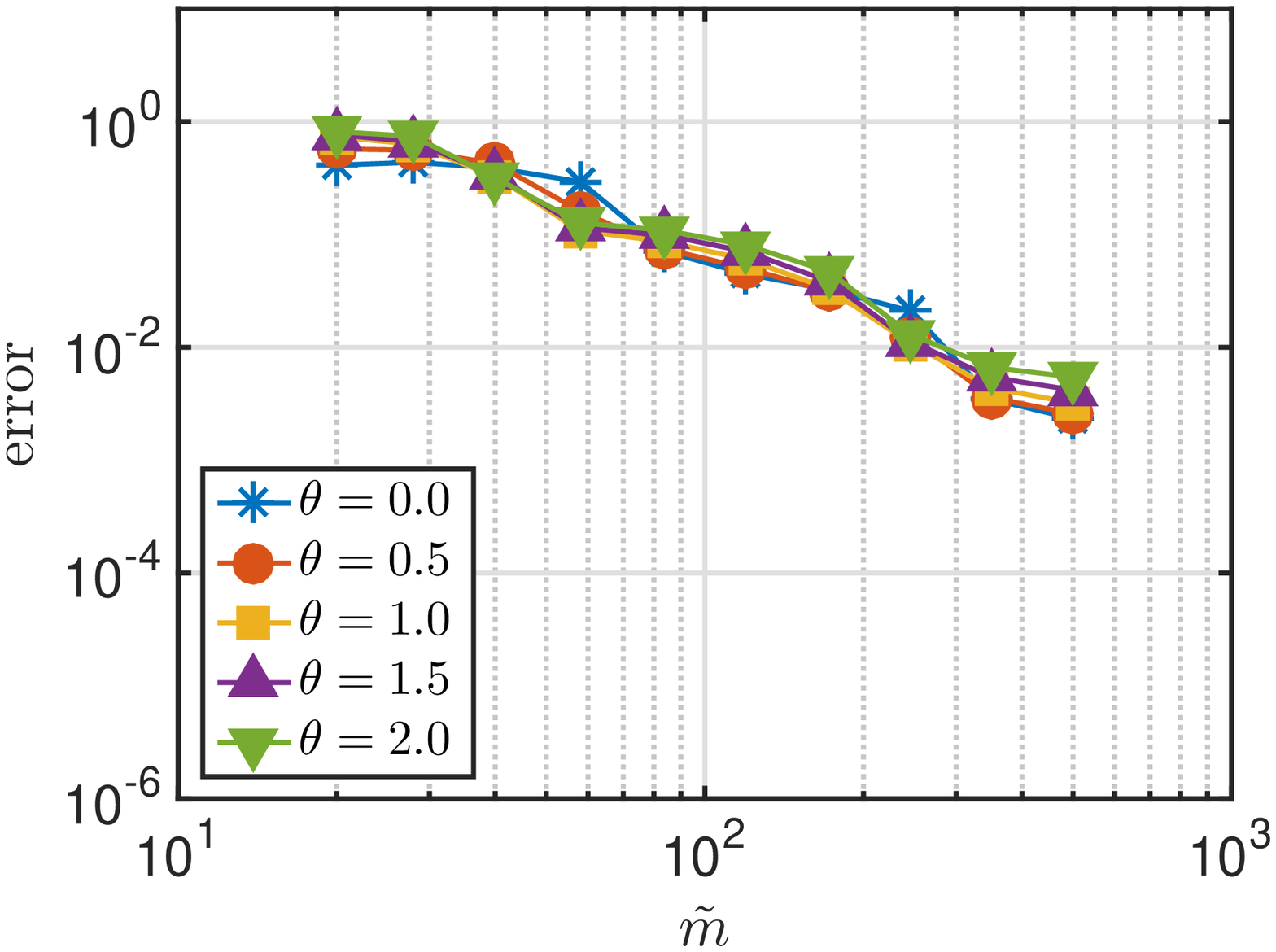}
     \includegraphics[width=1.8in]{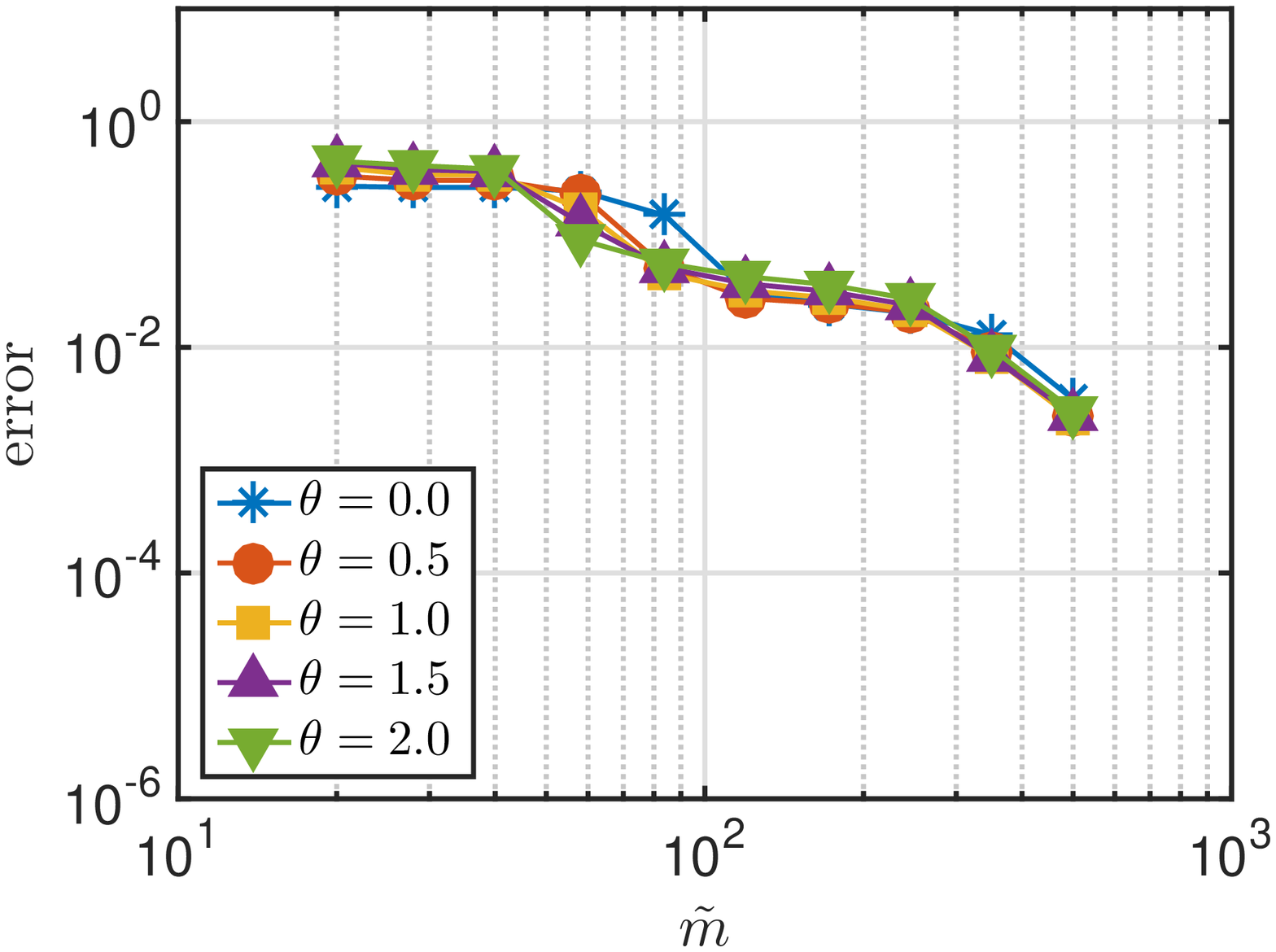}\\
      \includegraphics[width=1.8in]{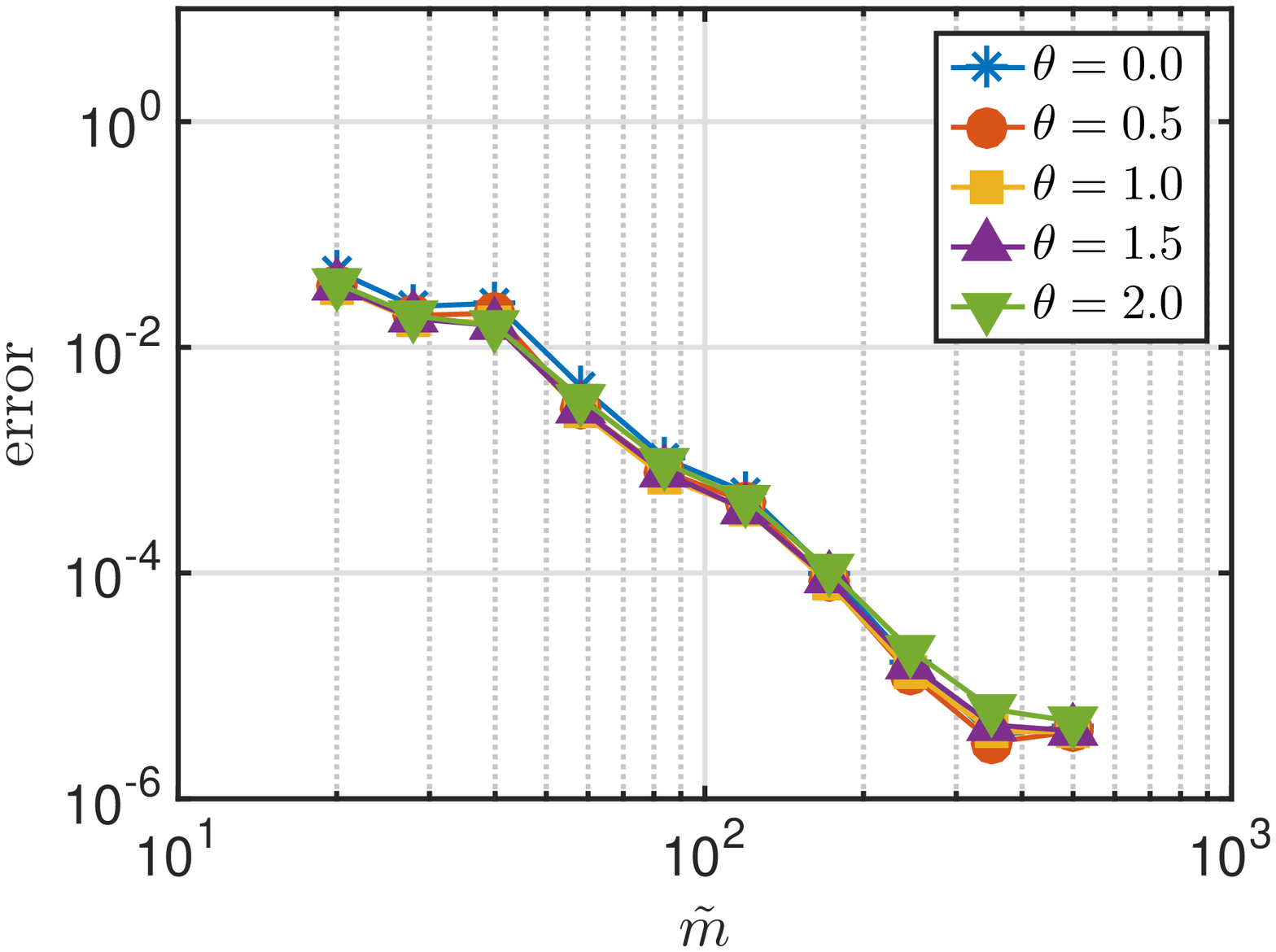} 
     \includegraphics[width=1.8in]{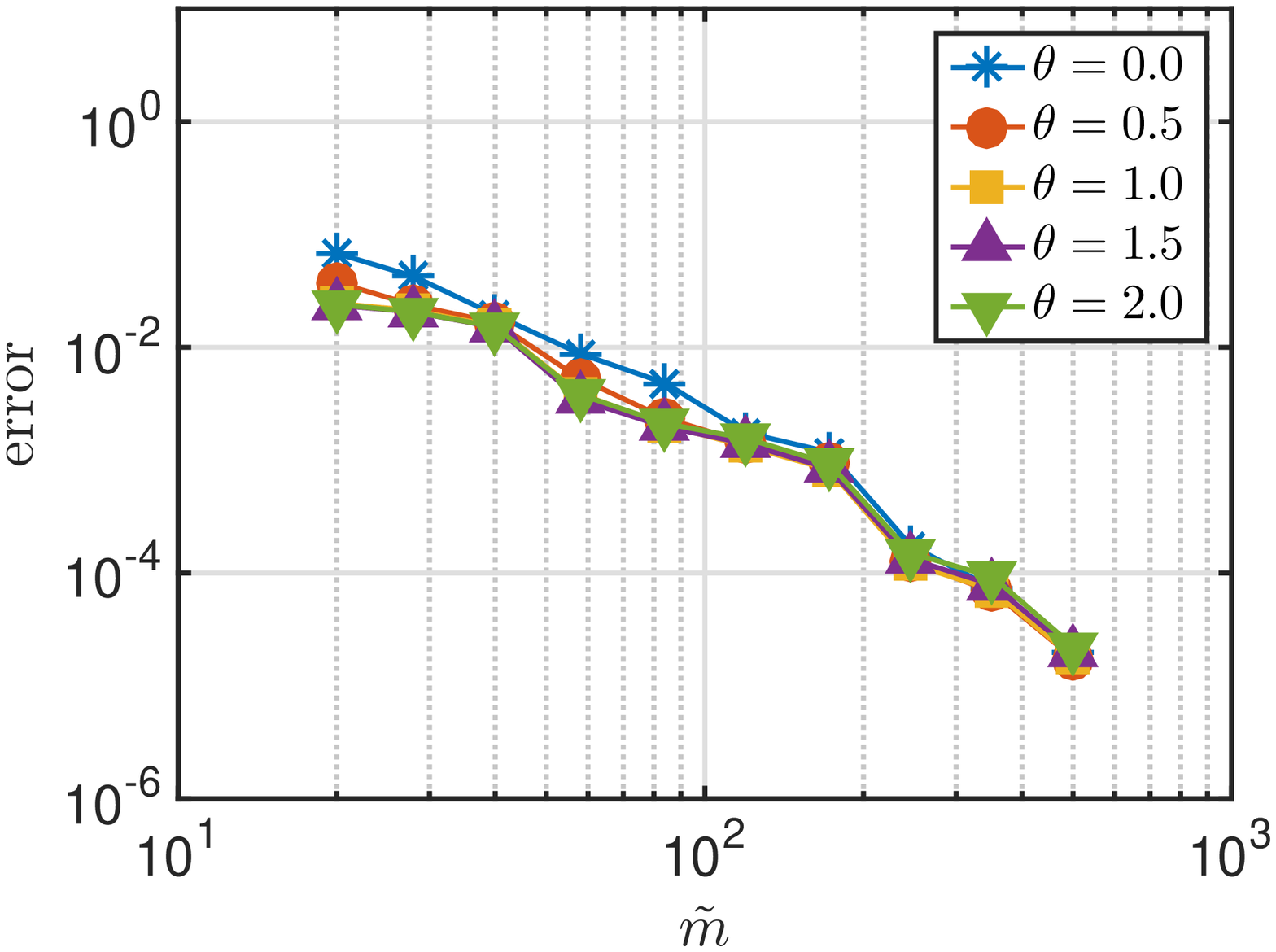}
      \includegraphics[width=1.8in]{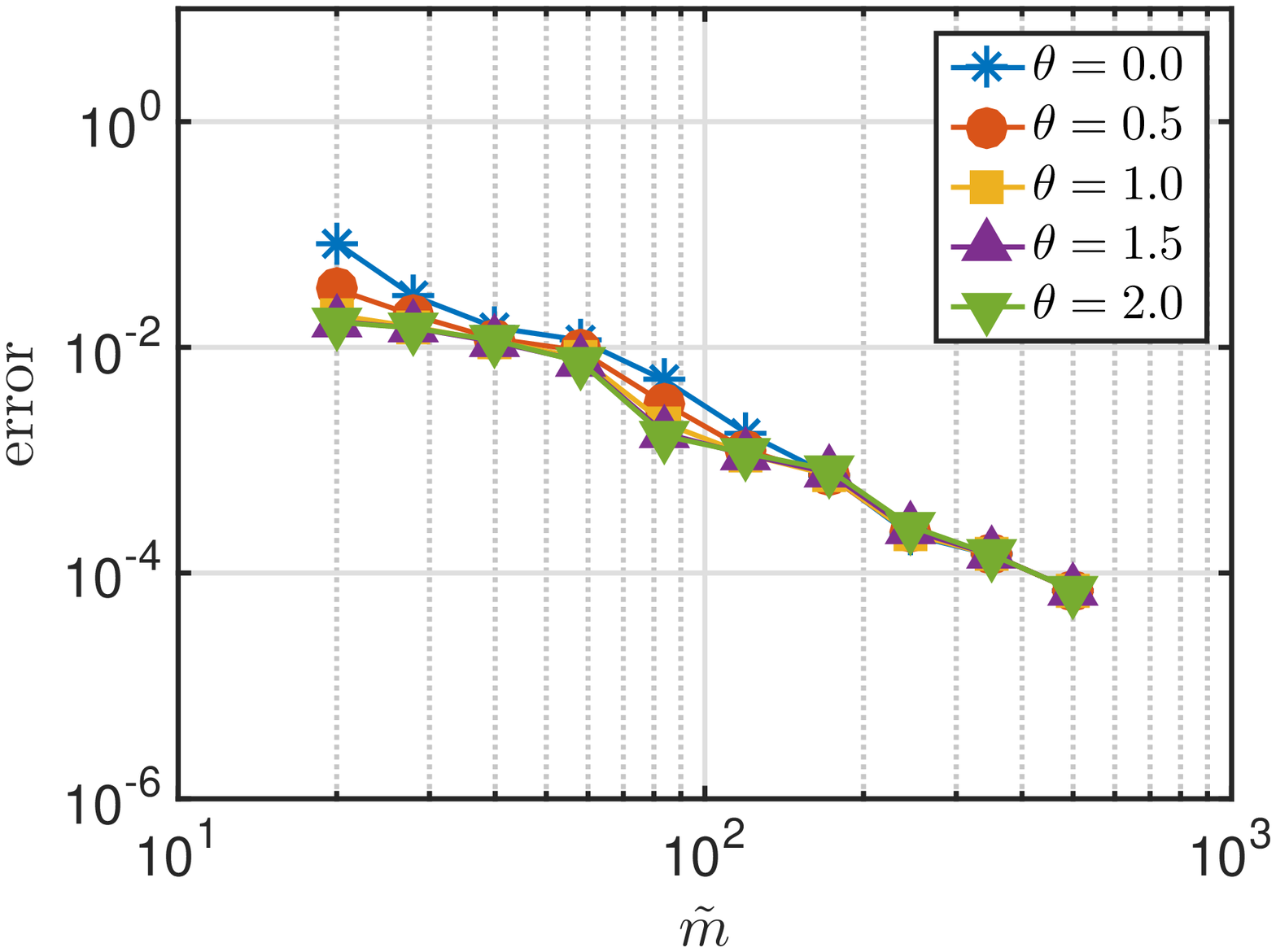}
    \end{tabular}
\caption{The same as Fig.\ \ref{LUoptexp} but for Chebyshev polynomials with points drawn 
from the Chebyshev density.}
\label{CCoptexp}
\end{figure}

\begin{figure}[ht]
\centering
  \begin{tabular}{cc}
     \includegraphics[width=2.1in]{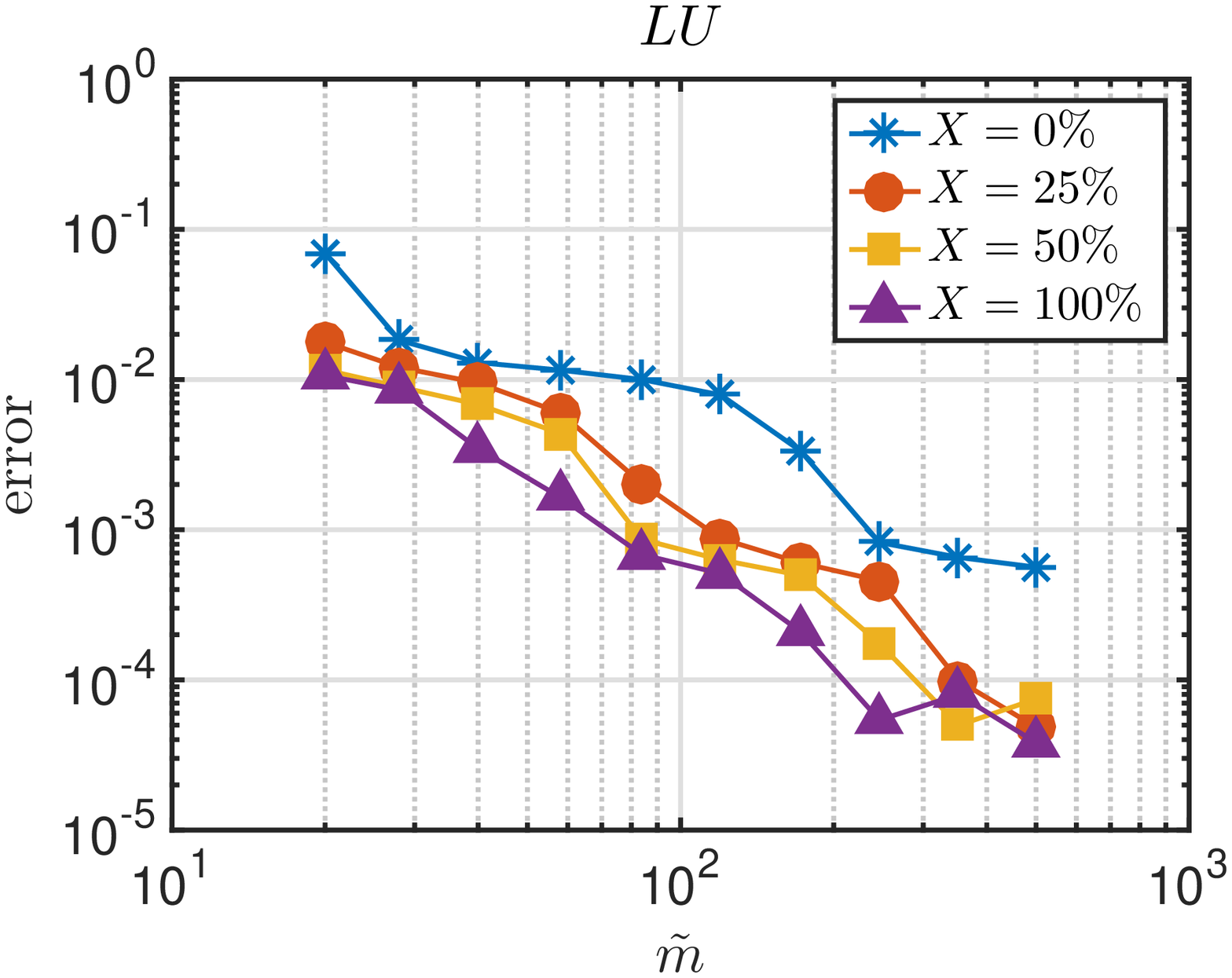} 
     \hspace{30pt}
     \includegraphics[width=2.1in]{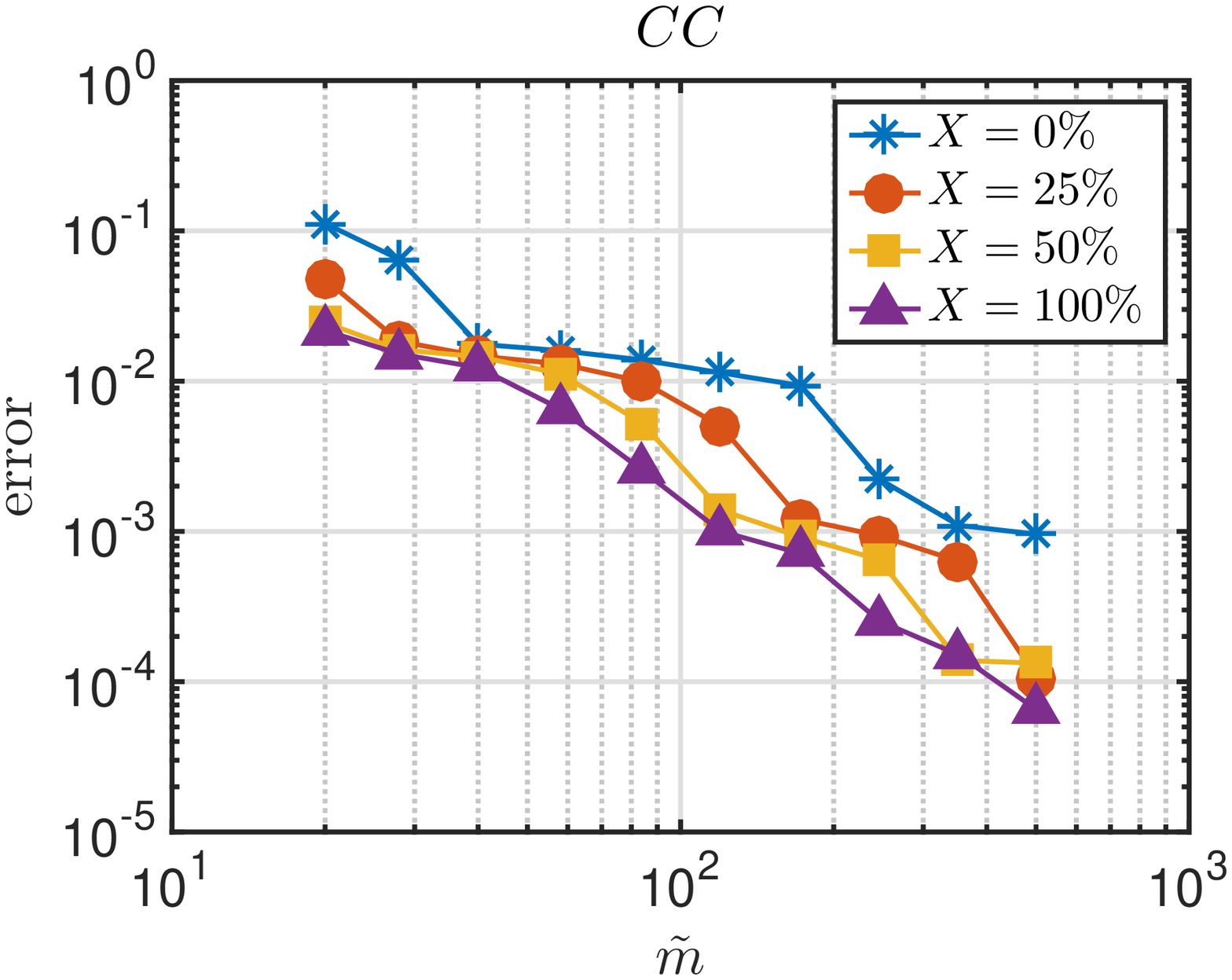}
    \end{tabular}
\caption{The error $\nmu{f_3-\tilde{f}_3}_{\tilde{H}^1(D)}$ against $\tilde{m}$ with a different percentage of gradient enhancement.  The values $(d, s) = (12, 14)$ were used.  The left plot shows the results for Legendre polynomials with uniform sampling and the right plot shows the results for Chebyshev polynomials with Chebyshev sampling.}
\label{d12optper}
\end{figure}

\begin{figure}[ht]
\centering
  \begin{tabular}{cc}
     \includegraphics[width=1.8in]{fig14dwith} 
     \includegraphics[width=1.8in]{fig18dwith}
     \includegraphics[width=1.8in]{fig112dwith}\\
      \includegraphics[width=1.8in]{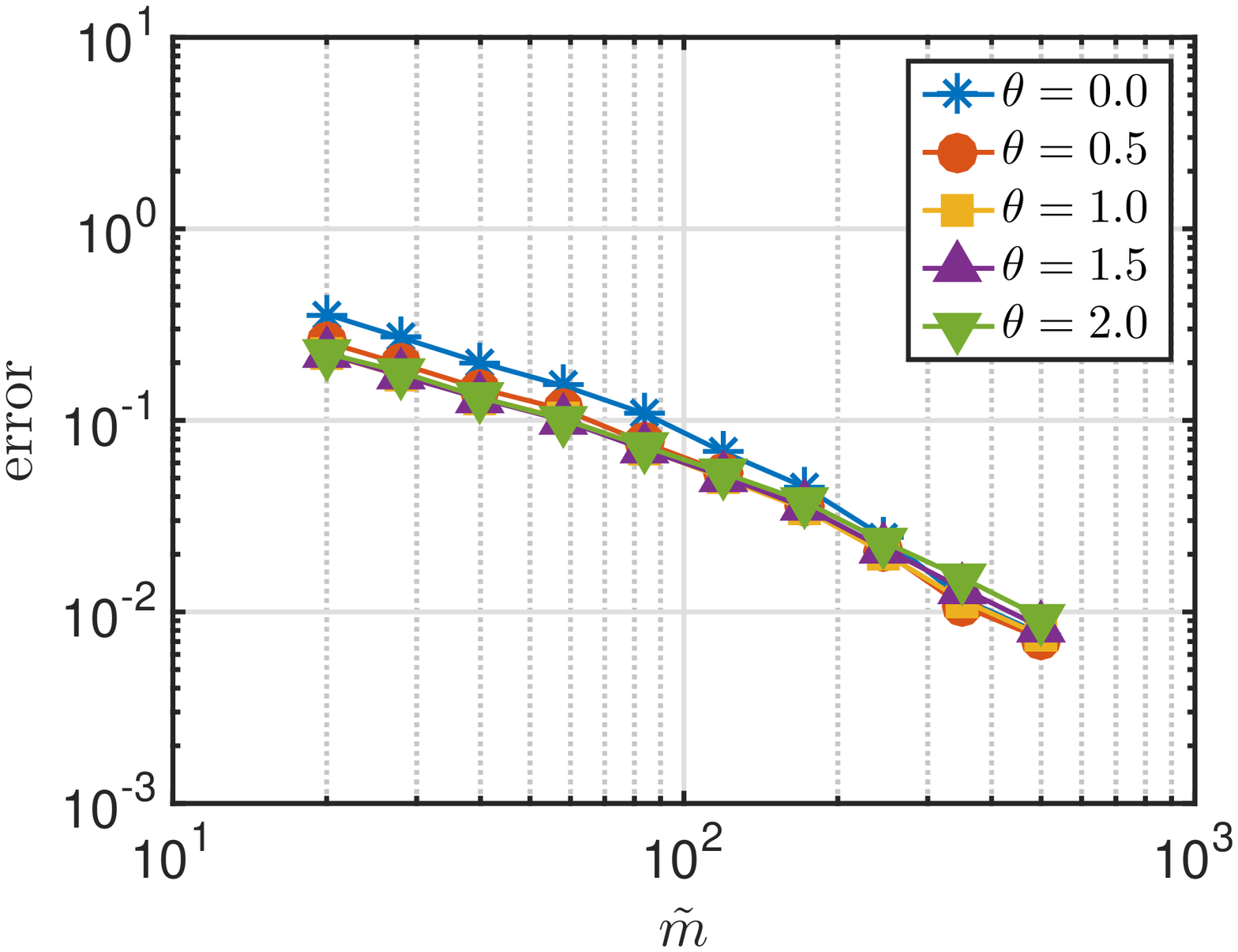} 
     \includegraphics[width=1.8in]{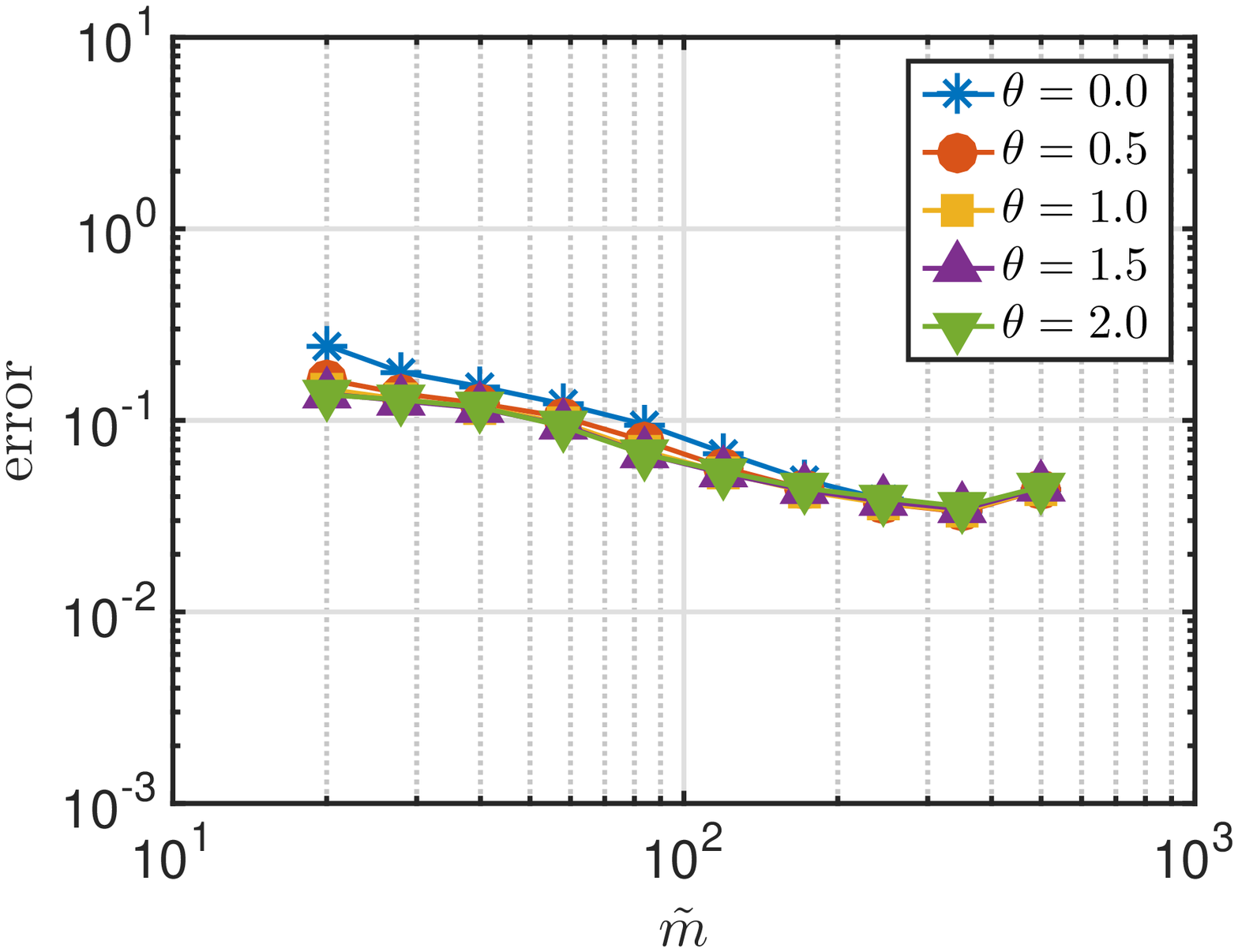}
      \includegraphics[width=1.8in]{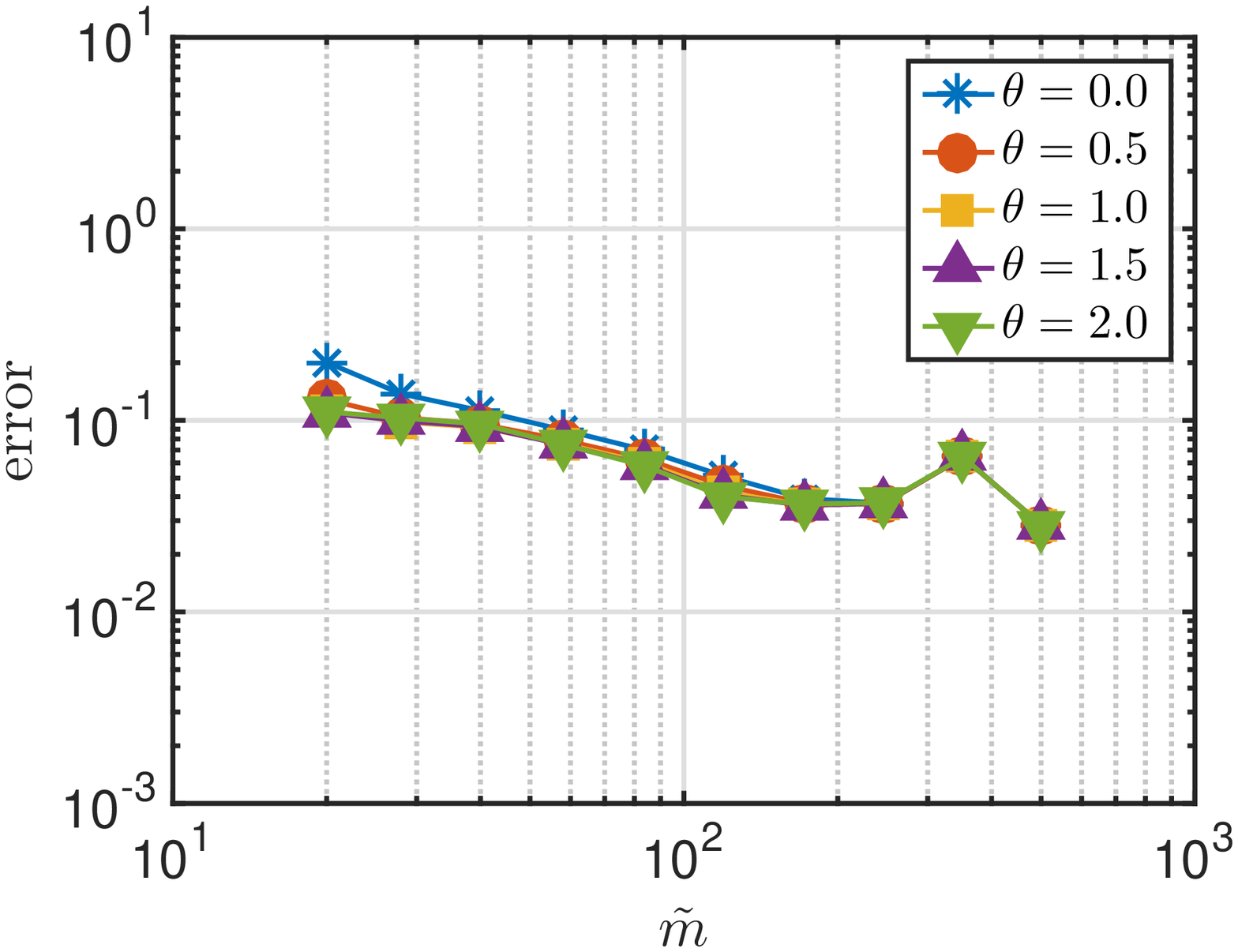}
    \end{tabular}
\caption{The error $\nmu{f_1-\tilde{f}_1}_{\tilde{H}^1(D)}$ against $\tilde{m}$ for Legendre polynomials with points drawn 
from the uniform density. From left to right, the values $(d, s) = (4, 72), (8, 23), (12, 14)$ 
were used.  The top row shows the original setup, and the bottom row shows independent gradient sampling.}
\label{LUoptpeakind}
\end{figure}

Throughout, we solve the weighted $\ell^1$ minimization problem using the SPGL1 package \cite{spgl1:2007, BergFriedlander} with a maximum number of 10,000 iterations and $\eta = 10^{-12}$. 
We choose the truncated index set $\Lambda$ as the hyperbolic cross index set of degree $s$.
 For Figs. \ref{LUoptpeak}--\ref{LUoptpeakind}, the $\tilde{H}^1$ norm error is computed on a fixed grid of 
$4|\Lambda|$ points drawn according to the uniform density for Legendre polynomials and the Chebyshev 
density for Chebyshev polynomials. The error is averaged over 10 trials.

In our first experiments, we take the weights as $w_{\bm{n}} = (u_{\bm{n}})^\theta$ for some $\theta \geq 0$.  We consider the following functions
\eas{
\mbox{Figs.\ \ref{LUoptpeak} \& \ref{CCoptpeak}:}& \qquad f_1(\bm{y}) = \prod_{j=1}^d \frac{d/4}{d/4 + (y_j - a_j)^2}, \quad  a_j = \frac{(-1)^j}{j+1},\qquad 
\\
\mbox{Figs.\ \ref{LUoptcos} \& \ref{CCoptcos}:}& \qquad f_2(\bm{y}) = \prod_{j=d/2+1}^d \cos(16 y_j/2^j)/\prod_{j=1}^{d/2}(1-y_j/4^j)
\\
\mbox{Figs.\ \ref{LUoptexp} \& \ref{CCoptexp}:}& \qquad f_3(\bm{y}) = \exp \left(-\sum_{j=1}^d y_j/(2d) \right)
}
In all dimensions and for all functions, we see that, with the same amount of computational cost $\tilde{m}$, a consistently smaller error is obtained by the
gradient-augmented recovery.  In other words, gradient samples are more beneficial than an equivalent 
number of function samples alone.

Figs.\ \ref{LUoptpeak}--\ref{CCoptexp} also compare different weighting strategies for the optimization problem. 
For the functions we tested, in most cases, the choice $\bm{w} = \bm{u}$ corresponding to $\theta = 1$ gives amongst the smallest, if not the smallest, error.  In particular, these weights often give an improvement over the unweighted case, which corresponds to $\theta = 0$.  This is in agreement with the theoretical results.  Note that larger values of $\theta$ can sometimes give a slightly smaller error depending on the function considered, but the difference is not substantial.

In our next experiment, Fig.\ \ref{d12optper}, we fix the weights as $\bm{w} = \bm{u}$ and consider the scenario where the gradient is measured at only a fixed percentage of the sample points.  A similar setup has also been considered in \cite{Penggradient}. We plot the error versus the effective cost $\tilde{m}$ defined in \R{cost}.  These results show a clear improvement with only $25\%$ gradient samples.  As this percentage increases, the error correspondingly decreases.

\rem{
We conjecture that our theoretical results can be extended to this case as follows.  If $p \in [0,1]$ is the fraction of gradient samples taken, then under the same sample complexity estimate, the error can be bounded in terms of a Sobolev-type norm where the partial derivative terms are weighted by $p$.  In other words, smaller $p$ (fewer gradient samples) corresponds to a weaker norm and larger $p$ (more gradient samples) corresponds to a stronger norm.  This is left as future work.
}

In Fig.\ \ref{LUoptpeakind} we investigate how the location of the gradient samples affects the approximation error.  Specifically, we compare the existing setup where $\nabla f$ is sampled at the same points as $f$ to the case of \textit{independent} gradient sampling locations, i.e.\ where $\nabla f$ is sampled at $m$ points $\bm{y}_{m+1},\ldots,\bm{y}_{2m}$ drawn independently and from the same density as $\bm{y}_1,\ldots,\bm{y}_m$.
As is evident, in all dimensions, independent gradient sampling gives similar recovery results to the original setup for the same computational cost (note we do not take into account here the fact that in practice sampling $\nabla f$ at distinct points may be more expensive).  Thus, there is apparently little benefit to sampling the gradient at a distinct set of sample points.

\rem{
We conjecture that all our theoretical results can all be adapted to the case of independent gradient sampling, with potentially only minor changes to the log factors.
}

\begin{figure}[ht]
\centering
  \begin{tabular}{cc}
     \includegraphics[width=1.8in]{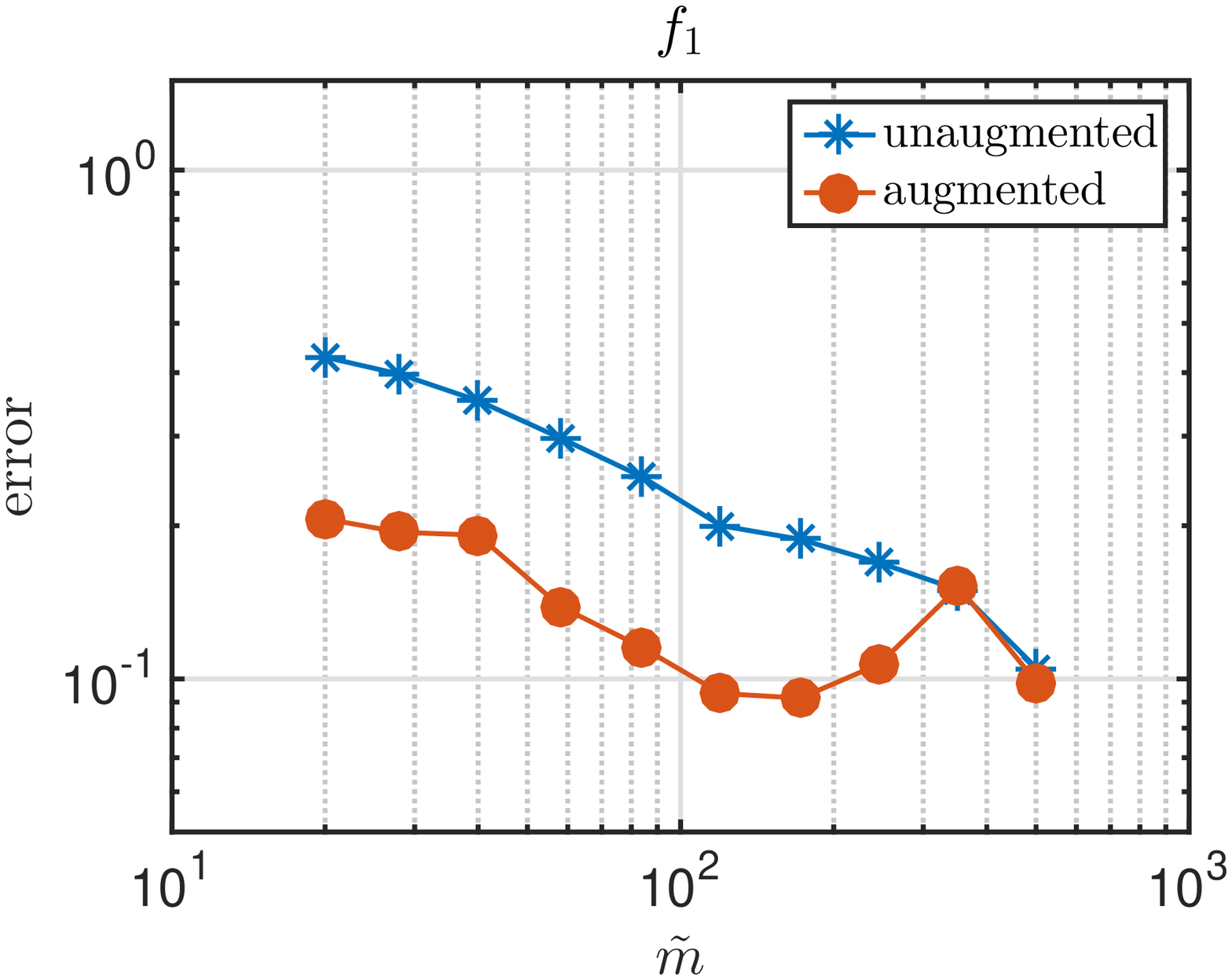} 
     \includegraphics[width=1.8in]{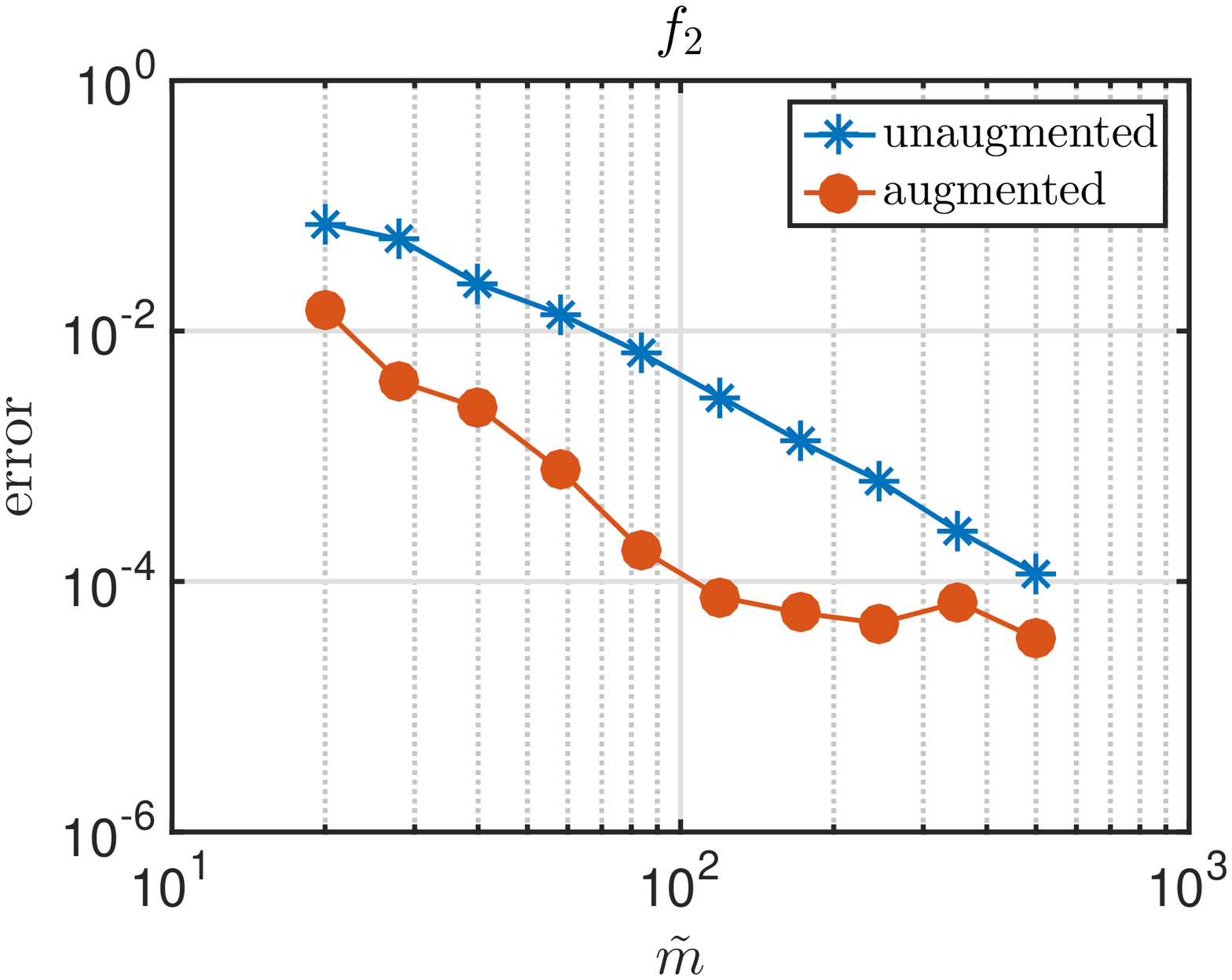}
     \includegraphics[width=1.8in]{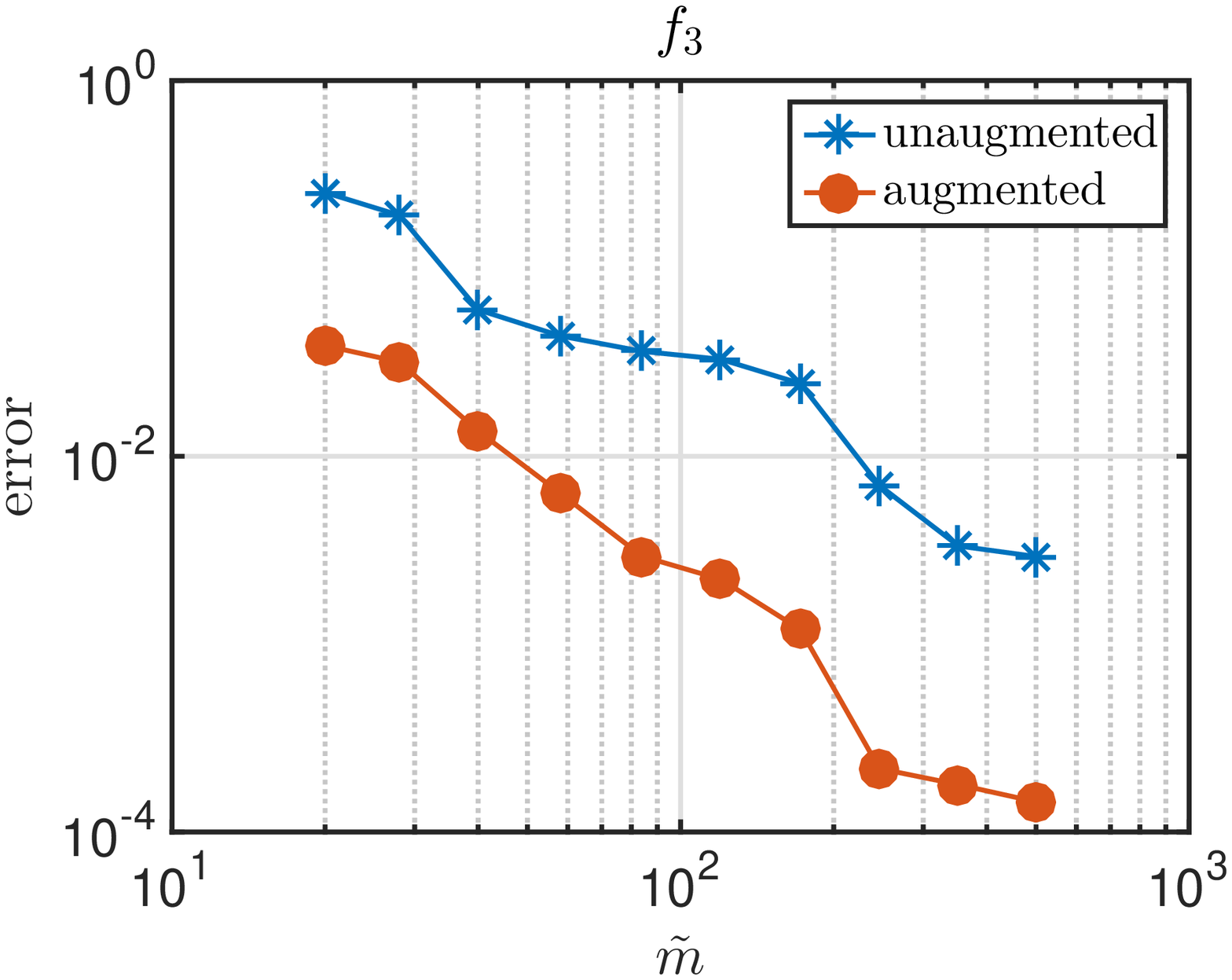}\\
      \includegraphics[width=1.8in]{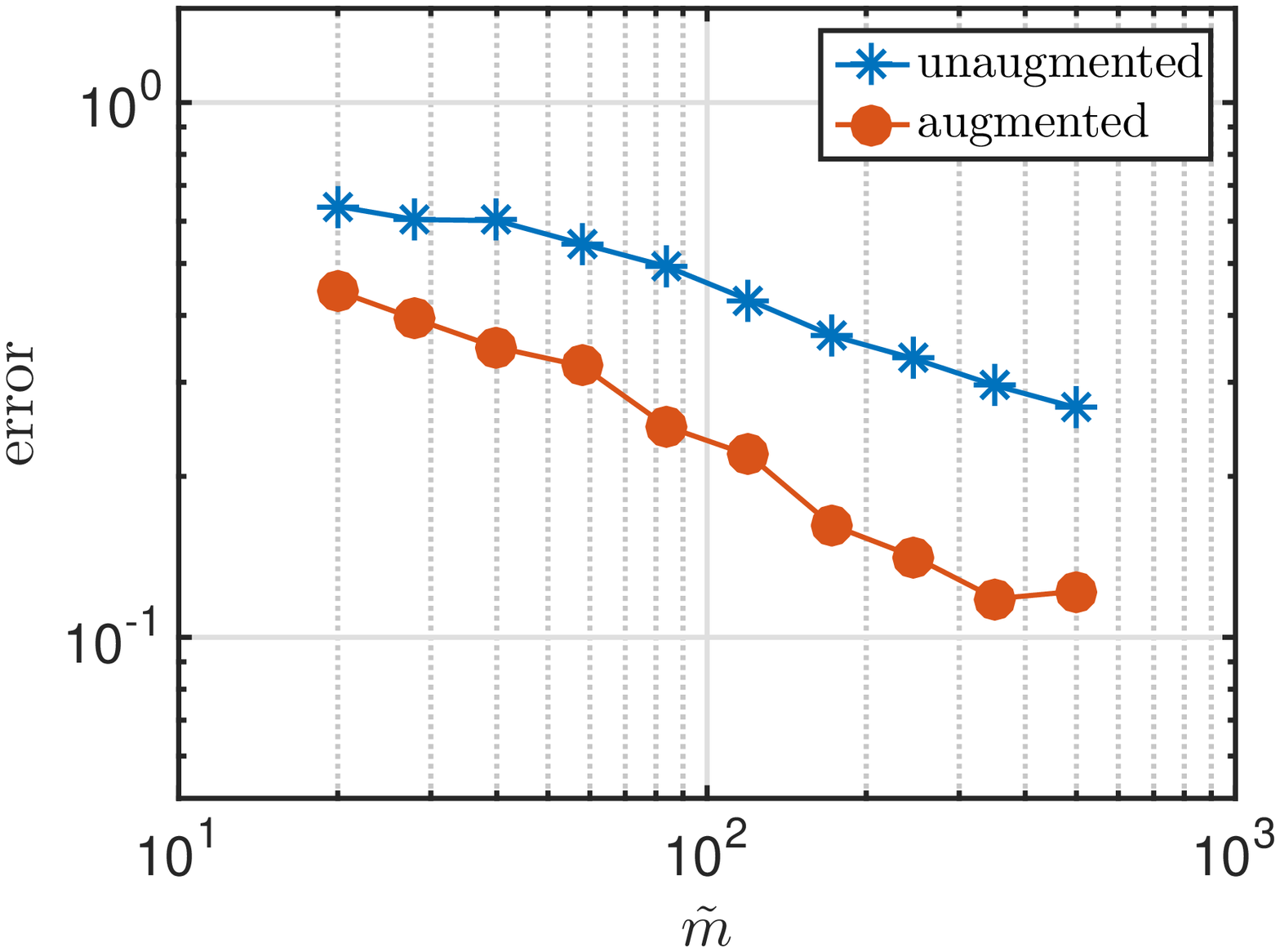} 
     \includegraphics[width=1.8in]{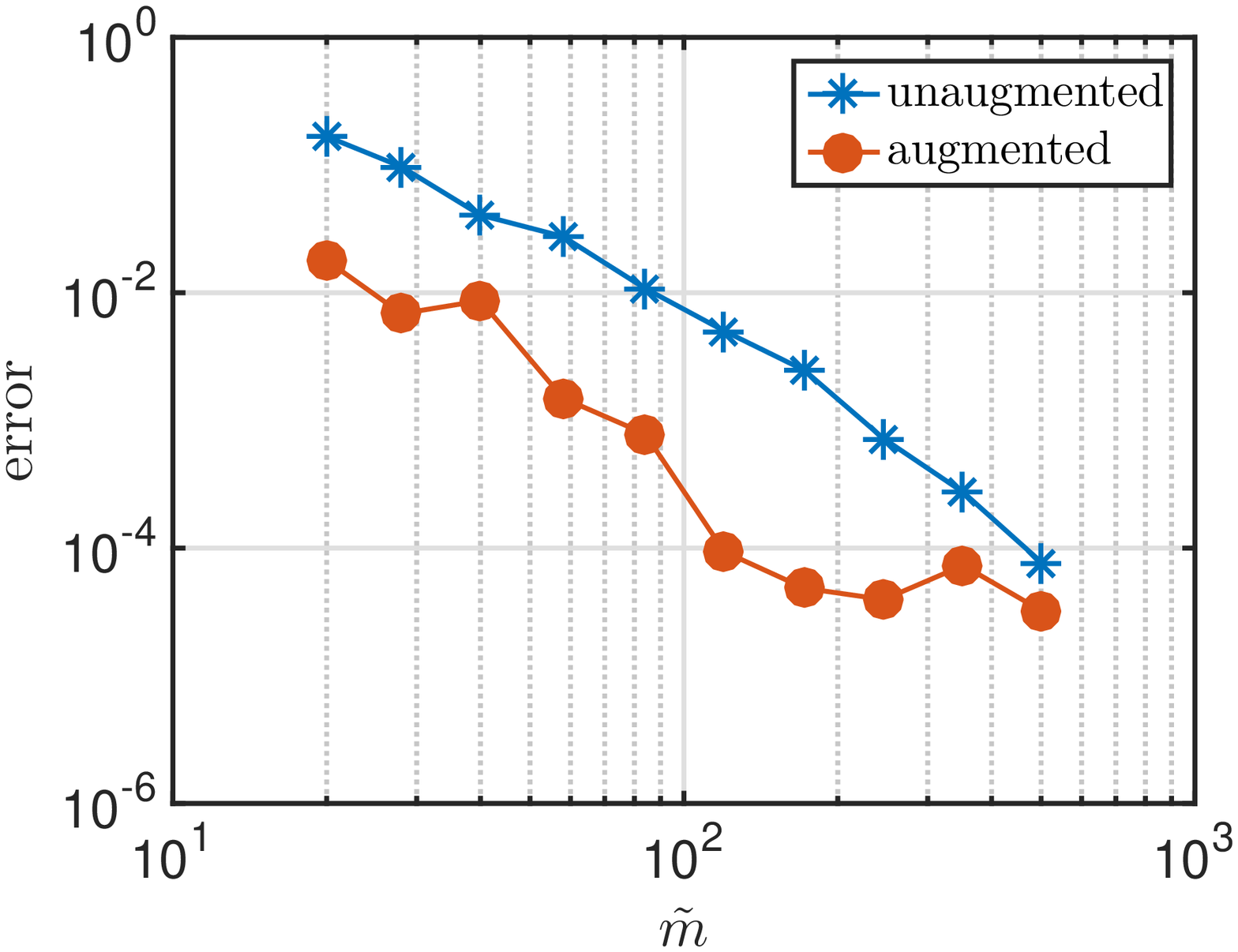}
      \includegraphics[width=1.8in]{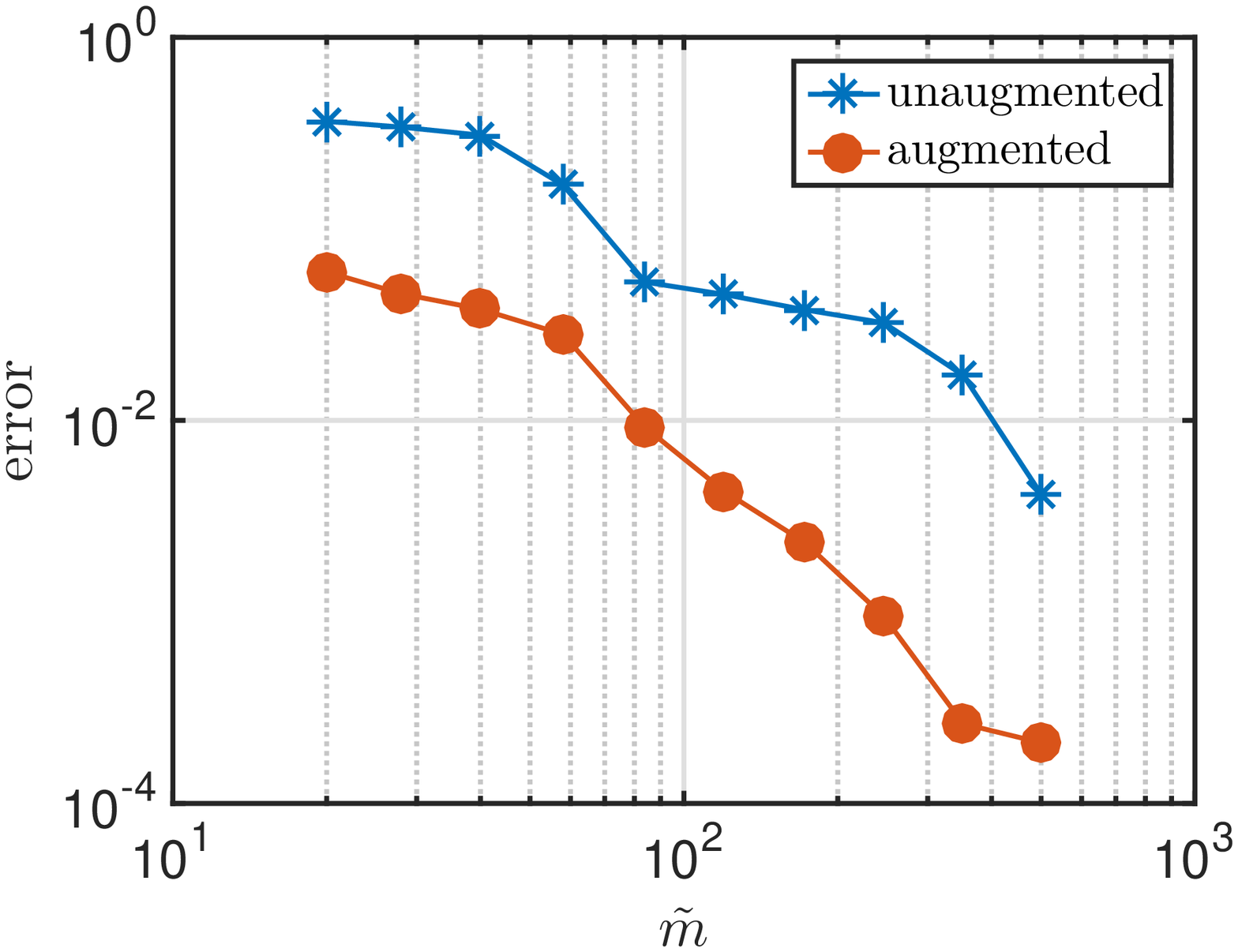}
    \end{tabular}
\caption{The error $\nmu{f_1-\tilde{f}_1}_{L^\infty}$ against $\tilde{m}$ for Legendre polynomials with 
points drawn from the uniform density (top) and Chebyshev polynomials with points drawn from the Chebyshev density (bottom). function $f_1$ to $f_3$ are shown from left to right.  The value $(d, s) = (12, 14)$ 
was used to generate the index set.}
\label{Linfty}
\end{figure}

Finally, in Fig. \ref{Linfty}, we compare the $L^\infty$-norm error for the unaugmented and gradient-augmented cases. Here, the error is computed on a fixed grid of $4|\Lambda|$ uniformly-distributed 
points and averaged over 10 trials. We fix the weights $\bm{w} = \bm{u}$. 
As in Figs.\ \ref{LUoptpeak}--\ref{CCoptexp}, we see that, with the same amount of 
computational cost,  the gradient-augmented recovery leads to a smaller error in the $L^\infty$ norm.

\section{Proofs}\label{s:proofs}

In this final section, we give the proofs of the main results.  To this end, we first show that the above problem can be reformulated as an instance of the general `parallel acquisition' compressed sensing framework of \cite{AdcockChunParallel} (see also \cite{BigotBlockCS,BoyerBlockStructured}).  This allows us to use the approach of \cite{AdcockChunParallel} (with modifications to take into account the weighted regularizer) to prove the recovery guarantees.

\subsection{The framework of \cite{AdcockChunParallel}}\label{ss:par_acq}
We follow the setup described in \cite[\S II-D]{AdcockChunParallel}.  For some 
 $\mathcal{D} \in \bbN$, let $F$ be a distribution on a set of $N \times   \mathcal{D}$ complex matrices.  We assume that $F$ is isotropic in the sense that
\bes{
\bbE (B B^*) = I, \qquad B \sim F.
}
Now let $\{ \bm{e}_i \}_{i=1}^{m}$ be the canonical basis of $\bbC^m$ and let $B_1,\ldots,B_m$ be a sequence of independent realizations of matrices from the distribution $F$.  Then we define the sampling matrix
\be{
\label{A_par_acq}
A = \frac{1}{\sqrt{m}} \sum^{m}_{i=1} \bm{e}_i \otimes B^*_i = \frac{1}{\sqrt{m}} \begin{bmatrix} B^*_1 \\ \vdots \\ B^*_m \end{bmatrix} \in \bbC^{\mathcal{D} m \times N}.
}
Note that this is an extension of the standard compressed sensing setup, which corresponds to the case $
\mathcal{D} = 1$, i.e.\ $A$ having independent rows.
The paper \cite{AdcockChunParallel} considered compressed sensing for this model of measurement matrices using $\ell^1$-minimization and proved a series of nonuniform recovery guarantees.  In what follows, we consider the generalization of this setup to the weighted $\ell^1$-minimization problem
\be{
\label{par_acq_l1umin}
\min_{\bm{z} \in \bbC^N} \nm{\bm{z}}_{1,\bm{w}}\ \mbox{subject to $\nm{A \bm{z} - \bm{y}}_{2} \leq \eta$},
}
where  $\bm{w} = (w_i)^{N}_{i=1} \in \bbR^N$ with $w_i \geq 1$, $\forall i$.
Here $\bm{y} = A \bm{x} + \bm{e}$ are noisy measurements of the unknown vector $\bm{x}$ (for ease of notation we write this rather than $\bm{x}_{\Lambda}$) and $\bm{e}$ is a vector satisfying $\nm{\bm{e}}_{2} \leq \eta$.

\subsection{Derivatives sampling as an instance of the parallel acquisition model}\label{ss:deriv_PA}
Consider the setup of \S \ref{s:gradaug}.  For the random variable $\bm{y}$ with probability density $\mu$ on $D$, define the random matrix
\be{
\label{B_deriv}
B = \left ( \sqrt{\tau_k(\bm{y})} \frac{\overline{\partial_{k} \phi_{\bm{n}_j}(\bm{y})} }{\sqrt{1+\lambda_{\bm{n}_j}}} \right )^{N,d}_{j=1,k=0} \in \bbC^{N \times (d+1)}.
}
This gives rise to a distribution $F$ on random matrices in $\bbC^{N \times  \mathcal{D}}$, 
where $ \mathcal{D} = (d+1)$.  Moreover, the corresponding matrix \R{A_par_acq} is (after a permutation of its rows) identical to the matrix defined in \R{A_deriv_prob}.  
Since the constraint $\nm{A \bm{z} - \bm{y}}_{2} \leq \eta$ is unaffected by row permutations, we deduce that the derivatives recovery problem \R{l1u_min_deriv} is a particular instance of the above framework, corresponding to choice $ \mathcal{D} = (d+1)$ and with $F$ being the distribution of matrices \R{B_deriv}.

\subsection{The parallel acquisition model with weighted $\ell^1$ minimization}

In order to prove our main result concerning derivative sampling, we first establish a general result for the model of \S \ref{ss:par_acq} with the weighted $\ell^1$ regularizer \R{par_acq_l1umin}, thereby generalizing the result of \cite{AdcockChunParallel}.  First, we require some notation.  
If $\Delta \subseteq \{1,\ldots,N\}$, then we use the notation $P_{\Delta}$ for the orthogonal projection $P_{\Delta} \in \bbC^{N \times N}$ onto $\spn \{ \bm{e}_j : j \in \Delta \}$.  We note in passing that $P_{\Delta} \bm{x} \in \bbC^{N}$ is isomorphic to a vector in $\bbC^{|\Delta|}$.  Also, given weights $\bm{w} \in \bbR^N$ we write $W = \diag(\bm{w})$.  Finally, we note that in this section we index over $\bbN$ where relevant, as opposed to $\bbN^d_0$ as in the original polynomial approximation problem.

Our first step, as in \cite{AdcockChunParallel}, is to define several notions of local coherence:

\defn{
\label{d:coh_rel_1}
Let $\Delta \subseteq \{1,\ldots,N\}$ and $F$ be as in  \S \ref{ss:par_acq}.  The local coherence of $F$ relative to $\Delta$ is the smallest constant $\Upsilon(F,\Delta)$ such that
\bes{
\nm{P_{\Delta} B B^* P_{\Delta} }_{2} \leq \Upsilon(F,\Delta),\qquad B \sim F,
}
almost surely.
}

\defn{
\label{d:coh_rel_2}
Let  $\bm{w} \in \bbR^N$ be a set of positive weights, $\Delta \subseteq \{1,\ldots,N\}$ and $F$ be as in 
\S \ref{ss:par_acq}.  The local coherence of $F$ relative to $\Delta$ with respect to the weights $\bm{w}$ is
\bes{
\Gamma(F,\bm{w},\Delta) = \max \left \{ \Gamma_1(F,\bm{w},\Delta) , \Gamma_2(F,\bm{w},\Delta)  \right \},
}
where $\Gamma_{1}(F,\bm{w},\Delta)$ and $\Gamma_{2}(F,\bm{w},\Delta)$ are the smallest quantities such that
\bes{
 \nm{W^{-1} B B^* P_{\Delta} W}_{\infty} \leq \Gamma_{1}(F,\bm{w},\Delta) ,\qquad B \sim F,
 }
almost surely, and
\bes{
 \sup_{\nm{\bm{z}}_{\infty} = 1} \max_{j=1,\ldots,N} \bbE | \ip{\bm{e}_j}{W^{-1} B B^* P_{\Delta} W \bm{z} } |^2 \leq \Gamma_{2}(F,\bm{w},\Delta).
}
}

By definition, if $j \in \Delta$ then
\bes{
\Gamma_1(F,\bm{w},\Delta) \geq \bbE | \ip{\bm{e}_j}{W^{-1} B B^* P_{\Delta} W \bm{e}_j } | \geq \left | \bbE \ip{\bm{e}_j}{W^{-1} B B^* P_{\Delta} W \bm{e}_j } \right | = | \ip{\bm{e}_j}{W^{-1} P_{\Delta} W \bm{e}_j } | = 1.
}
Hence we deduce that $\Gamma_1(F,\bm{w},\Delta) \geq 1$.  Similarly, we also have $\Gamma_2(F,\bm{w},\Delta) \geq 1$ and the same for the unweighted local coherence $\Upsilon(F,\Delta) \geq 1$.

Our main result for the abstract model of \S \ref{ss:par_acq} is now as follows:

\thm{
\label{t:par_acq_wl1min}
Let $0 < \epsilon < 1$, $\eta \geq 0$, $N \geq 2$, $\Delta \subseteq \{1,\ldots,N\}$ with $|\Delta | \geq 2$ and $\bm{w} \in \bbR^N$ be weights with $w_i \geq 1$, $\forall i$.  Fix $\bm{x} \in \bbC^N$ and construct $A \in \bbC^{m  \mathcal{D} \times N}$ as in \R{A_par_acq}.  Let $\bm{y} = A \bm{x} + \bm{e}$, where $\nm{\bm{e}}_{2} \leq \eta$.  If
\bes{
m \gtrsim \Upsilon(F,\Delta) \cdot \log(N/\epsilon) + \Gamma(F,\bm{w},\Delta) \cdot \left ( \log(N/\epsilon) + \log(|\Delta|_{\bm{w}}) \cdot  \log(|\Delta|_{\bm{w}}/\epsilon) \right ),
}
where $\Upsilon(F,\Delta)$ and $\Gamma(F,\bm{w},\Delta)$ are as in Definitions \ref{d:coh_rel_1} and \ref{d:coh_rel_2} respectively, then, with probability at least $1-\epsilon$, any minimizer $\bm{\hat{x}}$ of \R{par_acq_l1umin} satisfies
\bes{
\nm{\bm{x} - \bm{\hat{x}} }_{2} \lesssim \nm{\bm{x} - P_{\Delta} \bm{x}}_{1,\bm{w}} + \sqrt{|\Delta|_{\bm{w}}} \eta.
}

}

\subsection{Proof of Theorem \ref{t:par_acq_wl1min}}

The proof follows that of \cite[Thm. 12]{AdcockChunParallel}, making changes where necessary to account for the weighted regularizer.  Note that the particular case of the weighted regularizer with 
$ \mathcal{D} = 1$ (i.e.\ no derivatives in the case of function approximation) was essentially covered in \cite{AdcockCSFunInterp}.  The arguments we use next effectively combine those of \cite{AdcockChunParallel} and \cite{AdcockCSFunInterp} to yield Theorem \ref{t:par_acq_wl1min}.  For this reason, we only sketch the details, making references to the relevant parts of \cite{AdcockChunParallel} and \cite{AdcockCSFunInterp} wherever necessary.

We first require a series of technical lemmas:

\lem{
\label{l:tech1}
Let $0 < \epsilon < 1$, $\delta > 0$, $F$ and $A \in \bbC^{m \mathcal{D} \times N}$ be as in \S \ref{ss:par_acq} and suppose that $\Delta \subseteq \{1,\ldots,N\}$.  Then
\bes{
\| P_{\Delta} A^* A P_{\Delta} - P_{\Delta} \|_2 < \delta,
}
with probability at least $1-\epsilon$, provided
\bes{
m  \geq \Upsilon(F,\Delta) \cdot (2 \delta^{-2} + 2 \delta^{-1} / 3 ) \cdot \log(2 | \Delta | / \epsilon ),
}
where $\Upsilon(F,\Delta)$ is as in Definition \ref{d:coh_rel_1}.
}

This is identical to  \cite[Lem.\ 41]{AdcockChunParallel}, and hence its proof is omitted.  The following lemma is a straightforward extension of \cite[Lem.\ 42]{AdcockChunParallel} to the weighted setting:
\lem{
\label{l:tech2}
Let $0 < \epsilon < 1$, $\delta > 0$, $F$ and $A \in \bbC^{m  \mathcal{D} \times N}$ be as in \S \ref{ss:par_acq} and suppose that $\Delta \subseteq \{1,\ldots,N\}$ and $\bm{z} \in \bbC^N$.  Then
\bes{
\| W^{-1} (P_{\Delta} A^* A P_{\Delta} - P_{\Delta}) W \bm{z} \|_{\infty} < \delta \| \bm{z} \|_{\infty},
}
with probability at least $1-\epsilon$, provided
\bes{
m \geq \left ( 8 \Gamma_1(F,\bm{w},\Delta) \delta^{-1} /3+ 4 \Gamma_2(F,\bm{w},\Delta) \delta^{-2} \right ) \cdot \log(4 |\Delta| / \epsilon),
}
where $\Gamma_{1}(F,\bm{w},\Delta)$ and $\Gamma_{2}(F,\bm{w},\Delta)$ are as in Definition \ref{d:coh_rel_2}. 
}
\prf{
Let $\nm{\bm{z}}_{\infty} = 1$ without loss of generality.  Fix $j \in \Delta$ and observe that
\bes{
\ip{\bm{e}_j}{W^{-1} (P_{\Delta} A^* A P_{\Delta} - P_{\Delta} )W \bm{z} } = \frac{1}{m} \sum^{m}_{i=1} \ip{\bm{e}_j}{W^{-1}  \left ( B_i B^*_i - I \right ) P_{\Delta} W \bm{z}} = \frac{1}{m} \sum^{m}_{i=1} X_i,
}
where $X_i$ is the random variable $X_i =  \ip{\bm{e}_j}{W^{-1}  \left ( B_i B^*_i - I \right ) P_{\Delta} W \bm{z}}$.  Note that $\bbE(X_i) = 0$.  Also
\bes{
| X_i | = | \ip{\bm{e}_j}{W^{-1} (B_i B^*_i - I) P_{\Delta} W \bm{z}} | \leq \nm{W^{-1} P_{\Delta} B_i B^*_i P_{\Delta} W}_{\infty} +1 \leq \Gamma_{1}(F,w,\Delta)+1,
}
and
$
\bbE |X_i |^2 = \bbE | \ip{\bm{e}_j}{W^{-1} B_i B^*_i P_{\Delta} W \bm{z} } |^2 - | \ip{e_j}{P_{\Delta} 
 \bm{z} } |^2 \leq \Gamma_{2}(F,w,\Delta).
$
An application of Bernstein's inequality followed by the union bound now yields
\bes{
\bbP \left ( \| W^{-1} (P_{\Delta} A^* A P_{\Delta} - P_{\Delta}) W \bm{z} \|_{\infty}  \geq \delta \right ) \leq 4 | \Delta | \exp \left ( - \frac{m \delta^2/4}{\Gamma_2(F,\bm{w},\Delta) + 2 \Gamma_1(F,\bm{w},\Delta) \delta / 3} \right ).
}
Equating the right hand side with $\epsilon$ and rearranging gives the result.
}

\lem{
\label{l:tech3}
Let $0 < \epsilon < 1$, $\delta > 0$, $F$ and $A \in \bbC^{m  \mathcal{D} \times N}$ be as in \S \ref{ss:par_acq}.  Suppose that $\Delta \subseteq \{1,\ldots,N\}$ and  $\bm{w} \in \bbR^N$ is a vector of weights with $w_i \geq1$.  Then
\bes{
\max_{j \notin \Delta} \nm{P_{\Delta} A^* A W^{-1} \bm{e}_j}_{2}  \leq \delta,
}
with probability at least $1-\epsilon$, provided
\bes{
m \geq \left( \Upsilon(F,\Delta) + \Gamma_{1}(F,\bm{w},\Delta) \right )\cdot \left ( 8 \delta^{-2} + 28 \delta^{-1} / 3 \right ) \cdot \log(2 N/\epsilon).
}
where $\Upsilon(F,\Delta)$ and $\Gamma_{1}(F,\bm{w},\Delta)$ are as in Definitions \ref{d:coh_rel_1} and \ref{d:coh_rel_2} respectively. 
}
\prf{
Fix $j \notin \Delta$.  Then
\bes{
\nm{P_{\Delta} A^* A W^{-1} \bm{e}_j}_{2} = \frac{1}{m} \nm{\sum^{m}_{i=1} \bm{X}_i}_2,
}
where $\bm{X}_i = P_{\Delta} B_i B^*_i W^{-1} \bm{e}_j$ are independent copies of the random vector $\bm{X} = P_{\Delta} B B^* W^{-1} \bm{e}_j$.  We have $\bbE(X) = 0$ since $j \notin \Delta$.  Moreover,
\bes{
\nm{\bm{X}}_{2} = \| P_{\Delta} B B^* W^{-1} \bm{e}_j \|_{2} \leq \| W^{-1} B B^* P_{\Delta} W \|_{\infty} \| W^{-1} \|_{\infty} \leq \Gamma_{1}(F,\bm{w},\Delta),
}
since $w_i \geq 1$, and
\bes{
\bbE \nm{\bm{X}}^2_2 \leq \bbE \left( \nm{P_{\Delta} B}^2_{2} \nm{B^* W^{-1} \bm{e}_j }^2_2 \right ) \leq \Upsilon(F,\Delta) \| W^{-1} \bm{e}_j \|^2_2 \leq \Upsilon(F,\Delta).
}
We now argue as in \cite[Lem.\ 43]{AdcockChunParallel}.
}

The next lemma extends \cite[Lem.\ 44]{AdcockChunParallel}:
\lem{
\label{l:tech4}
Let $0 < \epsilon < 1$, $\delta > 0$, $F$ and $A \in \bbC^{m  \mathcal{D} \times N}$ be as in \S \ref{ss:par_acq}.  Suppose that $\Delta \subseteq \{1,\ldots,N\}$, $\bm{z} \in \bbC^N$ and 
 $\bm{w} \in \bbR^N$ is a vector positive weights.  Then
\bes{
\nm{P^{\perp}_{\Delta} W^{-1} A^* A P_{\Delta} W \bm{z} }_{\infty} \leq \nm{\bm{z}}_{\infty},
}
with probability at least $1-\epsilon$, provided
\bes{
m \geq \left ( 4 \Gamma_1(F,\bm{w},\Delta) \delta^{-1}/3 + 4 \Gamma_2(F,\bm{w},\Delta) \delta^{-2} \right ) \cdot \log(2 N / \epsilon),
}
where $\Gamma_{1}(F,\bm{w},\Delta)$ and $\Gamma_{2}(F,\bm{w},\Delta)$ are as in Definition \ref{d:coh_rel_2}. 
}
\prf{
We assume $\nm{z}_{\infty} = 1$ without loss of generality and fix $j \notin \Delta$.  Then
\bes{
\ip{\bm{e}_j}{W^{-1} A^* A P_{\Delta} W \bm{z}} = \ip{\bm{e}_j}{W^{-1} (A^* A-I) P_{\Delta} W \bm{z}} = \frac1m \sum^{m}_{i=1} X_i,
}
where $X_i$ is the random variable $X_i = \ip{\bm{e}_j}{W^{-1} (B_i B^*_i - I) W P_{\Delta } \bm{z}}$.  Note that
\bes{
| X_i | \leq \nm{ W^{-1} B_i B^*_i W P_{\Delta} z }_{\infty} + 1 \leq \Gamma_{1}(F,\bm{w},\Delta) + 1 \leq 2 \Gamma_1(F,\bm{w},\Delta).
}
Also,
$
\bbE |X_i |^2  = \bbE |  \ip{\bm{e}_j}{W^{-1} B_i B^*_i  P_{\Delta } W  \bm{z}} |^2 - | \ip{\bm{e}_j}{P_{\Delta} z } |^2 \leq \Gamma_{2}(F,\bm{w},\Delta).
$
The result now follows from Bernstein's inequality and the union bound.
}

Finally, we require the following lemma (see \cite[Lem.\ 8.1]{AdcockCSFunInterp}):

\lem{
\label{l:dual_certificate}
Let  $\bm{x} \in \bbC^N$ and $\bm{w} \in \bbR^N$ be weights with $w_i \geq 1$, $\Delta \subseteq \{1,\ldots,N\}$ and $A \in \bbC^{m \times N}$.  Suppose that
\bes{
(i): \| P_{\Delta} A^*  A P_{\Delta} - P_{\Delta} \|_2 \leq \alpha,\qquad (ii):  \max_{i \notin \Delta} \|  P_{\Delta} A^* A W^{-1} \bm{e}_i \|_{2}  \leq \beta,
}
and that there exists a vector $\bm{\rho} = W^{-1} A^*  \bm{\xi} \in \bbC^N$ for some $\bm{\xi} \in \bbC^m$ such that
\bes{
(iii): \| W(P_{\Delta} \bm{\rho} - \sgn(P_{\Delta} \bm{x})) \|_2 \leq \gamma,\quad (iv): \| P^{\perp}_{\Delta} \bm{\rho} \|_{\infty} \leq \theta,\quad (v): \| \bm{\xi} \|_2 \leq \lambda \sqrt{| \Delta |_{\bm{w}}}, 
}
for constants $0 \leq \alpha , \theta < 1$ and $\beta , \gamma, \lambda \geq 0$ satisfying $\frac{\sqrt{1+\alpha} \beta \gamma}{(1-\alpha)(1-\theta) } < 1$.  Let $\bm{y} = A \bm{x} + \bm{e}$ with $\| \bm{e} \|_2 \leq \eta$ and suppose that $\hat{\bm{x}}$ is a minimizer of the problem
\bes{
\min_{\bm{z} \in \bbC^N} \| \bm{z} \|_{1,\bm{w}}\  \mbox{subject to $\| A \bm{z} - \bm{y} \|_2 \leq \eta$.}
}
Then
\be{
\label{l1_error_est}
\| \bm{\hat{x}} - \bm{x} \|_2 \leq C_1 \left (1+\lambda \sqrt{|\Delta|_{\bm{w}}} \right ) \eta + C_2 \nm{\bm{x} - P_{\Delta} \bm{x} }_{1,\bm{w}},
}
where the constants $C_1$ and $C_2$ depend on $\alpha$, $\beta$, $\gamma$ and $\theta$ only.
}

\prf{[Proof of Theorem \ref{t:par_acq_wl1min}]
We follow the proof given in \cite[Thm.\ 12]{AdcockChunParallel}.  Our strategy is to use the so-called golfing scheme \cite{Gross} to construct a vector $\bm{\rho}$ so that Lemma \ref{l:dual_certificate} holds for appropriate parameters, which we arbitrarily take to be
\bes{
\alpha = 1/4,\quad \beta = 1,\quad \gamma = 1/4,\quad \theta = 1/2.
}
Recall that $| \Delta | \geq 2$.  In particular, $\log(|\Delta|) \geq \log(2) > 0$.  First, let $s^* = | \Delta |_{\bm{w}}$ and define
\be{
\label{L_choice}
L = 2+ \left \lceil \log_2(\sqrt{s^*}) \right \rceil \geq 3,
}
here we recall that $| \Delta |_{\bm{w}} \geq | \Delta | \geq 2$ since $w_i \geq 1$ for $\forall i$,
\be{
\label{a_choice}
a_1 = a_2 = \frac{1}{2 \sqrt{\log_2(\sqrt{s^*})} },\qquad a_l = 1/2,\quad l =3,\ldots,L,
}
\be{
\label{b_choice}
b_1 = b_2 = \frac{1}{4},\qquad b_l = \frac{\log_2(\sqrt{s^*})}{4} ,\quad l =3,\ldots,L,
}
and
\bes{
m_1 = m_2 = \left \lceil \frac14 m^* \right \rceil,\qquad m_l = \left \lceil \frac{1}{2(L-2)} m^* \right \rceil,\quad l=3,\ldots,L.
}
where $m^* = \lfloor m - L \rfloor$.  Observe that
\bes{
\sum^{L}_{l=1} m_l \leq m^* + L \leq m.
}
We now let
\bes{
A_{l} = \frac{1}{\sqrt{m_l}} \sum^{m_1+ \ldots + m_l}_{i=m_1 + \ldots + m_{l-1} + 1} \bm{e}_i \otimes B^*_i \in \bbC^{m \cD \times N},\qquad l=1,\ldots,L,
}
and notice that
\bes{
A = \sum^{L}_{l = 1} \sqrt{m_l/m} A_l.
}
The dual certificate is now constructed iteratively as follows.  Let $\bm{\rho}^{(0)} = 0$,
\bes{
\bm{\rho}^{(l)} =  W^{-1} (A_{l})^* A_l P_{\Delta} W \left ( \sgn(P_{\Delta}(\bm{x})) - P_{\Delta} \bm{\rho}^{(l-1)} \right ) + \bm{\rho}^{(l-1)},\qquad l=1,\ldots,L,
}
and set $\bm{\rho} = \bm{\rho}^{(L)}$.  

With this in hand, we define the vector $\bm{v}^{(l)}$ as 
\bes{
\bm{v}^{(l)} =  W \left ( \sgn(P_{\Delta} \bm{x}) - P_{\Delta} \bm{\rho}^{(l)} \right ),\qquad l=0,\ldots,L,
}
and consider the following events:
\eas{
A_l &:\quad \| W^{-1}(P_{\Delta} -  P_{\Delta} (A_{l})^* A_l P_{\Delta}) \bm{v}^{(l-1)} \|_{\infty} \leq a_l \| W^{-1} \bm{v}^{(l-1)} \|_{\infty},\qquad l=1,\ldots,L,
\\
 B_l &: \quad \|  P^{\perp}_{\Delta} W^{-1}  (A_{l})^* A_l P_{\Delta} \bm{v}^{(l-1)} \|_{\infty} \leq b_l \| W^{-1}\bm{v}^{(l-1)} \|_{\infty},\qquad l=1,\ldots,L.
\\
 C &: \quad \| P_{\Delta} A^* A P_{\Delta} - P_{\Delta} \|_2 \leq 1/4,
\\
 D &:\quad \max_{i \notin \Delta} \| P_{\Delta} A^* A W^{-1} \bm{e}_i \|_2  \leq 1,
\\
 E &: \quad A_1 \cap \cdots \cap A_L \cap B_1 \cap \cdots \cap B_L \cap C \cap D.
}
We now proceed in two steps: first, showing that event $E$ implies conditions (i)--(v) of Lemma \ref{l:dual_certificate}, and second, showing that event $E$ holds we high probability.

\pbk
\textbf{Step 1.}\ If event $E$ occurs, then events $C$ and $D$ give (i) and (ii) respectively.  Next consider (iii).  Observe that
\bes{
\bm{v}^{(l)} = W \sgn(P_{\Delta} \bm{x}) - P_{\Delta} (A_l)^* A_l P_{\Delta} \bm{v}^{(l-1)} - P_{\Delta} W \bm{\rho}^{(l-1)} = \left ( P_{\Delta} - P_{\Delta} (A_l)^* A_l P_{\Delta}  \right ) \bm{v}^{(l-1)}.
}
Hence
\ea{
\label{vl_bound}
\| \bm{v}^{(l)} \|_{2} \leq \sqrt{s^*} \nm{W^{-1} \bm{v}^{(l)} }_{\infty} \leq \sqrt{s^*} a_{l} \nm{W^{-1} \bm{v}^{(l-1)} }_{\infty} \leq \sqrt{s^*} \left( \prod^{l}_{j=1} a_j \right). 
}
This gives
\bes{
\nm{W(P_{\Delta} \bm{\rho} - \sgn(P_{\Delta} \bm{x} ) )}_{2} = \nm{\bm{v}^{(L)}}_{2} \leq \frac{\sqrt{s^*}}{2^L \log_2(\sqrt{s^*})} \leq \frac{1}{4},
}
and therefore (iii) holds.

Next consider (iv).  Using event $B_l$ and \R{vl_bound} we have
\eas{
\nm{P^{\perp}_{\Delta} \bm{\rho}^{(l)} }_{\infty} &\leq \nm{P^{\perp}_{\Delta} W^{-1} (A_l)^* A_l P_{\Delta} \bm{v}^{(l-1)} }_{\infty} + \nm{P^{\perp}_{\Delta} \bm{\rho}^{(l-1)} }_{\infty}
\\
& \leq b_{l} \nm{W^{-1} \bm{v}^{(l-1)} }_{\infty} +  \nm{P^{\perp}_{\Delta} \bm{\rho}^{(l-1)} }_{\infty}
\leq b_{l} \prod^{l-1}_{j=1} a_j + \nm{P^{\perp}_{\Delta} \bm{\rho}^{(l-1)} }_{\infty}.
}
Therefore
\bes{
\nm{P^{\perp}_{\Delta} \bm{\rho} }_{\infty} \leq \sum^{L}_{l=1} b_l \prod^{l-1}_{j=1} a_j \leq \frac{1}{4} \left ( 1 + \frac{1}{2} + \frac14 + \ldots \frac{1}{2^L} \right ) \leq \frac12,
}
which implies that (iv) holds.

Finally, consider condition (v).  Define $\bm{\xi}^{(0)} = 0$ and $\bm{\xi}^{(l)} = \sqrt{\frac{m}{m_{l}}} A_{l} \bm{v}^{(l-1)} + \bm{\xi}^{(l-1)}$,
so that $\bm{\rho}^{(l)} = A^* \bm{\xi}^{(l)}$.  Let $\bm{\xi} = \bm{\xi}^{(L)}$, which gives $\bm{\rho} = A^* \bm{\xi}$.  Then 
\be{
\label{wk_bound}
\| \bm{\xi}^{(l)} \|_2 \leq \sqrt{\frac{m}{m_{l}}} \| A_{l} \bm{v}^{(l-1)} \|_2 + \| \bm{\xi}^{(l-1)} \|_2.
}
Now
$
\| A_{l} \bm{v}^{(l-1)} \|^2_2 \leq \| \bm{v}^{(l-1)} \|^2_2 + \| \bm{v}^{(l-1)} \|_2 \| \bm{v}^{(l)} \|_2,
$
and therefore \R{vl_bound} gives
\bes{
\| A_{l}  \bm{v}^{(l-1)} \|^2_2 \leq s^* \left ( a_{l} + 1 \right ) \prod^{l-1}_{j=1} a^2_{j}.
}
We therefore deduce that
\bes{
\| \bm{\xi} \|_2 \leq  \sqrt{s} \sqrt{m} \sum^{L}_{l=1} \sqrt{\frac{a_{l} + 1}{m_{l}} } \prod^{l-1}_{j=1} a_{j},
}
Now
\bes{
\frac{a_{l} + 1}{m_{l}} = \frac{4}{m^*} \left ( 1 + \frac{1}{2 \sqrt{\log_2(\sqrt{s^*})}} \right ) \leq \frac{6}{m^*},\qquad l=1,2,
}
and
\bes{
\frac{a_{l} + 1}{m_{l}} \leq \frac{3(L-2)}{m^*},\qquad k=3,\ldots,L.
}
Hence we get
\eas{
\| \bm{\xi} \|_2  &\leq \sqrt{s^*} \sqrt{m} \left ( \sqrt{\frac{6}{m^*}} (1+1/2) + \sqrt{\frac{3(L-2)}{m^*}} \sum^{L}_{l=3} \frac{1}{2^{l-1} \log_2(\sqrt{s^*})} \right ) 
\\
& \leq \frac{4 \sqrt{s^*}\sqrt{m}}{\sqrt{m^*}}  \left ( 1 + \frac{\sqrt{L-2}}{\log_2(\sqrt{s^*})} \right ) 
\\
&\leq 8 \sqrt{s^*} \sqrt{m/m^*}.
}
We now recall that $m^* = \lfloor m-  L \rfloor \geq m - L -1$.  Since $m \geq 2 L + 2$ by assumption, we have $\sqrt{m/m^*} \leq \sqrt{2}$.  Hence condition (v) holds with $\lambda \leq 8 \sqrt{2}$.

\pbk
\textbf{Step 2.}\ We show that event $E$ holds with high probability.  By the union bound
\bes{
\bbP(E^c) \leq \sum^{L}_{l=1} \left ( \bbP(A^c_l) + \bbP(B^c_l) \right ) + \bbP(C^c) + \bbP(D^c).
}
Hence it suffices to show that
\eas{
\bbP(A^c_l) , \bbP(B^c_l) &\leq \epsilon / 16,\qquad l=1,2,
\\
\bbP(A^c_l), \bbP(B^c_l) &\leq \epsilon / (8(L-2)),\qquad l=3,\ldots,L,
\\
\bbP(C^c) , \bbP(D^c) & \leq \epsilon / 4.
}

For the events $A_l$ we apply Lemma \ref{l:tech2} to the matrices $A_l$ with the appropriate values for $\epsilon$ and $\delta$ to get, after recalling the definition of the $m_l$, the condition
\be{
\label{A_m}
m \gtrsim \Gamma(F,\bm{w},\Delta) \cdot \log(|\Delta |_{\bm{w}} ) \cdot \log(|\Delta|_{\bm{w}}/\epsilon)
}
For the events $B_l$, we apply Lemma \ref{l:tech4} to deduce, after some algebra, the condition
\be{
\label{B_m}
m \gtrsim \Gamma(F,\bm{w},\Delta) \cdot \log \left (N \log(|\Delta |_{\bm{w}} )  /\epsilon \right ).
}
Next, we note that Lemma \ref{l:tech1} implies that event $C$ holds with probability at least $1-\epsilon/4$ provided
\be{
\label{C_m}
m \gtrsim \Upsilon(F,\Delta) \cdot \log( |\Delta| / \epsilon),
}
and Lemma \ref{l:tech3} implies that event $D$ holds with probability at least $1-\epsilon/4$ provided
\be{
\label{D_m}
m \gtrsim \left(\Upsilon(F,\Delta) + \Gamma(F,\bm{w},\Delta) \right) \cdot \log(N/ \epsilon).
}
To complete the proof we note that \R{A_m}--\R{D_m} are all implied by the condition
\bes{
m \gtrsim \Upsilon(F,\Delta) \cdot \log(N/\epsilon) + \Gamma(F,\bm{w},\Delta) \cdot \left ( \log(N/\epsilon) + \log(|\Delta|_{\bm{w}}) \cdot  \log(|\Delta|_{\bm{w}}/\epsilon) \right ).
}
This gives the result.
}

\subsection{Proofs of Theorem \ref{t:main_res_derivs_identical} and Corollary \ref{c:main_cor}}
Theorem \ref{t:main_res_derivs_identical} will now follow as a corollary of the abstract recovery guarantee, Theorem \ref{t:par_acq_wl1min}, after estimating the local coherences $\Upsilon(F,\Delta)$ and $\Gamma(F,\bm{w},\Delta)$ for the derivative sampling problem.  This is done in the following two lemmas.  Note that in this section, we revert back to indexing over the multi-index set $\Lambda \subset \bbN^d_0$ (as was introduced in \S \ref{s:prelims}), rather than over the integers $\{1,\ldots,N\}$.

\lem{
\label{l:Upsilon_poly_case}
Let $\{ \phi_{\bm{n}} \}_{\bm{n} \in \bbN^d_0}$ be the orthonormal basis of tensor-product Sturm--Louiville eigenfunctions defined in \S \ref{s:SL}, $F$ be the distribution of matrices defined in \S \ref{ss:deriv_PA} for the derivative sampling problem, and suppose that $\Upsilon(F,\Delta)$ is as in Definition \ref{d:coh_rel_1}, where $\Delta \subset \bbN^d_0$ is a multi-index set.  Then
\bes{
\Upsilon(F,\Delta) \leq \max_{\bm{n} \in \Delta} \left \{ \frac{1+\kappa_{\bm{n}}}{1+\lambda_{\bm{n}}} \right \} | \Delta |_{\bm{u}},
}
where $\lambda_{\bm{n}}$, $\kappa_{\bm{n}}$ and $\bm{u}$ are as in \R{lambda_def}, \R{kappa_def} and \R{u_def} respectively
}
\prf{
Let $\bm{z} \in \bbC^N$ with $\nm{\bm{z}}_{2} = 1$ and let $B$ be as in \R{B_deriv}.  Then
\bes{
\nm{B^* P_{\Delta} \bm{z}}^{2}_{2} = \sum^{d}_{k=0} \left |  \sum_{\bm{n} \in \Delta}  \frac{\sqrt{\tau_k(\bm{y})} \ \overline{\partial_k \phi_{\bm{n}}(\bm{y})} {z_{\bm{n}}}  }{\sqrt{1+\lambda_{\bm{n}}}} \right |^2\leq \sum^{d}_{k=0} \sum_{\bm{n} \in \Delta} \frac{ \tau_k(\bm{y}) \left | \partial_k \phi_{\bm{n}}(\bm{y}) \right |^2 }{1+\lambda_{\bm{n}}}.
}
Observe that, when $k \neq 0$, 
\bes{
\tau_{k}(\bm{y})  \left | \partial_k \phi_{\bm{n}}(\bm{y}) \right |^2 
= \frac{\chi(y_k)\prod^{d}_{j=1, j \neq k}\nu(y_j)}{\mu(\bm {y})} | \phi'_{n_k}(y_k) |^2 \prod^{d}_{\substack{j=1 \\ j \neq k}} | \phi_{n_j}(y_j) |^2 
\leq \frac{\chi(y_k)}{\mu(y_k)} | \phi'_{n_k}(y_k) |^2 \prod^{d}_{\substack{j=1 \\ j \neq k}} u^2_{n_j} \leq \kappa_{n_k} u^2_{\bm{n}},
}
and therefore
\be{
\label{SL_deriv_bound}
\sum^{d}_{k=0} \tau_{k}(\bm{y})  \left | \partial_k \phi_{\bm{n}}(\bm{y}) \right |^2 \leq u^2_{\bm{n}} \left ( 1 + \sum^{d}_{k=1} \kappa_{n_k} \right ) =  u^2_{\bm{n}} (1+\kappa_{\bm{n}}).
}
Hence
\bes{
\nm{B^* P_{\Delta} \bm{z}}^{2}_{2} \leq \sum_{\bm{n} \in \Delta} \frac{1 + \kappa_{\bm{n}} }{1+\lambda_{\bm{n}} } u^2_{\bm{n}} \leq \max_{\bm{n} \in \Delta} \left \{ \frac{1+\kappa_{\bm{n}}}{1+\lambda_{\bm{n}}} \right \} | \Delta |_{\bm{u}}.
}
Since $\bm{z}$ was arbitrary we deduce the result.
}

\lem{
Let $\{ \phi_{\bm{n}} \}_{\bm{n} \in \bbN^d_0}$ and $F$ be as in Lemma \ref{l:Upsilon_poly_case} and $\Gamma(F,\bm{w},\Delta)$ be as in Definition \ref{d:coh_rel_2}.  Then
\bes{
\Gamma(F,\bm{w},\Delta) \leq |\Delta|_{\bm{u}} + \max_{\bm{n} \in \Lambda} \left \{ \frac{u^2_{\bm{n}} (1+\kappa_{\bm{n}})}{w^2_{\bm{n}} (1+\lambda_{\bm{n}})} \right \} | \Delta |_{\bm{w}}.
}
}
\prf{
Let $\bm{z} \in \bbC^N$ with $\nm{\bm{z}}_{\infty} = 1$ and $\bm{n'} \in \Lambda$.  Then
\be{
\label{green_mug}
\left | \ip{\bm{e}_{\bm{n'}}}{W^{-1} B B^* P_{\Delta} W \bm{z}} \right |= \frac{1}{w_{\bm{n'}}} \left | \sum^{d}_{k=0} \frac{\sqrt{\tau_k(\bm{y})}\partial_{k} \phi_{\bm{n'}}(\bm{y})}{\sqrt{1+\lambda_{\bm{n'}}} } \sum_{\bm{n} \in \Delta}  \frac{\sqrt{\tau_k(\bm{y})} \ \overline{\partial_{k} \phi_{\bm{n}}(\bm{y})}}{\sqrt{1+\lambda_{\bm{n}}} } w_{\bm{n}} z_{\bm{n}} \right |.
}
Hence
\eas{
\left | \ip{\bm{e}_{\bm{n'}}}{W^{-1} B B^* P_{\Delta} W \bm{z}} \right | &
\leq \frac{1}{w_{\bm{n'}}} \sum^{d}_{k=0} \sqrt{\frac{\tau_k(\bm{y}) | \partial_{k} \phi_{\bm{n'}}(\bm{y}) |^2}{1+\lambda_{\bm{n'}}} } \sum_{\bm{n} \in \Delta}  \sqrt{\frac{\tau_k(\bm{y}) |\partial_{k} \phi_{\bm{n}}(\bm{y}) |^2 }{1+\lambda_{\bm{n}} } } w_{\bm{n}}
\\
& \leq \frac{1}{w_{\bm{n'}}} \sqrt{\sum^{d}_{k=0} \frac{\tau_k(\bm{y}) | \partial_{k} \phi_{\bm{n'}}(\bm{y}) |^2}{1+\lambda_{\bm{n'}}} } \sqrt{\sum^{d}_{k=0} \left (\sum_{\bm{n} \in \Delta}  \sqrt{\frac{\tau_k(\bm{y}) |\partial_{k} \phi_{\bm{n}}(\bm{y}) |^2 }{1+\lambda_{\bm{n}} } } w_{\bm{n}}\right )^2 }
\\
& \leq \frac{1}{w_{\bm{n'}}} \sqrt{\sum^{d}_{k=0} \frac{\tau_k(\bm{y}) | \partial_{k} \phi_{\bm{n'}}(\bm{y}) |^2}{1+\lambda_{\bm{n'}}} } \sqrt{\sum_{\bm{n} \in \Delta} \frac{\sum^{d}_{k=0} \tau_k(\bm{y}) |\partial_{k} \phi_{\bm{n}}(\bm{y}) |^2  }{1+\lambda_{\bm{n}}} } \sqrt{|\Delta|_{\bm{w}} }.
}
We now apply \R{SL_deriv_bound} to get
\eas{
\left | \ip{\bm{e}_{\bm{n'}}}{W^{-1} B B^* P_{\Delta} W \bm{z}} \right | & \leq   \frac{u_{\bm{n'}}}{w_{\bm{n'}}} \sqrt{\frac{1+\kappa_{\bm{n'}} }{1+\lambda_{\bm{n'}}} } \sqrt{\sum_{\bm{n} \in \Delta} \frac{1+\kappa_{\bm{n}} }{1+\lambda_{\bm{n}}} u^2_{\bm{n}} } \sqrt{|\Delta|_{\bm{w}}}
\\
& \leq \sqrt{\max_{\bm{n} \in \Lambda} \left \{ \frac{u^2_{\bm{n}} (1+\kappa_{\bm{n}})}{w^2_{\bm{n}} (1+\lambda_{\bm{n}})} \right \}} \sqrt{|\Delta|_{\bm{u}} |\Delta|_{\bm{w}} }.
}
Since $\bm{z}$ and $\bm{n'}$ were arbitrary, after an application of the inequality $ab \leq a^2/2+b^2/2$, we obtain
\be{
\label{Gamma1_deriv}
\Gamma_1(F,\bm{w},\Delta) \leq  \frac12 | \Delta |_{\bm{u}} + \frac12 \max_{\bm{n} \in \Lambda} \left \{ \frac{u^2_{\bm{n}} (1+\kappa_{\bm{n}})}{w^2_{\bm{n}} (1+\lambda_{\bm{n}})} \right \} | \Delta |_{\bm{w}}.
} 
We now consider $\Gamma_2(F,\bm{w},\Delta)$.  From \R{green_mug} and \R{SL_deriv_bound} we have
\eas{
\bbE \left | \ip{\bm{e}_{\bm{n'}}}{W^{-1} B B^* P_{\Delta} W \bm{z}} \right |^2 & \leq \frac{1}{w^2_{\bm{n'}}} 
 \bbE  \left ( \sum^{d}_{k=0} \frac{\tau_k(\bm{y}) | \partial_k \phi_{\bm{n'}}(\bm{y}) |^2}{1+\lambda_{\bm{n'}}} 
 \left( \sum^{d}_{k=0} \left | \sum_{\bm{n} \in \Delta}  \frac{\sqrt{\tau_k(\bm{y})}\partial_{k} \phi_{\bm{n}}(\bm{y})}{\sqrt{1+\lambda_{\bm{n}}} } w_{\bm{n}} z_{\bm{n}}    \right |^2 \right ) \right) 
\\
& \leq \frac{u^2_{\bm{n'}} (1+\kappa_{\bm{n'}})}{w^2_{\bm{n'}} (1+\lambda_{\bm{n'}})} \bbE  \sum^{d}_{k=0} \tau_k(\bm{y}) \left | \sum_{\bm{n} \in \Delta}  \frac{\partial_{k} \phi_{\bm{n}}(\bm{y})}{\sqrt{1+\lambda_{\bm{n}}} } w_{\bm{n}} z_{\bm{n}}    \right |^2
\\
& = \frac{u^2_{\bm{n'}} (1+\kappa_{\bm{n'}})}{w^2_{\bm{n'}} (1+\lambda_{\bm{n'}})} \sum^{d}_{k=0} \int_{D} \left | \sum_{\bm{n} \in \Delta}  \frac{\partial_{k} \phi_{\bm{n}}(\bm{y})}{\sqrt{1+\lambda_{\bm{n}}} } w_{\bm{n}} z_{\bm{n}}    \right |^2 \nu_k(\bm{y}) \D \bm{y}.
}
Recall that the functions $\partial_{k} \phi_{\bm{n}}$ are orthogonal with respect to the weight function $\nu_k$, and that $\int_{D} | \partial_{k} \phi_{\bm{n}}(\bm{y}) |^2 \nu_k(\bm{y}) \D \bm{y} = \lambda_{n_k}$.  Therefore, by Parseval's identity, we get 
\eas{
\bbE \left | \ip{\bm{e}_{\bm{n'}}}{W^{-1} B B^* P_{\Delta} W \bm{z}} \right |^2 & \leq \frac{u^2_{\bm{n'}} (1+\kappa_{\bm{n'}})}{w^2_{\bm{n'}} (1+\lambda_{\bm{n'}})} \sum^{d}_{k=0} \sum_{\bm{n} \in \Delta} \frac{|w_{\bm{n}} z_{\bm{n}} |^2 \lambda_{n_k}}{1+\lambda_{\bm{n}}}
\\
& = \frac{u^2_{\bm{n'}} (1+\kappa_{\bm{n'}})}{w^2_{\bm{n'}} (1+\lambda_{\bm{n'}})} \sum_{\bm{n} \in \Delta} |w_{\bm{n}} z_{\bm{n}} |^2
 \leq \max_{\bm{n} \in \Lambda} \left \{ \frac{u^2_{\bm{n}} (1+\kappa_{\bm{n}})}{w^2_{\bm{n}} (1+\lambda_{\bm{n}})} \right \} | \Delta |_{\bm{w}},
}
where in the last step we recall that $\nm{\bm{z}}_{\infty} = 1$.  Since $\bm{z}$ and $\bm{n'}$ were arbitrary, we deduce that 
\bes{
\Gamma_{2}(F,\bm{w},\Delta) \leq \max_{\bm{n} \in \Lambda} \left \{ \frac{u^2_{\bm{n}} (1+\kappa_{\bm{n}})}{w^2_{\bm{n}} (1+\lambda_{\bm{n}})} \right \} | \Delta |_{\bm{w}}.
}
Combining this with \R{Gamma1_deriv} now completes the proof.
}

\prf{[Proof of Theorem \ref{t:main_res_derivs_identical}]
With the previous two lemmas in hand, we now apply Theorem \ref{t:par_acq_wl1min}.  Note that this gives the error estimate
$
\nm{\bm{z}_{\Lambda} - \bm{\hat{z}} }_2 \lesssim \nm{\bm{z}_{\Lambda} - \bm{z}_{\Delta} }_{1,\bm{w}} + \sqrt{|\Delta|_{\bm{w}}} \eta.
$
We now recall that $\bm{z}_{\Lambda} = Q \bm{x}_{\Lambda}$, $\bm{\hat{z}} = Q\bm{\hat{x}}$ and $\hat{f} = \sum_{\bm{n} \in \Lambda} x_{\bm{n}} \phi_{\bm{n}}$.  Hence
\bes{
\nmu{f_{\Lambda} - \hat{f}}_{\tilde{H}^1(D)} = \nm{Q(\bm{x}_{\Lambda} - \bm{\hat{x}})}_2 = \nm{\bm{z}_{\Lambda} - \bm{\hat{z}} }_2 \lesssim \nm{\bm{x}_{\Lambda} - \bm{x}_{\Delta} }_{1,\bm{v}} + \sqrt{|\Delta|_{\bm{w}}} \eta,
}
where $v_{\bm{n}} = \sqrt{1+\lambda_{\bm{n}}} w_{\bm{n}}$, $\bm{n} \in \bbN^d_0$.  The result now follows from the triangle inequality.
}

We may now also prove Corollary \ref{c:main_cor}:

\prf{[Proof of Corollary \ref{c:main_cor}]
Given $s \geq 1$, let $\Delta$ be a lower set with $| \Delta | \leq s$ such that
$
\nm{\bm{x}_{\Lambda} - \bm{x}_{\Delta}}_{1,\bm{v}} = \sigma_{s,L}(\bm{x}_{\Lambda})_{1,\bm{v}}.
$
Note that $\Delta \subseteq \Lambda$ since $\Lambda = \Lambda^{\mathrm{HC}}_{s}$ is the union of all lower sets of size at most $s$.  Furthermore, it is known that
$
N = | \Lambda^{\mathrm{HC}}_s | \leq \min \left \{ 2 s^3 4^d , \E^2 s^{2 + \log_2(d)} \right \}.
$
See, for example, \cite[Eqn.\ (10)]{BASBCWMatheon}.  In particular,
$
\log(N/\epsilon) \lesssim \min \left \{ \log(s/\epsilon) + d  ,  \log(s/\epsilon)  \log(2d)\right \}.
$
We now apply Theorem \ref{t:main_res_derivs_identical} with $\bm{w} = \bm{u}$, noting that $|\Delta|_{\bm{u}} \leq K(s)$.
}

\subsection{Proofs of Corollaries \ref{c:main_res_jacobi_identical_lower} and \ref{c:main_res_leg_lower}}\label{ss:corproof}

We first require some further background on Jacobi polynomials.  For $\alpha,\beta > -1$ and $n \in \bbN_0$, let $P^{(\alpha,\beta)}_{n}$ be the Jacobi polynomial of degree $n$.  These polynomials are orthogonal on $(-1,1)$ with respect to the weight function $\omega^{(\alpha,\beta)}(y) = (1-y)^{\alpha}(1+y)^{\beta}$, and satisfy
\bes{
\ip{P^{(\alpha,\beta)}_{n} }{P^{(\alpha,\beta)}_{m}}_{L^2_{\omega^{(\alpha,\beta)}}} = \kappa^{(\alpha,\beta)}_{n} \delta_{n,m},
}
where
\bes{
\kappa^{(\alpha,\beta)}_{n} = \frac{2^{\alpha+\beta+1}}{2 n + \alpha + \beta + 1} \frac{\Gamma(n+\alpha+1) \Gamma(n+\beta+1)}{\Gamma(n+1) \Gamma(n+\alpha+\beta+1)}.
}
These polynomials are normalized so that
$
P^{(\alpha,\beta)}_{n}(1) = \left ( \begin{array}{c} n+ \alpha \\ n \end{array} \right ),
$
Moreover, if $\alpha,\beta \geq -1/2$ then 
\be{
\label{Jacobi_max}
\sup_{y \in (-1,1)} | P^{(\alpha,\beta)}_{n}(y) | =  \left ( \begin{array}{c} n+ q \\ n \end{array} \right ) \sim \frac{n^{q}}{\Gamma(q+1)},\quad n \rightarrow \infty,
}
where $q = \max \{ \alpha , \beta \}$.  See, for example, \cite[Thm.\ 7.32.1]{SzegoOrthPolys}.  We also note the reflection property
\be{
\label{Jacobi_reflection}
P^{(\alpha,\beta)}_{n}(y) = (-1)^n P^{(\beta,\alpha)}_{n}(-y).
}
Let
$
c^{(\alpha,\beta)} = \int^{1}_{-1} \omega^{(\alpha,\beta)}(y) \D y ,
$
and define the probability density function
$
\nu^{(\alpha,\beta)}(y) = \frac{\omega^{(\alpha,\beta)}(y)}{c^{(\alpha,\beta)}}.
$
Then the corresponding orthonormal polynomials with respect to this density are given by
\be{
\label{Jacobi_phin}
\phi_{n}(y) =  \frac{P^{(\alpha,\beta)}_n(y)}{\sqrt{\kappa^{(\alpha,\beta)}_n c^{(\alpha,\beta)}}} ,\qquad n \in \bbN_0.
}

\prf{[Proof of Corollary \ref{c:main_res_jacobi_identical_lower}]
In view of Corollary \ref{c:main_cor}, it suffices to show that
$\kappa_{\bm{n}} \lesssim \lambda_{\bm{n}}$, $\forall \bm{n} \in \bbN^d_0$.
Since $\lambda_{\bm{n}} = \sum^{d}_{k=1} \lambda_{n_k}$ and $\kappa_{\bm{n}} = \sum^{d}_{k=1} \kappa_{n_k}$ (see \R{lambda_def} and \R{kappa_def} respectively), it is enough to show
$\kappa_{n} \lesssim \lambda_{n}$, $\forall n \in \bbN_0$.
Using the definition of $\kappa_{n}$ (see \R{kappa_def_1D}), and fact that $\chi(y) = \frac{1}{c_{\alpha,\beta}}(1-y)^{\alpha+1}(1+y)^{\beta+1}$ and $\nu(y) = (1-y)^{\alpha}(1+y)^{\beta} / c_{\alpha,\beta}$ in the Jacobi case (see \R{Jacobi_pw}), this is equivalent to
\bes{
\sup_{y \in (-1,1)} \left \{ \sqrt{1-y^2} | \phi'_n(y) | \right \} \lesssim u_n \sqrt{\lambda_n}.
}
Furthermore, using \R{lambda_Jacobi}, \R{Jacobi_phin} and \R{Jacobi_max}, we see that it is sufficient to show that
\be{
\label{Jacobi_required}
\sup_{y \in (-1,1)} \left \{ \sqrt{1-y^2} \left | \left ( P^{(\alpha,\beta)}_{n}(y) \right )' \right | \right \} \lesssim n^{1+q}.
}
Note that from this equation onwards we allow the constant implied by the expression $\lesssim$ to depend on $\alpha$ and $\beta$.  The derivatives of the Jacobi polynomials satisfy the following bound:
\be{
\label{Jacobi_deriv}
\left | \frac{\D P^{(\alpha,\beta)}_{n}(y)}{\D y} \Bigg |_{y = \cos(\theta)} \right | \lesssim \left \{ \begin{array}{cc} \theta^{-\alpha-3/2} n^{1/2} & c n^{-1} \leq \theta \leq \pi/2 \\ n^{2+\alpha} & 0 \leq \theta \leq c n^{-1} \end{array} \right . ,
}
(see \cite[Thm.\ 7.32.4]{SzegoOrthPolys}).
Using this and the fact that $\sin(\theta) \leq \theta$ for $0 \leq \theta \leq \pi/2$, we deduce that
\eas{
\sup_{0 \leq y \leq 1} \sqrt{1-y^2} \left | \left ( P^{(\alpha,\beta)}_{n}(y) \right )' \right | &= \sup_{0 \leq \theta \leq \pi/2} \sin(\theta) \left | \frac{\D P^{(\alpha,\beta)}_{n}(y)}{\D y} \Bigg |_{y = \cos(\theta)} \right | 
\\
& \lesssim \max \left \{ \sup_{c n^{-1} \leq \theta \leq \pi/2} \theta^{-\alpha-1/2} n^{1/2}  , \sup_{0 \leq \theta \leq c n^{-1}} n^{2+\alpha} \theta  \right \}
\lesssim n^{\alpha + 1}.
}
Now suppose that $-1 \leq y \leq 0$.  Using \R{Jacobi_reflection} and replacing $\alpha$ with $\beta$ in the above arguments, we deduce that 
\bes{
\sup_{-1 \leq y \leq 0} \sqrt{1-y^2} \left | \left ( P^{(\alpha,\beta)}_{n}(y) \right )' \right | \lesssim n^{\beta+1},
}
Therefore \R{Jacobi_required} follows immediately, completing the proof.
}

\prf{[Proof of Corollary \ref{c:main_res_leg_lower}]
As in the proof of Corollary \ref{c:main_res_jacobi_identical_lower}, we first need to show that $\kappa_n \lesssim \lambda_n, \ \forall n \in \mathbb{N}_0$, which is equivalent to 
\be{
\sup_{y \in (-1,1) } (1-y^2)^ {3/4}|\phi_n'(y)| \lesssim u_n \sqrt{\lambda_n}, \label{u_eqn}
}
where
\bes{
u_n = \sup_{y \in (-1,1)} \left({\pi}/{2} \right)^{1/2}(1-y^2)^{1/4}|\phi_n(y)|.
}
We first seek a lower bound for $u_n$.  The classical Legendre polynomials $P_n = P^{(0,0)}_n$ satisfy
\bes{
P_n^{(0,0)}(\cos \theta) = 2^{1/2}(\pi n \sin \theta)^{-1/2}\cos\left(\left(n+1/2\right)\theta -\pi/4 \right)+ \ordu{n^{-3/2}},  \qquad 
0 < \theta < \pi.  
}
See  \cite[Thm.\ 8.21.2]{SzegoOrthPolys}.  This formula holds uniformly in the interval $\epsilon \leq \theta \leq \pi - \epsilon$.  When $n$ is even, Legendre polynomials have extrema at $\cos(\theta) = 0$, i.e. $\theta = \frac{\pi}{2}$.
Then, we have 
\begin{align*}
P_n^{(0,0)}(0) = 2^{1/2} (\pi n)^{-1/2} \cos \left(n \pi/2 \right) + \ordu{n^{-3/2}} = \left(2/\pi\right)^{1/2}n^{-1/2}(-1)^{n/2} + \ordu{n^{-3/2}}.
\end{align*}
When $n$ is odd, we consider the point $\theta = \frac{\pi}{2} + \epsilon_n$, where $\epsilon_n = \pi/(2n+1)$.  Then
\begin{align*}
P_n^{(0,0)}\left( \cos \left( \pi/2 + \epsilon_n \right)\right ) &= 2^{1/2} (\pi n)^{-1/2} 
(\cos \left(\epsilon_n \right))^{-1/2} \cos \left(n\pi/2 + \left( n + 1/2\right)\epsilon_n\right)+ \ordu{n^{-3/2}}
\\
& = (2/\pi)^{1/2} n^{-1/2} (-1)^{(n-1)/2} + \ordu{n^{-3/2}}
\end{align*}
Therefore, for both even and odd $n$, we have
\bes{
\sup_{y \in (-1,1)} (1-y^2)^{1/4}|P_n^{(0,0)}(y)| \gtrsim n^{-1/2},
} 
and since $\phi_n(y) = \sqrt{2n+1} P_n^{(0,0)} (y)$, we deduce that $u_n \gtrsim 1$.
Since $\lambda_n = n(n+1)$ , we then see that \R{u_eqn} is now implied by
\be{
\sup_{y \in (-1,1) } (1-y^2)^ {3/4}|\phi_n'(y)| \lesssim n.
\label{phi_prime_bnd}
}
Using \R{Jacobi_deriv} with $\alpha = \beta = 0$ and arguing as in Corollary \ref{c:main_res_jacobi_identical_lower} we obtain
\eas{
\sup_{ 0\leq y \leq 1 }  (1-y^2)^ {3/4}
\left | \frac{\D P^{(0,0)}_{n}(y)}{\D y}  \right |
\lesssim 
\max \left\{  \sup_{c n^{-1} \leq \theta \leq \pi/2} n^{1/2}  , 
\sup_{0 \leq \theta \leq c n^{-1}}n^{2} \theta^{3/2} \right\}
\lesssim n^{1/2}.
}
Using the reflection property and the fact that $\phi_n(y) = \sqrt{2n+1} P_n^{(0,0)} (y)$, we now deduce \R{phi_prime_bnd}.

With this in hand, we apply Corollary \ref{c:main_cor} to get that the conclusions of Corollary \ref{c:main_res_leg_lower} hold under the condition
\bes{
m \gtrsim K(s) \cdot L,
}
where $L = \left ( \min \{ d+ \log(s/\epsilon) , \log(2d) \log(s/\epsilon) \} + \log(K(s)) \cdot \log(K(s) / \epsilon ) \right )$.  It remains to estimate $K(s)$ and $L$.  In \cite[Cor.\ 7.7]{AdcockCSFunInterp}, it was shown that
$
K(s) \lesssim \min \left\{2^ds, (\pi/2)^d s^{\log(1+4/\pi)/\log(2)} \right \}.
$
From this, we also observe that $\log(K(s)) \lesssim d + \log(s)$ and $\log(K(s)/\epsilon) \lesssim d + \log(s/\epsilon)$, and therefore
$
L \lesssim (d+\log(s)) (d+\log(s/\epsilon)),
$
which completes the proof.
}

\section{Conclusions}

In this paper we have studied the sparse polynomial approximation of a high-dimensional function from measurements of both the function and its gradient.  Our main results show that gradient-augmented measurements permit an error bound in a stronger Sobolev norm as opposed to a $L^2$-norm, for the same sample complexity.  Numerically, we observe recovery from gradient-augmented measurements gives smaller errors (when measured in a fixed norm) than the case of function samples only, under a reasonable model of computational cost.

There are several areas for future work.  First, in high dimensions the Sobolev norm is weaker than in low dimensions (see, for instance, the Sobolev embedding theorem).  This might suggest the improvement due to gradient samples lessens in higher dimensions, yet this is seemingly at odds with our numerical experiments.  In particular, Fig.\ \ref{Linfty} shows a consistent improvement even though the error is measured in the $L^{\infty}$-norm.  Second, as noted in \S \ref{ss:discussion}, our recovery guarantees are nonuniform, and correspondingly the error bounds are worse than those obtained from uniform recovery guarantees.  Deriving uniform recovery guarantees in the case of gradient-augmented measurements (for example, extending the work of \cite{Chkifalower}), is an open problem.  Finally, as mentioned in \S \ref{s:introduction}, Hermite interpolation as pursued in this paper is not the only way gradient information could be used to enhance the approximation.  A thorough comparison of this with other approaches is a topic for future work.

\section*{Acknowledgements}
This work is supported in part by the NSERC grant 611675 and an Alfred P.\ Sloan Research Fellowship.  Yi Sui also acknowledges support from an NSERC PGSD scholarship.

\bibliographystyle{abbrv}
\small
\bibliography{FuncApproxDerivativesRefs}

\end{document}